\documentclass[a4paper,12pt]{article}

\begin{document}

\title{Normal families of functions and groups of pseudoconformal
diffeomorphisms of quaternion and octonion variables.}
\author{Ludkovsky S.V.}
\date{25 October 2005}
\maketitle

\begin{abstract}

This paper is devoted to the specific class of pseudoconformal
mappings of quaternion and octonion variables. Normal families of
functions are defined and investigated. Four criteria of a family
being normal are proven. Then groups of pseudoconformal
diffeomorphisms of quaternion and octonion manifolds are
investigated. It is proven, that they are finite dimensional Lie
groups for compact manifolds. Their examples are given. Many
characteristic features are found in comparison with commutatiive
geometry over $\bf R$ or $\bf C$.

\end{abstract}

\section{Introduction.}
As it is well-known quaternions were defined and investigated by
W.R. Hamilton in the sixties of the 19-th century
\cite{hamilt,waerd}. Few years later on Cayley had made their
generalization and then Cayley-Dickson algebras ${\cal A}_r$ were
studied in details \cite{baez,kansol,kurosh,ward}. For $r=1$ they
coinside with the complex field $\bf C$, for $r=2$ with the skew
field of quaternions $\bf H$, for $r=3$ with the algebra of
octonions $\bf O$, which is already nonassociative, but alternative.
The Cayley-Dickson algebras for $r\ge 4$ have not the division
property or multiplicative norm. \par Groups of diffeomorphisms of
connected manifolds are important not only in mathematics, but also
in theoretical physics, for example, in gauge theories, theory of
superstrings, quantum field theory, quantum gravity
\cite{bogsch,bolokto,solov,subfai}. There are substantial
differences in the structure of groups of diffeomorphisms of real
and complex manifolds \cite{kobtrgr,pontr}. In the first case they
are infinite dimensional and do not satisfy the Campbell-Hausdorff
formula even for compact manifolds, but in the second case they are
finite dimensional Lie groups, consequently, they satisfy the
Campbell-Hausdorff formula locally
\cite{balafrat,boch44}-\cite{bochmon47,hcart35l,lurim1}. To groups
of diffeomorphisms of real manifolds numerous papers were devoted,
for example, \cite{balafrat,ebmars,lurim1} and references therein,
but it is impossible to refer all here. Groups of diffeomorphisms of
complex manifolds were studied in classical works and recent also
\cite{hcart34,hcart35l,boch44,bochmon45, bochmon47,myestin38,land04}
and originate from publications of H. Cartan. On the other hand,
they certainly are based on classical works on complex analysis of
Cauchy, Hartogs, Levi, Montel, Picard, Riemann, Rouch\'e
\cite{julia,lavrsch,schab}. Therefore, present studies of groups of
pseudoconformal diffeomorphisms of quaternion and octonion variables
have demanded investigations of the corresponding aspects of
function theory over quaternions and octonions in this paper in
addition to that of previous.
\par Recently the theory of superdifferentiable functions of quaternion
and octonion variables is being developed
\cite{luladfcdv,ludoyst,ludfov}. Shortly we shall write
differentiable functions instead of superdifferentiable functions,
because here are concrete algebras quaternion and octonion of low
dimensions over $\bf R$. Quaternion and octonion manifolds and their
transformation groups can be used in quantum field theory and in
problems with spin, as well as in the aforementioned theories, since
there are embeddings of the group $SU(2)$ in the skew field $\bf H$
of quaternions. Yang and Mills have used in their theory functions
of quaternion variables and their followers have used functions of
octonion variables, but they felt drawback of insufficient
development of their mathematical theory in their times
\cite{guetze}. Some aspects of the noncommutative geometry over
quaternions and octonions including physical applications can be
found in \cite{emch,lawmich,guetze,oystaey,rothe}. But in previous
works the noncommutative geometry was considered mainly over
Grassman algebras or in the framework of $C^*$-algebra theory over
complex numbers \cite{connes,dewitt,khren}.
\par  This paper is devoted to the investigation of structures
of groups of diffeomorphisms of quaternion and octoniuon manifolds,
which were not studied earlier. In this work the quaternion and
octonion analogs of pseudoconformal complex mappings are introduced.
It is proved, that the families of pseudoconformal mappings are
sufficiently large, though they do not form vector spaces. It is
proved further, that groups of pseudoconformal diffeomorphisms of
compact quaternion and octonion manifolds are finite dimensional
analytic Lie groups over the real field. Then examples of
evaluations of such groups for certain domains are given.
\par The second section is devoted to the investigation of families
of functions, which bear the property of being pseudoconformal
besides some points of a common domain. They are called normal
families of functions. Four criteria of a family being normal are
proved in Theorems 13, 20, 27 and 49. Their properties are described
and examples are given. Among main results of the second section are
also noncommutative versions of the Montel, Riemann, Schwarz
theorems. In theorem 21 it is proven, that each zero of a
pseudoconformal function on an open domain in $\bf K$ is isolated.
The argument principle for pseudoconformal mappings is given in
Theorem 23 and there is proved the noncommutative analog of Rouch\'e
theorem 24. \par Theorem 30 states that an image of an open domain
in $\bf H$ or $\bf O$ under a pseudoconformal mapping is an open
domain. In particular, fractional $\bf R$-linear pseudoconformal
mappings of quaternion and octonion variables are studied. In
Theorem 33 it is proved their property to transform a hypersphere
into a hypersphere contained in the one point compactification of
$\bf H$ or $\bf O$. Symmetry properties of pseudoconformal functions
are presented in Theorem 34. \par In Theorems 39 and 42 the symmetry
principle and extension theorems for pseudoconformal mappings of
noncommutative geometry over quaternions and octonions are proved.
In the complex case they were proven by Riemann and Schwarz. The
noncommutative version of the monodromy is studied in Theorem 45.
The nonassociative noncommutative analog over octonions of the
Riemann mapping theorem in terms of pseudoconformal mappings is
investigated in Theorem 47. Its quaternion version is contained in
\cite{lusmldg05}. The quaternion and octonion analogs of modular
functions are constructed in Example 48. Properties of families of
functions being normal in open domains besides some points or
surfaces are described at the end of the second section, in
particular, in Theorems 54 and 63. Finally the example is given.
\par It is necessary to note, that quaternion and octonion holomorphic
functions take an intermediate place between complex holomorphic and
real locally analytic functions, having their own specific features.
More narrow class compose pseudoconformal mappings of quaternion and
octonion variables, because they are holomorphic functions subjected
to additional conditions, but in some respect they are not so rigid
as complex holomorphic or pseudoconformal functions of several
variables. Properties of quaternion and octonion pseudoconformal
mappings have many specific features in comparison with the complex
case due to the noncommutativity of $\bf H$ and $\bf O$ and
nonassociativity of $\bf O$.
\par The third section is devoted to groups of pseudoconformal
diffeomorphisms.  It is proven in Theorems 23, 24 and 25 that groups
of pseudoconformal diffeomorphisms of compact manifolds over
quaternions or octonions are finite dimensional Lie groups and the
estimation of their dimension is done. Their examples are given in
Theorems 37, 38 and 40. This has demanded to prove Theorems 30 and
34 related with derivatives and the maximum principle of
pseudoconformal functions. In Propositions 28 and 29 the property
for a domain being perfect and an interplay with its group of
pseudoconformal diffeomorphisms is studied. There is a tight
relation between groups of pseudoconformal diffeomorphisms of
compact domains and locally compact domains having one point
(Alexandroff) compactification, since it is possible to consider the
condition of a function being pseudoconformal at the infinity.
Though for certain unbounded quaternion or octonion domains groups
of pseudoconformal diffeomorphisms need not be locally compact. In
addition groups of holomorphic diffeomorphisms of holomorphic
noncommutative manifolds are investigated in Theorems 2, 13, 15, 18
and 21. For this the noncommutative analog of the Kovalevsky theorem
of the Cauchy problem in partial differential equations in the class
of holomorphic functions is proved in Theorem 4. With the help of it
noncommutative holomorphic manifolds are studied in Proposition 11
and Theorem 12. In particular, it is proved in Theorem 10 that
compact pseudoconformal manifolds are orientable.
\par The results of
sections 2 and 3 show that there are substantial differences between
the commutative case over $\bf C$ and the noncommutative case over
quaternions and octonions. All results of this paper are proved for
the first time. It is possible to develop this theme further
considering groups of pseudoconformal diffeomorphisms with specific
behaviour near a boundary of a domain.
\section{Normal families of functions.}
\par {\bf 1. Definition.} Let $U$ be an open domain in $\bf K$
the skew field of quaternions ${\bf K}=\bf H$ or in the algebra of
octonions ${\bf K}=\bf O$, then the (super)differentiable function
$f: U\to \bf K$ of quaternion or octonion variable $z\in U$ we call
pseudoconformal at a point $\xi \in U$, if it satisfies Conditions
$(1-3)$: \par $(1)$ $\partial f(z)/\partial {\tilde z}=0$ for $z=\xi
$; \par $(2)$ $Re \{ [(\partial f(z)/
\partial z).h_1] [(\partial f(z))/ \partial z).h_2]^{*} \} |h_1|
|h_2| = |(\partial f(z)/\partial z).h_1| |(\partial f(z)/\partial
z).h_2| Re (h_1{\tilde h}_2)$ for $z=\xi $ for each $h_1, h_2\in \bf
K$,
\par $(3)$ $(\partial f(z)/\partial z)|_{z=\xi }\ne 0$, where ${\tilde
z}=z^*$ denotes the adjoint quaternion or octonion $z$ such that $z
{\tilde z} = |z|^2$, $Re (z) := (z+{\tilde z})/2$; for $f$ it is
used either the shortest phrase compatible with these conditions or
in the underlying real space (shadow)  $\bf R^4$ or $\bf R^8$
nonproper rotations associated with $f'$ are excluded. For short
$f(z,{\tilde z})$ is written as $f(z)$ due to the bijectivity
between $z\in \bf K$ and $\tilde z$. If $f$ is pseudoconformal at
each point $\xi \in U$, then it is called pseudoconformal in the
domain $U$.
\par {\bf 2. Lemma.} {\it If a differentiable function $f(z)$
of a quaternion or octonion variable  $z$ on an open domain $U$ in
$\bf K$ satisfies Conditions $(1-7)$: $(1)$ $C^2\ni f({\tilde
z})=(f(z))^*$ in a neighborhood (open) $V$ of a point $\xi \in U$;
$(2)$ $f(z)f({\tilde z})=\phi (a(z,{\tilde z}))$, $C^2\ni \phi :
W\to \bf R$, where $W$ is an open domain in $\bf R$; $(3)$
$a=a(z,{\tilde z})\in \bf R$ for each $z\in V$, $W\supset \{ x:
x=a(z,{\tilde z}): z\in V \} $; $(4)$ $(\partial a(z,{\tilde
z})/\partial z).h= [(\partial a(z,{\tilde z})/\partial {\tilde
z}).h]^*$ for $z=\xi $ for each $h\in \bf K$; $(5)$ $(\partial
^2a(z,{\tilde z})/\partial z\partial {\tilde z}).(h_1,h_2)=
(\partial ^2a(z,{\tilde z})/\partial {\tilde z}\partial
z).(h_1,h_2)=\psi (z,{\tilde z}) {\tilde h}_1 h_2$ in the
neighborhood $V$ for each $h_1, h_2\in \bf K$, where $\psi $ is some
real valued function; $(6)$ $[(d^2\phi (a)/da^2)|\partial a/\partial
z|^2+(d\phi (a)/da)\partial ^2a(z,{\tilde z})/\partial z\partial
{\tilde z}]\ne 0$ for $z=\xi $; $(7)$ $\partial f(z)/\partial
{\tilde z}=0$ for $z=\xi $; then $f$ is pseudoconformal at the point
$\xi $.}
\par {\bf Proof.} Since $f$ is differentiable by $z$, then due to
Condition $(7)$ it is holomorphic, hence it is locally analytic. In
view of Conditions $(1,7)$ and the alternativity of the algebra $\bf
K$ there is the equality:
\par $(i)$ $(\partial ^2(f(z)f(\tilde z))/\partial z \partial
{\tilde z}).( h_1,h_2)=[(\partial f(z)/
\partial z).h_1] [(\partial f(z))/
\partial z).h_2]^{*},$ in particular, \par
$(ii)$ $(\partial ^2(f(z)f(\tilde z))/\partial z\partial {\tilde
z}).( h,h)=|(\partial f(z)/ \partial z).h|^2> 0$ for $h=h_1=h_2\ne
0$, hence Conditions $(1,3)$ of Definition 1 are satisfied. Since
$\phi (a)\in \bf R$ for each $a\in W$, then $d\phi (a)/da$ and
$d^2\phi (a)/da^2\in \bf R$ on $W$. Due to Conditions $(2-5)$ the
second mixed derivative of the composite function $\phi (a)$ has the
form: $(\partial ^2\phi (a(z,{\tilde z}))/\partial z\partial {\tilde
z}).(h_1,h_2)=(d^2\phi (a)/da^2)|\partial a/\partial z|^2{\tilde
h}_1h_2$ \\ $+(d \phi (a)/da)(\partial ^2a(z,{\tilde z})/\partial
z\partial {\tilde z}).(h_1,h_2)$, consequently, \\ $Re [(\partial
^2\phi (a(z,{\tilde z}))/\partial z\partial {\tilde
z}).(h_1,h_2)]=[(d^2\phi (a)/da^2)|\partial a/\partial z|^2+(d\phi
(a)/da)\psi (z,{\tilde z})] Re ({\tilde h}_1h_2)$, since $Re
(a{\tilde b})=Re ({\tilde a}b)=Re ({\tilde b}a)$ for each $a, b\in
\bf K$. Then from $(i,ii)$ and Condition $(6)$ it follows Equality
$(2)$ of Definition 1.
\par {\bf 3. Corollary.} {\it The functions $f_1(z):=z+b$,
$f_2(z):=(az)q$ or $f_2(z)=a(zq)$ on $\bf K$ for constants $a, b,
q\in \bf K$ with $aq\ne 0$; $(1-z)^{\alpha }$ with $\alpha \in {\bf
R}\setminus \bf Z$ on ${\bf K}\setminus \{ 1 \} $, $e^z$ on $\bf K$
and $z^n$ on ${\bf K}\setminus \{ 0 \} $ are pseudoconformal, where
$n=1,2,3,...$.}
\par {\bf Proof.} If $\xi , \eta \in \bf K$ with $\xi - {\tilde \xi }
=\beta (\eta - {\tilde \eta })$ for some $\beta \in \bf R$, then
$\exp (\xi )\exp (\eta )=\exp (\xi +\eta )$ due to Proposition 3.2
\cite{ludfov}. For $f(z)=(1-z)^{\alpha }$ put $\phi (a(z,{\tilde
z}))= f(z)f({\tilde z})$, where $f(z)f({\tilde z})=(1+a)^{\alpha }$,
$\alpha \in {\bf R}\setminus \bf Z$. Then the function $a$ has the
form: $a(z,{\tilde z})= -z -{\tilde z} + z {\tilde z}$, since
$(1-z)^{\alpha }=\exp (\alpha Ln (1-z))$. Therefore, $d\phi (a)/da
=\alpha (1+a)^{\alpha -1}$, $d^2\phi (a)/da^2 = \alpha (\alpha -1)
(1+a)^{\alpha -2}$, hence $Re (\partial ^2\phi (a(z,{\tilde
z})/\partial z\partial {\tilde z}).(h_1,h_2) = (d\phi (a)/da
+d^2\phi (a)/da^2 |1-z|^2)Re (h_1{\tilde h}_2)$.
\par For $f(z)=e^z$ the function $a$ is $a(z,{\tilde z})=
z+{\tilde z}$, $\phi (a) = e^a$, $(\partial a/\partial z).h=h$,
$(\partial a/\partial {\tilde z}).h={\tilde h}$, $\partial ^2
a/\partial z\partial {\tilde z}=0$. Thus $Re [(\partial ^2\phi
(a(z,{\tilde z}))/\partial z\partial {\tilde z}).(h_1,h_2)]=e^a Re
({\tilde h}_1h_2)$, hence by Lemma 2 $e^z$ is pseudoconformal.
\par {\bf 4. Theorem.} {\it A function $f$ is pseudoconformal at
a point $\xi $ or on an open domain $U$ in $\bf H$ if and only if,
there exist nonzero quaternions $a, b\in \bf H$ or quaternion
holomorphic functions $a, b: U\to \bf H$ such that $f'(\xi ).h =
ahb$ for each $h\in \bf H$ or $f'(z).h = a(z)hb(z)$ for each $z\in
U$ and $h\in \bf H$ respectively.}
\par {\bf Proof.} Mention that the skew field of quaternions
$\bf H$ has the real shadow which is the linear space $\bf R^4$ over
$\bf R$, the scalar product in which can be written in the form
$(a,b):=Re ({\tilde a}b)$ for each $a, b\in \bf H$. To each
quaternion $b\ne 0$ it can be provided two functions $b_L$ and $b_R$
by the formulas: $b_L(z) := bz$ and $b_R(z) := zb$ for each $z\in
\bf H$. On the other hand, for each $h\in \bf K$,
$h=h_0i_0+...+h_{2^r-1}i_{2^r-1}$ with $h_j\in \bf R$ there are
identities: \\
$h_j=(-hi_j+ i_j(2^r-2)^{-1} \{ -h
+\sum_{k=1}^{2^r-1}i_k(hi_k^*) \} )/2$ for each $j=1,2,...,2^r-1$, \\
$h_0=(h+ (2^r-2)^{-1} \{ -h +\sum_{k=1}^{2^r-1}i_k(hi_k^*) \} )/2$,
where $r=2$ for quaternions and $r=3$ for octonions.
\par In view of Proposition 10.30 \cite{port69} the mapping
$(a,b)\mapsto a_Lb_R$ from $S^3\times S^3$ on the special orthogonal
group $SO(4)$ is the group surjection with the kernel $\{ (1,1);
(1,-1) \} $, where $S^3 := \{ z\in {\bf H}: |z|=1 \} $ denotes the
multiplicative subgroup in $\bf H$ isomorphic with $SU(2)$. When $f$
is presented by a shortest phrase with imposed Conditions $(1-3)$ of
Definition 1, then a nonproper rotation  associated with $f'$ is
impossible, since to it there corresponds a transformation from the
orthogonal group $O(4)$ with negative determinant, that, in
particular, corresponds to the composition of a proper rotation with
the positive determinant with the functions of the form $z\mapsto
i_k{\tilde z}i_k$, but no any phrase $\nu $ associated with it for
$f'$ while expression through $z$ becomes either four times longer
or it does not satisfy the condition $\partial \nu /\partial {\tilde
z}=0$. On the other hand, $\partial f(z)/\partial z$ satisfying
Condition $(2)$ of Definition 1 belongs to the direct product
$SO(4)\times (0,\infty )$ of the group $SO(4)$ and the
multiplicative group $(0,\infty )$ of positive numbers. From this
the first part of the statement of this theorem follows. \par The
second statement follows from (super)differentiability of $f$ and
from the fact, that due to Condition $(1)$ of Definition 1 the
function $f$ is quaternion holomorphic, consequently, it is locally
analytic by the quaternion variable $z$, hence $a(z)$ and $b(z)$ are
locally analytic on $U$, which is equaivalent to the quaternion
holomorphicity due to Theorems 2.16, 3.10 and Corollary 2.13
\cite{ludoyst}.
\par {\bf 5. Theorem.} {\it A holomorphic function $f$ is pseudoconformal
at a point $\xi $ or on an open domain $U$ in $\bf O$ if and only if
the operator $f'(\xi )$ of its derivative at $\xi $ or at each $\xi
\in U$ correspondingly can be presented as $f'(\xi ).h =
(a((bh)c))d$ for each $h\in \bf O$, where $a(\xi )$, $b(\xi )$,
$c(\xi )$ and $d(\xi )$ are octonion nonzero constants or
holomorphic functions respectively.}
\par {\bf Proof.} At first prove, that $f(\xi )$ is pseudoconformal
at $\xi $ or on $U$ if and only if $f'(\xi )$ is the composition of
the real dilation operator $h\mapsto \lambda h$, where $\lambda >0$,
and $\bf R$-homogeneous $\bf O$-additive operators $Y_{k,m}$ of the
form $Y_{k,m}=\exp (t_{k,m}X_{k,m})h=h+(-1+\cos
(t_{k,m}))(h_ki_k+h_mi_m) +\sin (t_{a,b})(h_mi_k-h_ki_m)$ for each
$h\in \bf O$, $t_{a,b}\in \bf R$ and $0\le a<b\le 7$, where $ \{
i_0,i_1,...,i_7 \} $ are standard generators of $\bf O$. Moreover,
if $f$ is pseudoconformal on $U$, then $\lambda $ and $t_{k,m}$ are
locally analytic functions of $\xi $.
\par In view of the Conditions $(i-iii)$ of Definition 1
and since $Re (h{\tilde v})$ induces the standard scalar product
$(h,v)=h_0v_0+...+h_7v_7$ on the shadow (underlying) Euclidean space
$\bf R^8$ at each point $\xi $, where $f$ is pseudoconformal, the
operator of its derivative $f'(\xi )=(\partial f(z)/\partial
z)|_{z=\xi }$ is uniquelly characterized as the composition of
dilation and a corresponding operator $A\in SO({\bf R},8)$ on $\bf
R^8$. Indeed, when $f$ is presented by a shortest phrase with the
imposed conditions, then nonproper rotations associated with $f'$
are impossible, since the corresponding to them operators from
$O({\bf R},8)$ with negative real determinants, that in particular
corresponds to a composition of a proper rotation with a mapping of
the form $z\mapsto i_k{\tilde z}i_k$ for some generator $i_k$, but
an associated with its phrase $\nu $ for $f'$ while expression
through $z$ becomes longer eight times or it does not satisfy the
condition $\partial \nu /\partial {\tilde z}=0$.
\par The Lie group $SO({\bf R},8)$ consists of operators
$A=\exp (B)$, where $B$ are $8\times 8$ skew-symmetric real matrices
(see, for example, \cite{gogros}). That is, the basis of the Lie
algebra $so({\bf R},8)$ consists of generators $B_{k,m}=-B_{k,m}^T$,
where $0\le k<m\le 7$, $B_{k,m}=E_{k,m}-E_{m,k}$, $B^T$ denotes the
transposed matrix $B$, $E_{k,m}$ is the $8\times 8$ matrix with $1$
at the position $(k,m)$ on the cross of a row with number $k$ and a
column with number $m$ and all its others elements are zero.
Therefore, to $B_{k,m}$ over $\bf R$ there corresponds the generator
$X_{k,m}=(E_{k,m}i_ki_m-E_{m,k}i_mi_k)$ over $\bf O$. Then
$E_{k,m}E_{l,p}=\delta _{m,l}E_{k,p}$, where $\delta _{k,p}=1$ for
each $k=p$ and $\delta _{k,p}=0$ for each $k\ne p$ is the Kronecker
delta-function. Thus $X_{k,m}^2=E_{k,k}+E_{m,m}$ for each $k<m$.
Therefore, $X_{k,m}^3=X_{k,m}$ and $\exp(t_{k,m}X_{k,m})h=h+(-1+\cos
(t_{k,m}))(h_ki_k+h_mi_m) +\sin (t_{a,b})(h_mi_k-h_ki_m)$, where $z=
h_0i_0+...h_7i_7$, $h_0,...,h_7\in \bf R$. \par There are projection
operators $\pi _j(h):=h_j=(-hi_j+ i_j(2^r-2)^{-1} \{ -h
+\sum_{l=1}^{2^r-1}i_l(hi_l^*) \} )/2$ for each $j=1,2,...,2^r-1$, \\
$\pi _0(h):=h_0=(h+ (2^r-2)^{-1} \{ -h
+\sum_{l=1}^{2^r-1}i_l(hi_l^*) \} )/2$ on ${\cal A}_r$, in
particular, for ${\bf O}={\cal A}_3$. If $f$ is holomorphic on $U$,
then it is locally analytic, consequently, $f'(z).h$ is locally
analytic and inevitably $\lambda $ and $t_{k,m}$ for each $0\le
k<m\le 7$ are locally analytic functions of $h_j$ and hence of $h$.
\par From this follows, that if a holomorphic function $f$ has a
derivative $f'(\xi )$ at a point $\xi $ of an open domain $U$ in
$\bf K$ such that $f'(\xi ).h=(a((bh)c)d$ for each $h\in \bf O$,
where $a\ne 0$, $b\ne 0$, $c\ne 0$ and $d\ne 0$ are octonions, then
$f$ is pseudoconformal at a point $\xi $. If $a(\xi )$,...,$d(\xi )$
are holomorphic, then by Theorems 3.9 and 3.10 \cite{ludfov} $f(\xi
)$ is also holomorphic and hence locally analytic by $\xi $.
\par Consider products of octonions:
\par $(i)$ $(a(z_0 + z_ll))b = az_0b + (z_la{\tilde b})l$,
\par $(ii)$ $((al)(z_0+z_ll))(bl) = - ({\tilde b}a{\tilde z}_0) -
(b{\tilde z}_la)l$, \par $(iii)$  $(a(z_0+z_ll))(bl) =  - ({\tilde
b}z_la) + (baz_0)l$ and \par $(iv)$ $((al)(z_0+z_ll))b = - ({\tilde
z}_lab) + (a{\tilde z}_0{\tilde b})l $ \\
for each $a$, $b$, $z_0$ and $z_l\in \bf H$, where $l=i_4$ is the
generator of the doubling procedure of $\bf O$ from $\bf H$. Let $M$
and $N\in \bf O$ with $Re (M)=0$ and $Re (N)=0$, $|M|=1$, $|N|=1$,
$Re (M{\tilde N})=0$, that is, $M\perp N$, then take $K=MN$ and the
algebra over $\bf R$ with generators $\{ 1, M, N, MN \} $ is
isomorphic with $\bf H$ and contained in $\bf O$. An arbitrary
octonion $z$ can be written uniquelly in the form $z=z_0+z_ll$,
where $z_0$ and $z_l$ are quaternions. From Formulas $(i-iv)$,
Theorem 4 and the preceding proof it follows, that $f'(\xi )$ can be
written in the form $f'(\xi ).h = (a((bh)c))d$ for each $h\in \bf
O$, where $a(\xi )$, $b(\xi )$, $c(\xi )$ and $d(\xi )$ are nonzero
octonions, since each one-parameter pseudoconformal transformation
group with a real parameter can be written in such form. If $f$ is
pseudoconformal on $U$, then $f$ is locally analytic and nonzero
octonion valued functions $a(\xi )$, $b(\xi )$, $c(\xi )$ and $d(\xi
)$ can be chosen locally analytic.
\par {\bf 5.1. Note.} The form $<ah_1b,ah_2b>$ for each
$h_1, h_2\in \bf H$ with quaternions $a\ne 0$ and $b\ne 0$ is $\bf
R$-homogeneous of degree $4$ relative to the substitutions $a\mapsto
\lambda a$ and $b\mapsto \lambda b$, where $\lambda \in \bf R$,
while the form $<(a((bh_1)c)d, (a((bh_2)c)d>$ with $h_1\in \bf O$
and $h_2\in \bf O$ is $\bf R$-homogeneous of degree $8$ relative to
the substitutions $a\mapsto \lambda a$, $b\mapsto \lambda b$,
$c\mapsto \lambda c$ and $d\mapsto \lambda d$. By the totality of
all variables the first form is $\bf R$-homogeneous of degree $6$
while the second of degree $10$. The restrictions on $\bf R^6$ and
$\bf R^{10}$ give multilinear forms from them, that is, forms of
Finsler geometry (see also about Finsler geometry \cite{asanov}).
Thus pseudoconfromal mappings over $\bf H$ and $\bf O$ can serve as
models for noncommutative analogs of Finsler geometry.
\par {\bf 6. Theorem.} {\it Let functions $f: U\to \bf K$ and
$g: W\to \bf K$ be quaternion or octonion holomorphic in open
domains $U$ and $W$ in $\bf K$, where $W\supset f(U)$. \par  1. If
$f$ is pseudoconformal at a point $\xi $ and $g$ is pseudoconformal
at a point $\zeta = f(\xi )$, then their composition $g\circ f$ is
pseudoconformal at $\xi $. \par 2. If $f$ is pseudoconformal on $U$
and $g$ is pseudoconformal on $W$, then $g\circ f$ is
pseudoconformal on $U$. \par 3. If there exists an inverse function
$f^{-1}: P\to U$, where $P=f(U)$ is an open domain in $\bf K$, then
from the pseudoconformity of $f$ at a point $\xi \in U$ or on a
domain $U$ it follows the pseudoconformity of the inverse function
$f^{-1}$ at a point $f(\xi )$ or on a domain $P$ respectively.}
\par {\bf Proof.} The composition of quaternion or
octonion holomorphic functions is quaternion or octonion
holomorphic, since $\partial g(\zeta )/\partial {\tilde {\zeta
}}=0$, $\partial f(\xi )/\partial {\tilde {\xi }}=0$ and hence
$\partial g(f(\xi ))/\partial {\tilde {\xi }}= (\partial
g(y)/\partial y)|_{y=\zeta }.(\partial f(\xi )/\partial {\tilde {\xi
}})=0$. On the other hand, if there exists the inverse mapping
$f^{-1}$, then \par $(\partial f^{-1}(\zeta )/\partial {\tilde
{\zeta }}).h = - (f^{-1}(\zeta ))[(\partial f(\xi )/\partial {\tilde
{\xi }}).h)(f^{-1}(\zeta ))]=0$ \\ for each $h\in \bf K$, where
$\zeta =f(\xi )$. Thus $f^{-1}$ is holomorphic at $\zeta $. \par
There are identities: $(\partial g(f(\xi ))/\partial {\xi }).h =
(\partial g(y)/\partial y)|_{y=\zeta }.[(\partial f(\xi )/\partial
{\xi }).h]$ and $(\partial f^{-1}(\zeta )/\partial  {\zeta }).h = -
(f^{-1}(\zeta ))[(\partial f(\xi )/\partial {\xi }).h)(f^{-1}(\zeta
))]$ for each $h\in \bf K$. Then from Condition $(2)$ of Definition
1 for $f$ and $g$ it follows, that it is accomplished for $g\circ f$
and $f^{-1}$ correspondingly.
\par {\bf 7. Corollary.} {\it Each branch of the quaternion or
octonion function of the logarithmic function $Ln$ is
pseudoconformal besides zero. If a function $f(z)$ is
pseudoconformal at a point $\xi $ and $f(\xi )\ne 0$, then $1/f(z)$
is also pseudoconformal at $\xi $.}
\par {\bf Proof.} The logarithmic function is the inverse function
of the exponential function and $Ln$ has the countable family of
branches such that $Ln$ is holomorphic on ${\bf K}\setminus \{ 0 \}
$ (see \S 3.7 \cite{ludfov}). Then $Ln$ is pseudoconformal on ${\bf
K}\setminus \{ 0 \} $ due to Theorem 6. \par The function
$g(z):=1/z$ is pseudoconformal on ${\bf K}\setminus \{ 0 \} $ due to
Corollary 3. Since $1/f(z)$ is the composition $g\circ f$, then it
is pseudoconformal by Theorem 6.
\par {\bf 8. Definitions.}  For an open
domain $U$ in $\bf K^n$ with $n\in \bf N$, $n>1$, a function
$f=(\mbox{ }_1f,...,\mbox{ }_mf): U\to \bf K^m$ is called
pseudoconformal at a point or on the domain $U$, if $\mbox{ }_sf$ is
pseudoconformal by each quaternion or octonion variable $\mbox{
}_jz$ at a point or on the domain, where $z=(\mbox{ }_1z,...,\mbox{
}_nz)\in U$, $\mbox{ }_jz\in \bf K$, $s=1,...,m$, $j=1,...,n$. Let
$M$ be a connected compact manifold over $\bf K$ with the atlas $At
(M)=\{ (U_j,\phi _j): j=1,...,q \} $, such that the connecting
mappings $\phi _i\circ \phi _j^{-1}$ are pseudoconformal for charts
$(U_i,\phi _i)$ and $(U_j,\phi _j)$ with $U_i\cap U_j\ne \emptyset
$, $U_j$ is open in $M$, $\phi _j(U_j)$ is open in $\bf {\hat K}^n$
or in ${\bf {\hat K}^{n-1}}\times \bf K_+$, $\bigcup_{j=1}^qU_j=M$,
$1\le n\in \bf N$, where ${\bf K_+} := \{ z\in {\bf K}: Re (z)\ge 0
\} $, $\bf \hat K$ denotes the one point Alexandroff
compactification of $\bf K$ which is homeomorphic with the unit
sphere $S^{2^r}$ of dimension $2^r$ over $\bf R$.
\par Suppose that $M$ and $N$ are two manifolds with pseudoconformal
connecting mappings of their atlaces $At (M)=\{ (U_j,\phi _j):
j=1,...,q \} $ and $At (N)= \{ (W_j, \psi _j) \} $. Shortly such
manifolds we shall call pseudoconformal. A function $f: M\to N$ we
call pseudoconformal, if $\psi _i\circ f\circ \phi _j^{-1} =
:f_{i,j}$ is pseudoconformal for each $f(U_j)\cap W_j\ne \emptyset $
on its domain $Q_{i,j}$ of definition, moreover, for boundary points
we suppose, that each $f_{i,j}$ has a pseudoconformal extension on a
neigborhood of $Q_{i,j}$ in $\bf K^n$. On the family $P(M,N)$ of all
such pseudoconformal mappings we introduce the metric $$\rho
_M(f,g):=\rho (f,g):=\max_{z\in M} [d(f(z),g(z))+ \max_{({h\in \bf
K^n}, |h|=1; i, j)} |f_{i,j}'(z).h- g_{i,j}'(z).h|],$$ where $d$
denotes the metric on $N$, where $h\in \bf K^n$ and the tangent
space $T_zM$ is isomorphic with $\bf K^n$, $f_{i,j}'(z) : T_zM\to
T_{f(z)}N$. In particular, for $M=N$ denote $P(M) := P(M,M)$.
\par {\bf 9. Definition.} A family $\cal F$ of functions $f$
on an open domain $U$ in $\bf K^n$ is called normal, if each $f\in
\cal F$ is $\bf K$ holomorphic on $U$, moreover, at each point $\xi
\in U$ it is by each variable $\mbox{ }_jz$ either pseudoconformal
or $(\partial f(z)/\partial \mbox{ }_jz)|_{\xi =z}=0$, while for
each infinite sequence $\{ f_n: n\in {\bf N} \} \subset \cal F$ for
each canonical closed bounded subdomain $V$ in $U$ there exists a
subsequence $ \{ f_{n_k}: k\in {\bf N} \} $ converging uniformly
relative to the metric $\rho $ from Definition 8 for $M=V$ to some
limit function, such that it may also be either finite or infinite
constant, where $n_1<n_2<n_3...$. If a family $\cal F$ is normal in
some open ball with the centre at a point $z\in \bf K^n$, then it is
called normal at $z$.
\par {\bf 10. Lemma.} {\it If a sequence
$ \{ f_n: n\in {\bf N} \} $ of $\bf K$ pseudoconformal functions is
uniformly converging relative to the metric $\rho $ on open subsets
$W$ which form a covering family of $U$ to some finite function $f$,
where $U$ is open in $\bf K$, then $f$ is at each point $\xi \in U$
either $\bf K$ pseudoconformal or $f'(\xi )=0$.}
\par {\bf Proof.} Since a sequence
$ \{ f_n: n\in {\bf N} \} $ is uniformly converging relative to the
$C^1$ uniformity generated by the $C^1$ compact-open topology
relative to $(z,{\tilde z})$-superdifferentiation, then due to
Theorems 4 and 5 $f:=\lim_{n\to \infty }f_n$ belongs to the space
$C^1_{(z,{\tilde z})}(W,{\bf K})$ of all $(z,{\tilde
z})$-differentiable functions on $W$ with values in $\bf K$. At the
same time either $f'(\xi ).h=a(\xi )hb(\xi )$ or $f'(\xi ).h=(a(\xi
)((b(\xi )h)c(\xi ))e(\xi )$ for each $h\in \bf K$ and each $\xi \in
U$ for ${\bf K}=\bf H$ or ${\bf K}=\bf O$ respectively, where $a(\xi
)$, $b(\xi )$, $c(\xi )$, $e(\xi )$ are continuous $\bf K$-valued
functions by $\xi $. Thus $\partial f(z)/\partial {\tilde z}=0$ for
each $z\in W$, consequently, due to Theorems 2.16 and 3.10
\cite{ludoyst,ludfov} $f$ is $\bf K$ holomorphic and locally
analytic by $z$. Hence $a(\xi )$, $b(\xi )$, $c(\xi )$, $e(\xi )$
can be chosen also holomorphic by $\xi $.
\par {\bf 11. Lemma.} {\it If $\cal F$ is a normal family and a sequence
$ \{ f_n: n\in {\bf N} \} $ converges to a function $f$ from
Definition 9 uniformly relative to the metric $\rho $ on open
subsets $W$ in $U$ which form a covering of $U$ to a function $f$,
where $U$ is open in $\bf K^n$, then $f$ is holomorphic and for each
$j=1,...,n$ at each $z\in U$ a function $f$ is either
pseudoconformal by $\mbox{ }_jz$ or $\partial f(z)/\partial \mbox{
}_jz=0$.}
\par {\bf Proof.} In view of Lemma 10 a function $f$ is $\bf K$
holomorphic. For the proof it is sufficient to consider a behavour
of $f$ by each variable $\mbox{ }_jz$. That is, we reduce the final
part of the proof to the case $n=1$. If $f$ is constant in a
neighborhood of a point $z$, then $f'(z)=0$. In view of Theorems 4
and 5 its derivative $f_n'(z)$ can be written with the help of
functions $a_n(z), b_n(z), c_n(z), e_n(z)$, moreover, $a_n(z)\ne
0,..., e_n(z)\ne 0$, if $f_n'(z)\ne 0$, where each function
$a_n(z),...,e_n(z)$ is holomorphic. Since $f_n'.h$ converges to
$f'.h$ on $W$ uniformly relative to the $C^0$-uniformity, generated
by the $C^0$ compact-open topology, then for each $f'$ there exist
functions $a(z),...,e(z)$ such that $f'(z).h=a(z)hb(z)$ in the
quaternion case and $f'(z).h=\{ a(z)b(z)hc(z)e(z) \} _{q(5)} $ in
the octonion case for each $h\in \bf K$.
\par {\bf 12. Remark.} In accordance with Theorem 1.6.7 \cite{span}
a topological space $X$ is $n$-connected if and only if it is
arcwise connected and $\pi _k(X,x)$ is trivial for each base point
$x\in X$ and each $k$ such that $1\le k\le n$.
\par Denote by $Int (U)$ the interior of a subset $U$
in a topological space $X$, denote also $cl (U)=\bar U$ the closure
of $U$ in $X$. For a subset $U$ in ${\bf K}= {\bf K}_r$ let $\pi
_{s,p,t}(U):= \{ u: z\in U, z=\sum_{v\in \bf b}w_vv,$ $u=w_ss+w_pp
\} $ for each $s\ne p\in \bf b$, where $t:=\sum_{v\in {\bf
b}\setminus \{ s, p \} } w_vv \in {\bf K}_{r,s,p}:= \{ z\in {\bf
K}_r:$ $z=\sum_{v\in \bf b} w_vv,$ $w_s=w_p=0 ,$ $w_v\in \bf R$
$\forall v\in {\bf b} \} $, where ${\bf b} := \{
i_0,i_1,...,i_{2^r-1} \} $ is the family of standard generators of
the algebra $\bf K$, $i_0=1$, $r=2$ for ${\bf K}_2=\bf H$, $r=3$ for
${\bf K}_3=\bf O$. That is, geometrically $\pi _{s,p,t}(U)$ is the
projection on the complex plane ${\bf C}_{s,p}$ of the intersection
of $U$ with the plane ${\tilde \pi }_{s,p,t}\ni t$, ${\bf C}_{s,p}
:= \{ as+bp:$ $a, b \in {\bf R} \} $, since $sp^*\in {\hat b}:={\bf
b}\setminus \{ 1 \} $. Recall that in \S \S 2.5-7 \cite{ludfov} for
each continuous function $f: U\to \bf K$ there was defined the
operator $\mbox{ }_j{\hat f}$ by each variable $\mbox{ }_jz\in \bf
K$.
\par The following criteria permit to characterize normal families
of functions. Henceforth it is supposed that a domain $U$ in $\bf
K^n$ has the property that each projection $p_j(U)=:U_j$ is
$2^r-1$-connected; $\pi _{s,p,t}(U_j)$ is simply connected in $\bf
C$ for each $k=0,1,...,2^r-1$, $s=i_{2k}$, $p=i_{2k+1}$, $t\in {\bf
K}_{r,s,p}$ и $u\in {\bf C}_{s,p}$ for which there exists $z=u+t\in
U_j$, where $e_j = (0,...,0,1,0,...,0)\in \bf K^n$ is the vector
with $1$ on the $j$-th place, $p_j(z) = \mbox{ }_jz$ for each $z\in
\bf K$, $z=\sum_{j=1}^n\mbox{ }_jz e_j$, $\mbox{ }_jz\in \bf K$ for
each $j=1,...,n$, $n\in {\bf N} := \{ 1,2,3,... \} $.
\par {\bf 13. Theorem.} {\it Let a domain $U$ in ${\bf K}^n$
satisfy Conditions of Note 12, while a family $\cal F$ consists of
functions $\bf K$ holomorphic on $U$, such that each $f\in \cal F$
for each $\xi \in Int (U)$ by the variable $\mbox{ }_jz$ is either
pseudoconformal or $(\partial f(z)/
\partial \mbox{ }_jz)_{\xi =z}=0$. Let in addition on
$U$ the family $\cal F$ be uniformly equibounded:
\par $(i)$ $\sup_{(z\in Int (U); f\in {\cal F})}
|f(z)| =:c<\infty $, \\  then the family $\cal F$ is normal.}
\par {\bf Proof.} If a family $\cal F$ is finite, then it is nothing
to prove. Suppose that $\cal F$ is infinite. Let $V$ be a canonical
closed bounded subset in $U$, that is, $cl (Int (V)) = V$. For each
internal point $z_0$ in $V$ due to Theorems 3.9 and 3.10
\cite{ludfov} there exists $0<R<\infty $, such that each $f(z)\in
\cal F$ has a decomposition in the converging power series in the
ball with centre $z_0$ and radius $R$:
\par $(ii)$ $f(z) = \sum_{|m|=0}^{\infty } \phi _{f,|m|} (z-z_0)$,
where
\par $(iii)$ $\phi _{f,c}(z-z_0) =
\sum_{|m|=c} \{ (a_{f,m}, (z-z_0)^m) \} _{q(2v)}$, \\
$c=0,1,2,...$, $a_{f,m,j}\in \bf K$, $m = (m_1,...,m_v)$,
$m_j=0,1,2,...$, $a_{f,m}=(a_{f,m,1},..., a_{f,m,v})$, $|m| :=
m_1+...+m_v$, $v\in \bf N$, $\{ (a_{f,m},z^m) \} _{q(2v)}:= \{
a_{f,m,1}z^{m_1}...a_{f,m,v}z^{m_v} \} _{q(2v)}$. \\
There exists  $R>0$ one for all $f\in \cal F$ due to Condition
$(i)$, since the family $\phi _{f,c}$ satisfies the inequality:
\par $(iv)$ $|\phi _{f,c}(z-z_0)|\le J_0|z-z_0|^c/R^c$
for each $|z-z_0|<2R/3$, where $J_0$ is a constant independent from
$f, c$, since
$$(v)\quad f(z)= \sum_{k=0}^{\infty }\phi _k(z),$$
$$\mbox{where }\phi _k(z):= (2\pi )^{-1}
(\int_{\psi }f(\zeta )((\zeta -a)^{-1}(z-a))^k(\zeta -a)^{-1}d\zeta
)M^*$$ for each $\bf K$ holomorphic function $f$, where a point $z$
is encompassed by a rectifiable loop $\psi $ with a directing purely
imaginary quaternion or octonion $M\in {\bf K}$, $Re (M)=0$,
$|M|=1$. \par Into an infinite open covering by interiors of balls
$B({\bf K^n},z_0,R/2)$ of the set $V$ it can be refined a finite
subcovering due to compactness of $V$. Centres of balls of this
finite covering we denote by $z_1,...,z_w$. Each function $\phi
_{f,c}(z-z_j)$ is the homogeneous polynomial, containing in its
decompostion $(iii)$ a finite number of terms.
\par Each ball $B({\bf K^p},z,R)$ is compact for each $p\in \bf N$
and $0<R<\infty $, while the product of an arbitrary family of
compact spaces is compact in the product topology. The topological
space $X$ is called countably compact, if from each open countable
covering of $X$ there can be extracted a finite subcovering. In view
of Theorem 3.10.3 \cite{eng} the space $X$ is countably compact if
and only if each countable infinite subset $A$ in $X$ has in $X$ a
strictly limit point, that is, in each its neighborhood lies an
infinite number of points from $A$. Due to Theorem 3.10.31
\cite{eng} the sequential compactness is equivalent to the
compactness in the class of $T_1$-spaces with the first countability
axiom. On the other hand, due to Theorem 3.10.35 \cite{eng} the
product of an arbitrary countable family of sequentially compact
spaces is sequentially compact. \par From the proof of Theorems 2.7
and 3.10 \cite{ludfov} and Condition $(i)$ it follows, that the
family $\{ f'(z).h : f\in {\cal F} \} $ is uniformly equibounded on
$V\times B({\bf K^n},0,1)$. Due to conditions on the domain $U$ and
the family $\cal F$ the family of functions $\cal F$ is uniformly
equibounded on each compact subset $V$ in $Int (U)$. Then due to the
separability of $\bf K$ as the metric space the consideration of the
family $\{ a_{f,m,v} \} $ gives that the family of functions $\cal
F$ contains a convergent subsequence $ \{ f_n: n\in {\bf N} \} $
relative to the metric $\rho $ for $V$, since the family of indices
$\{ (f,m,v) \} $ is countable, while derivatives of the considered
series also converge uniformly in each ball of the finite covering.
Due to Lemma 11 for each point $z\in W$ the function $f:=\lim_{n\to
\infty }f_n$ is either $\bf K$ pseudoconformal or holomorphic with
$f'(z)=0$.
\par For the susequent proceedings there are necessary
noncommutative analogs of the main theorem of mathematical analysis.
For this there are given below theorems. For the extracting of a
univalent concrete branch of a path integral introduce the following
notations and definitions.
\par {\bf 14. Definitions and Notations.} 1. Let $\bf 1$ denotes
the unit operator on ${\cal A}_r$, that is, ${\bf 1}(h)=h$ for each
$h\in {\cal A}_r$, where ${\cal A}_r$ denotes the Cayley-Dickson
algebra with $2^r$ standard generators $ \{ i_0,i_1,...,i_{2^r-1} \}
$, $i_0=1$, $i_p^2=-1$ for each $p\ge 1$, $i_pi_s=-i_si_p$ for each
$1\le s\ne p\ge 1$. Consider also the operator of conjugation
${\tilde {\bf 1}}(h)={\tilde h}$ for each $h\in {\cal A}_r$.
\par 2. We shall distinguish symbols: $e = {\bf 1}(1)$ и ${\tilde e} =
{\tilde {\bf 1}}(1)$, where $1\in {\cal A}_r$.
\par 3. $\bf 1$ and ${\tilde {\bf 1}}$ we shall not identify
or mutually suppress with each other or with symbols of the type
${\bf 1} {\tilde {\bf 1}}$, or with others basic symbols; while
symbols $e$ and ${\tilde e}$ we shall not identify or suppress with
each other or with expressions of the form $e{\tilde e}$ or with
others basic symbols.
\par 4. $(z^pz^q)$ is identified with $z^{p+q}$, while
$({\tilde z}^p{\tilde z}^q)$ is identified with ${\tilde z}^{p+q}$,
where $p$ and $q$ are natural numbers, $p, q=1,2,3,...$. In each
word we suppose that in each fragment of the form $...(z^p{\tilde
z}^q)...$ the symbol $z^p$ is displayed on the left of the symbol
${\tilde z}^q$.
\par 5. If $w$ is some word or
a phrase, $c\in \bf R$, then the phrases $cw$ and $wc$ are
identified; if $b, c\in \bf R$, then the phrases $b(cw)$ and $(bc)w$
are identified.
\par 6. While the conjugation we shall use the identifications
of words: $(a({\bf 1}b))^*=({\tilde b}{\tilde {\bf 1}}){\tilde a}$,
$(a(eb))^*=({\tilde b}{\tilde e}){\tilde a}$ ; while the inversion
gives the identifications: $(a({\bf 1}b))^{-1}=(b^{-1}{\bf
1})a^{-1}$, $(a({\tilde {\bf 1}}b))^{-1}=(b^{-1}{\tilde {\bf
1}})a^{-1}$, $(a(eb))^{-1}=(b^{-1}e)a^{-1}$, $(a({\tilde
e}b))^{-1}=(b^{-1}{\tilde e})a^{-1}$, where $ab\ne 0$, $a, b\in
{\cal A}_r$.
\par 7. A word $\{ w_1...w_{k-1} (\{ a_1...a_l \} )_{q(l)} w_{k+1}...w_j
\} _{q(j)}$ is identified with a word \\ $\{
w_1...w_{k-1}w_kw_{k+1}...w_j \} _{q(j)}$, as well as with words $
с\{ w_1...w_{k-1}bw_{k+1}...w_j \} _{q(j)}$, \\ $ \{ (с
w_1)...w_{k-1}bw_{k+1}...w_j \} _{q(j)}$, where $w_k= (\{ a_1...a_l
\} )_{q(l)}$, $a_1,...,a_l\in {\cal A}_r$, $l, k, j\in \bf N$,
$cb=w_k$, $c\in \bf R$, $b\in {\cal A}_r$.
\par 8. Phrases are sums of words. Sums may be finite or infinite
countable. A phrase $aw+bw$ is identified with a word $(a+b)w$, a
phrase $wa+wb$ is identified with a word $w(a+b)$ for each words $w$
and constants $a, b \in {\cal A}_r$, since $c=a+b\in {\cal A}_r$.
\par 9. Symbols of functions, that is, corresponding to them phrases
finite or infinite, are defined with the help of the following
initial symbols: constants from ${\cal A}_r$, $e$, ${\tilde e}$,
$z^p$, ${\tilde z}^p$, where $p\in \bf N$. Symbols of operator
valued functions are composed from the symbols of the set: $\{ {\bf
1}, {\tilde {\bf 1}}; e, {\tilde e}; z^p, {\tilde z}^p: p\in {\bf N}
\} $, where the symbol ${\bf 1}$ or ${\tilde {\bf 1}}$ is present in
a phrase. At the same time a phrase corresponding to a function or
an operator valued function can not contain words consisting of
constants only, that is, each word $w_1w_2...w_k$ of a phrase must
contain not less than one symbol $w_1, w_2,...,w_k$ from the set $\{
e, {\tilde e}; z^p, {\tilde z}^p: p\in {\bf N} \} $. For a phrase
$\nu $ and a void word $\emptyset $ put $\nu \emptyset =\emptyset
\nu =\nu $.
\par 10. For symbols from the set $\{ {\bf 1}, {\tilde
{\bf 1}}; e, {\tilde e}; z^p, {\tilde z}^p: p\in {\bf N} \} $ define
their lengths: $l(0)=0$, $l(a)=1$ for each $a\in {\cal A}_r\setminus
\{ 0 \} $, $l({\bf 1})=1$, $l({\tilde {\bf 1}})=1$, $l(e)=1$,
$l({\tilde e})=1$, $l(z^p)=p+1$, $l({\tilde z}^p)=p+1$ for each
natural number $p$. A lenght of a word is the sum of lengths of
composing it symbols. Analogously there are considered words and
phrases for several variables $\mbox{ }_1z,...,\mbox{ }_nz$, also
the function of a lenght of a word, where symbols corresponding to
different indices $v=1,...,n$ are different. \par Henceforth, we
consider phrases subordinated to the following conditions. Either
all words $w$ of a given phrase $\nu $ contain each of the symbols
$\mbox{ }_ve$, $\mbox{ }_v{\tilde e}$, $\mbox{ }_v\bf 1$ or $\mbox{
}_v{\tilde {\bf 1}}$ with the same finite constant multiplicity,
which may be dependant on the index $v=1,...,n$ or on a symbol
itself; or each word $w$ of $\nu $ contains neither $\mbox{ }_ve$
nor $\mbox{ }_v{\tilde e}$ besides words $w$ with one $\mbox{ }_ve$
and without any $\mbox{ }_vz^l$ with $l\ne 0$ or $w$ with one
$\mbox{ }_v{\tilde e}$ without any $\mbox{ }_v{\tilde z}^l$ with
$l\ne 0$. If a phrase $\nu $ in each its word $w$ contains $\mbox{
}_v\bf 1$ with multiplicity $p$ or $\mbox{ }_v{\tilde {\bf 1}}$ with
multiplicity $k$, then it is supposed to be a result of partial
differentiation $\partial ^{p+k}f(z,{\tilde z})/\partial \mbox{
}_vz^p \partial \mbox{ }_v{\tilde z}^k$; while the case of $\mbox{
}_ve$ with multiplicity $p$ or $\mbox{ }_v{\tilde e}$ with
multiplicity $k$ is supposed to be arising from $(\partial
^{p+k}f(z,{\tilde z})/\partial \mbox{ }_vz^p \partial \mbox{
}_v{\tilde z}^k).(1^{\otimes (p+k)})$. A presence of at least one of
the symbols $\mbox{ }_v\bf 1$ or $\mbox{ }_v{\tilde {\bf 1}}$ in a
phrase characterize an operator valued function apart from a
function.
\par 11. Supply the space $\cal P$ of all phrases subordinated
to the restrictions given above over ${\cal A}_r$ with the metric by
the formula:
$$ d(\nu ,\mu )=\sum_{j=0,1,2,...} l(\nu _j,\mu _j)b^j,$$
where $b$ is a fixed number, $0<b<1$, $\nu , \mu \in \cal P$,
\par $\nu =\sum_{j=0,1,2,...} \nu _j$, a function or an operator
valued function $f_j(z,{\tilde z})$ corresponding to $\nu _j$ is
homogeneous of a degree $j$, that is, $f_j(tz,t{\tilde
z})=t^jf(z,{\tilde z})$ for each $t\in \bf R$ and $z\in {\cal A}_r$,
$\nu _j=\sum_k \nu _{j,k}$, where $\nu _{j,k}$ are words with
corresponding homogeneous functions or operator valued functions
$f_{j,k}$ of degree $j$, $l(\nu _j,\mu _j):=\max_{k,s} l(\nu _{j,k};
\mu _{j,s})$, where $l(\nu _{j,k}; \mu _{j,s})=0$ for $\nu
_{j,k}=\mu _{j,s}$, $l(\nu _{j,k}; \mu _{j,s})=\max [l(\nu _{j,k});
l(\mu _{j,s})]$ for $\nu _{j,k}$ not equal to $\mu _{j,s}$.
\par 12. To each phrase $\nu $ there corresponds a continuous
function $f$, if a phrase is understood as a sequence of finite sums
of composing it words, such that the corresponding sequence of
functions converges. Consider the space $C^0_{\cal P}(U,X)$,
consisting of pairs $(f, \nu )$, $f\in C^0(U,X)$, where $X={\cal
A}_r^n$ or $X=X_t$ is the space of $\bf R$-polyhomogeneous
polyadditive operators on ${\cal A}_r^m$, $m\in \bf N$, $\nu \in
\cal P$ such that to a phrase $\nu $ there corresponds a continuous
function or operator valued function $f$ on a domain $U$ in ${\cal
A}_r^n$ with values in $X$, where $t=(t_1,t_2)$, $0\le t_1, t_2 \in
\bf Z$, $X_t$ is the ${\cal A}_r$-vector space generated by
operators of the form $A= \mbox{ }_1A\otimes ... \otimes \mbox{
}_{t_1+t_2}A$, $\mbox{ }_jA\in span_{{\cal A}_r}\{ \mbox{ }_1{\bf
1},...,\mbox{ }_n{\bf 1} \} $ for $j=1,...,t_1$; $\mbox{ }_jA\in
span_{{\cal A}_r} \{ \mbox{ }_1{\tilde {\bf 1}},...,\mbox{
}_n{\tilde {\bf 1}} \} $ for $j=t_1+1,...,t_1+t_2$. Introduce on
$C^0_{\cal P}(U,X)$ the metric:
$${\cal D}((f,\nu ); (g,\mu )):=
\sup_{z\in U} \| f(z,{\tilde z})-g(z,{\tilde z}) \| +d(\nu ,\mu ),$$
where to a phrase $\nu $ there corresponds a function $f$, while to
a phrase $\mu $ there corresponds a function $g$, where the norm of
a $m$-times $\bf R$-polyhomogeneous polyadditive operator is given
by the formula $\| S \| := \sup_{(\| h_1 \| =1,...,\| h_m \| =1)} \|
S(h_1,...,h_m) \| /[\| h_1 \| ... \| h_m \| ]$, where for a
function, that is, for $m=0$, there is taken the absolute value of
the function $|f|=\| f \| $ instead of the operator norm.
\par {\bf 15. Note.} For $r=\infty $ or $r=\Lambda $
with $card (\Lambda )\ge \aleph _0$ the variables $z$ and $\tilde z$
are algebraically independent over ${\cal A}_r$ due to Theorem 3.6.2
\cite{ludfov}. Since $A_m$ with $m<r$ has an algebraic embedding
into ${\cal A}_r$, then it is possible to consider functions
$f(z,{\tilde z})$ on domains $U$ in ${\cal A}_m$, which are
restrictions $f=g|_U$ of functions $g_f(z,{\tilde z})$ on domains
$W$ in ${\cal A}_r$, for which $W\cap {\cal A}_m=U$. \par In view of
Definition 14 symbols of functions, that is, corresponding to them
phrases finite or infinite are defined with the help of the initial
symbols. From these phrases are defined more general analytic
functions, in particular, $\exp (z)=e^z$, $Ln (z)$, $z^a$, $a^z$,
etc. For this local decompositions into power series are used for
locally analytic functions, as well as the Stone-Weierstrass
theorem, stating that the polynomials are dense in the space of
continuous functions on a compact canonical closed subset in $\bf
R^n$, hence also in ${\cal A}_r$ with $0<r<\infty $ (see
\cite{ludoyst,ludfov}).
\par For nonassociative algebras with $m\ge 3$ an order of
multiplication is essential, which is prescribed by an order of
brackets or by a vector $q(s)$ indicating on an order of
multiplications in a word. For example, words $(az^p)(z^qb)$ and
$a(z^{p+q}b)$ or $(az^{p+q})b$ for $a, b \in {\cal A}_m\setminus \bf
R$, $m\ge 3$, are different, where $p\ne 0$ and $q\ne 0$ are natural
numbers.
\par {\bf 16. Lemma.} {\it If two functions or operator valued functions
$f$ and $g$ are bounded on $U$ (see Definition 14.12), then ${\cal
D}((f,\nu ); (g,\mu ))<\infty $. In particular, if $U$ is compact,
then ${\cal D}((f,\nu );(g,\mu ))<\infty $ for each $(f,\nu ),
(g,\mu )\in C^0_{\cal P}(U,X)$.}
\par {\bf Proof.} For each degeree of homogeneouity $j$ each word
$\nu _{j,k}$ consists no more, than of $b(j+\omega )$ symbols, where
$b=2$ for $\bf H$, $b=3$ for ${\cal A}_r$ with $r\ge 3$,
consequently, a number of such words is finite and $l(\nu _{j,k})\le
s(j+\omega )<\infty $, since a multiplicity $s_v$ or $u_v$, or
$p_v$, or $y_v$ of an apperance of $\mbox{ }_ve$ or $\mbox{
}_v{\tilde e}$, or $\mbox{ }_v\bf 1$, or $\mbox{ }_v{\tilde {\bf
1}}$ is finite and constant in all words of a given phrase may be
besides finite number of minor words associated with constants,
though it may depend on the index $v=1,...,n$ or on a symbol itself,
while each symbol may be surrounded on both sides by two constants
from ${\cal A}_r$, moreover, $\bf H$ is associative apart from
${\cal A}_r$ with $r\ge 3$, where $\omega =\omega (\nu
):=s_1+...+s_n+u_1+...+u_n+p_1+...+p_n+y_1+...+y_n$. Then $l(\nu
_{j,k};\mu _{j,s})\le b(j+y)$, where $y=\max (\omega (\nu ),\omega
(\mu ))$. Since the series $\sum_{j=0}^{\infty } [b(j+y)]b^j$
converges, then $d(\mu ,\nu )<\infty $. On the other hand,
$\sup_{z\in U}|f(z,{\tilde z})-g(z,{\tilde z})|<\infty $ for bounded
functions $f$ and $g$ on $U$, consequently, ${\cal D}((f,\nu );
(g,\mu ))<\infty $. At the same time it is known that a continuous
function on a compact set is bounded, that finishes the proof of the
second statement of this lemma.
\par {\bf 17. Theorem.} {\it Let $U$ be an open region in ${\cal
A}_r$, $r\ge 2$. Then there exists a continuous operator $S:
C^0_{\cal P}(U,X)\times \Gamma \to X$, where $\Gamma $ is a space of
all rectifiable paths in $U$ such that $S((f,\nu );\gamma
):=\int_{\gamma }(f,\nu )(z,{\tilde z})dz$.}
\par {\bf Proof.} On a space of rectifiable paths there exists a
natural metric $d_r(\gamma ,\eta )$ induced by a metric between
arbitrary subsets $A$ and $B$ in ${\cal A}_r$: $d_r(A,B):=\max (\psi
(A,B), \psi (B,A))$, $\psi (A,B):= \sup_{z\in A}\inf_{\xi \in
B}|z-\xi | $. To a function $f$ and its phrase representation $\nu $
there corresponds a unique function $g$ and its phrase
representation $\mu $, which is constructed by the following. A
function $g$ is characterized by the conditions.
\par Let $f$ be defined by a continuous function
$\xi : U^2\to {\cal A}_r$ such that
\par $(i)$ $\xi (\mbox{ }_1z,\mbox{ }_2z)|_{\mbox{ }_1z=z,
\mbox{ }_2z=\tilde z}=f(z,{\tilde z})$, \\
where $\mbox{ }_1z$ and $\mbox{ }_2z\in U$. Let also $g: U^2\to
{\cal A}_r$ be a continuous function, which is $\mbox{
}_1z$-superdifferentiable such that
\par $(ii)$ $(\partial g(\mbox{ }_1z,\mbox{ }_2z)/\partial \mbox{ }_1z).1=
\xi (\mbox{ }_1z,\mbox{ }_2z)$ on $U^2$. Then put
\par $(iii)$ ${\hat f}(z,{\tilde z}).h:={\hat f}_z(z,{\tilde z}).h:=
[(\partial g(\mbox{ }_1z,\mbox{ }_2z)/
\partial \mbox{ }_1z).h]|_{\mbox{ }_1z=z,\mbox{ }_2z=\tilde z}$
for each $h\in {\cal A}_r$. Shortly it can be written as $(\partial
g(z,{\tilde z})/ \partial z).1=f(z,{\tilde z})$ and ${\hat
f}_z(z,{\tilde z}).h:={\hat f}(z).h :=(\partial g(z,{\tilde z})/
\partial z).h$.
\par A phrase $\mu $ is constructed by the algorithm: at first in each
word $\nu _{j,k}$ of the phrase $\nu $ substitute each $e$ one time
on $z$, that gives a sum of words $\lambda _{j+1,\alpha (k,i)}$,
$\alpha =\alpha (k,i)\in \bf N$, $j, k, i \in \bf N$, that is, $\nu
_{j,k}\mapsto \sum_i\lambda _{j+1,\alpha (k,i)}$, where $i$
enumerates positions of $e$ in the word $\nu _{j,k}$, $i=1,...,s$,
where $s=s_1$ (see \S 12). Therefore, to $\nu $ there corresponds
the phrase $\lambda =\sum_{j,\alpha }\lambda _{j,\alpha }$. The
partial differentiation $\partial f/\partial z$ was defined for
functions and their phrase reprsentations. Consider the remainder
$(\partial \lambda /\partial z).1 - \nu =:\zeta $. If $\zeta =0$,
then put $\mu =\lambda $. If $\zeta \ne 0$, then to each word $\zeta
_{j,k}$ of this phrase $\zeta $ apply a left or right algoritm from
\S 2.6 \cite{ludfov}. The first aforementioned step is common for
both algorithms. For this consider two $z$-locally analytic
functions $f_1$ and $q$ on $U$ such that $f_1$ and $q$ noncommute as
well as corresponding to them phrases $\psi $ and $\sigma $. Let
$\psi ^0:=\psi $, $\sigma ^0:=\sigma $, $\sigma ^{-n}:=\sigma
^{(n)}$, $(\partial (\sigma ^n)/\partial z).1 =: \sigma ^{n-1}$ and
$\sigma ^{-k-1}=0$ for some $k\in \bf N$ using the same notation
$\sigma _{j,k}^n$ for each word of the phrase $\sigma $. Then
\par $(iv)$ $(\psi \sigma )^1=\psi ^1\sigma -\psi ^2\sigma ^{-1}+
\psi ^3\sigma ^{-2}+...+(-1)^k\psi ^{k+1}\sigma ^{-k}$. In
particular, if $\psi =(a_1z^na_2)$, $\sigma =(b_1z^kb_2)$, with
$n>0$, $k>0$, $a_1, a_2, b_1, b_2\in {\cal A}_r\setminus {\bf R}I$,
then $\psi ^p=[(n+1)...(n+p)]^{-1}(a_1z^{n+p}a_2)$ for each $p\in
{\bf N}$, $\sigma ^{-s}=k(k-1)...(k-s+1)(b_1z^{k-s}b_2)$ for each
$s\in \bf N$. Apply this left algorithm by induction to appearing
neighbor subwords of a given word going from the left to right. Then
apply this to each word of a given phrase. Symmetrically the right
algorithm is:
\par $(v)$ $(\psi \sigma )^1=\psi \sigma ^1-\psi ^{-1}\sigma ^2+
\psi ^{-2}\sigma ^3+...+(-1)^n\psi ^{-n}\sigma ^{n+1}$, when $\psi
^{-n-1}=0$ for some $n\in \bf N$. Then apply this right algorithm
going from the right to the left by neighborhood subwords in a given
word.  For each word both these algorithms after final number of
iterations terminate, since a lenght of each word is finite. These
algorithms apply to solve the equation $(\partial \xi (z,{\tilde
z})/\partial z).1= \zeta (z,{\tilde z})$ for each $z\in U$ in terms
of phrases. Use one of these two algorithms to each word of the
phrase $\zeta $ that to get a unique phrase $\xi $ and then put $\mu
=\lambda +\xi $.
\par If $f_1$ and $q$ have series converging in $Int (B({\cal A}_r,0,R))$,
then these formulas demonstrate that there exists a $z$-analytic
function $(f_1q)^1$ with series converging in $Int (B({\cal
A}_r,z_0,R))$, since $\lim_{n\to \infty }(nR^n)^{1/n}=R$, where
$0<R<\infty $.  Since $f$ is locally analytic, then $g$ is also
locally analytic. Therefore, for each locally $z$-analytic function
$f$ and its phrase $\nu $ there exists the operator $\hat f$ and its
phrase $\partial \mu /\partial z$. Considering a function $G$ of
real variables corresponding to $g$ we get that in view of to Lemma
2.5.1 \cite{ludfov} all solutions $(g,\mu )$ differ on constants in
${\cal A}_m$, since $\partial g/\partial w_s+(\partial g/\partial
w_p).(s^*p)=0$ for each $s=i_{2j}$, $p=i_{2j+1}$,
$j=0,1,...,2^{r-1}-1$ and $\partial g/\partial w_1$ is unique, hence
$\hat f$ is unique for $f$.
\par Denote the described above mapping ${\cal P}\ni \nu \mapsto
\mu \in {\cal P}$ by $\phi (\nu )=\mu $. Thus each chosen algorithm
between two these algorithms gives $\phi (\nu _1+\nu _2)=\phi (\nu
_1)+\phi (\nu _2)=\mu =\mu _1+\mu _2$, if $\nu =\nu _1+\nu _2$,
moreover, $\omega (0)=0$. Therefore, this procedure gives the unique
branch of $S$. Recall, that if the following limit exists
$$(vi)\quad \int_{\gamma }(f,\nu )(z,{\tilde z})dz:=
\lim_P I((f,\nu );\gamma ;P),\mbox{ where}$$
$$(vii)\quad I((f,\nu );\gamma ;P):=\sum_{k=0}^{t-1}
({\hat f},{\hat {\nu }})(z_{k+1},{\tilde z}_{k+1}).(\Delta z_k),$$
where $\Delta z_k:=z_{k+1}-z_k$, $z_k:=\gamma (c_k)$ for each
$k=0,...,t$, then we say that $(f,\nu )$ is line integrable along
$\gamma $ by the variable $z$.
\par From Equation $2.7.(4)$ \cite{ludfov} it follows, that
$$(viii)\quad |S((f,\nu ) - (y,\psi );\gamma )|\le {\cal D}((f,\nu );
(y,\psi ))V(\gamma )C_1 \exp (C_2 R^{2^m+2})$$ for each $(f,\nu ),
(y,\psi )\in C^0_{\cal P}(U,X)$ for $U\subset {\cal A}_m$ with
finite $m\in \bf N$, where $C_1$ and $C_2$ are positive constants
independent of $R$, $(f,\nu )$ and $(y,\psi )$, $0<R<\infty $ is
such that $\gamma \subset B({\cal A}_m,z_0,R)$. The space $X_t$ is
complete relative to the norm of polyhomogeneous operators as well
as the Cayley-Dickson algebra ${\cal A}_r$ relative to its norm.
Therefore, $\{ S((f^v,\nu ^v);\gamma ): v\in {\bf N} \} $ is the
Cauchy sequence in $X$ for each Cauchy sequence $\{ (f^v, \nu ^v):
v\in {\bf N} \} \subset C^0_{\cal P}(U,X)$, hence $\{ S((f^v,\nu
^v);\gamma ): v\in {\bf N} \} $ converges.
\par For each rectifiable path $\gamma $ in $U$ there corresponds a
canonical closed bounded subset $V$ in ${\cal A}_r$ such that
$\gamma \subset Int (V)\subset V\subset U$. Then $V_m:=V\cap {\cal
A}_m$ is compact for each $m\in \bf N$ and a continuous function
$f|_V$ is uniformly continuous. Such that $\lim_{d(\eta ,0e)\to
0}{\cal D}((y|_{V_m},\eta ),(0,0e))=0$ for each $m\in \bf N$. For
finite $r$ take $m=r$. For infinite $r$ there exists a sequence of
rectifiable paths $\gamma _m$ in $V_m$ for suitable algebraic
emebeddings of ${\cal A}_m$ into ${\cal A}_r$ such that $\lim_{m\to
\infty } d_r(\gamma _m,\gamma )=0$. In the algebra ${\cal A}_r$ the
union of all algebraically embedded subalgebras ${\cal A}_m$, $m\in
\bf N$ is dense. Therefore, for each $\epsilon _k=1/k$ there exists
a continuous functions $f_m$ on $V_m$ and a phrase $\nu _m$ over
${\cal A}_m$ for $m=m(k)$ with $m(k)< m(k+1)$ for each $k\in \bf N$
such that $\lim_{m\to \infty }{\cal D}((f_m,\nu _m),(f|_{V_m},\mu
))=0$, where each $\nu _m$ has a natural extension over ${\cal A}_m$
and hence $f_m$ as the corresponding to this phrase function has the
natural extension on a neighborhood of $V_m$ in $U$. This gives
$S((f,\nu ),\gamma )=\lim_{m\to \infty } S((f_m,\nu _m),\gamma _m)$.
\par In accordance with the proof of Theorem 2.7 \cite{ludfov} the
operator $S$ is continuous for functions and analogously for
operator valued functions and their phrases relative to the metric
${\cal D}$.
\par {\bf 18. Theorem.} {\it Let a domain $U$ in ${\cal
A}_r$ satisfies the condition of Note 12 and $(f,\nu )$ be integral
holomorphic in $U$. Then $\partial S((f,\nu );\gamma )/\partial
z=({\hat f},{\hat {\nu }})(z)$ for each $z_0, z\in Int (U)$ and
rectifiable paths $\gamma $ in $Int (U)$ such that $\gamma (0)=z_0$
and $\gamma (1)=z$. In particular, if $f=(\partial g/\partial z).1$
and $\nu =(\partial \mu /\partial z).1$, where $\partial f/\partial
{\tilde z}=0$ and $\partial \nu /\partial {\tilde z}=0$ on $U$, then
$S((f,\nu );\gamma )=(g,\mu )(z)-(g,\mu )(z_0)$.}
\par {\bf Proof.} In view of the condition that $(f,\nu )$ is
integral holomorphic, that is, $S((f,\nu );\eta )=0$ for each
rectifiable loop $\eta $ in $U$, it follows, that $S((f,\nu );\gamma
)$ depends only on $z_0=\gamma (0)$ and $z=\gamma (1)$ and is
independent of rectifiable paths with the same ends. Since
$(\partial g(z,{\tilde z})/\partial z).1=f(z)$ and $(\partial \mu
(z,{\tilde z})/\partial z).1=\nu $, such that ${\hat f}(z,{\tilde
z})=\partial g(z,{\tilde z})/\partial z$ and ${\hat {\nu
}}(z,{\tilde z})=\partial \mu (z,{\tilde z})/\partial z$ exist in
the sense of distributions (see also \cite{luladfcdv}), then from
Formulas 17(vi,vii) the first statement of this theorem follows.
\par If $(f,\nu )=(\partial (g,\mu )/\partial z).1$, then
there exists $(y,\psi )\in C^0_{\cal P}(U,X)$ such that $(\partial
(y,\psi )/\partial z).1=(g,\mu )$, hence $(\partial ^2(y,\psi
)/\partial z^2)(h_1,h_2)=(\partial ^2(y,\psi )/\partial
z^2)(h_2,h_1)$ for each $h_1, h_2\in {\cal A}_r$, since $(f,\nu )$
and consequently $(g,\mu )$ and $(y,\psi )$ are holomorphic,
particularly, infinite continuous superdifferentiable in $z$ on $U$
(see Theorems 2.11, 2.16, 3.10 and Corollary 2.13 \cite{ludfov}).
Mention, that the first step of the algorithm of the proof of
Theorem 17 gives $\mu =\lambda $. For $h_1=1$ and $h_2=\Delta z$
this gives $({\hat f},{\hat {\nu }}).\Delta z=(\partial (g,\mu
)/\partial z).\Delta z$ and we can substitute the latter expression
in the integral sums instead of $({\hat f},{\hat {\nu }}).\Delta z$.
The second statement therefore follows from the fact that integral
holomorphicity of $(f,\nu )$ is equivalent to $\partial (f,\nu
)/\partial {\tilde z}=0$ on $U$, where $U$ satisfies conditions of
Remark 12.
\par {\bf 18.1. Corollary.} {\it Let an open domain $U$ in $\bf K$
satisfies conditions of Note 12 and $f$ be a function on $U$ such
that either $f'(\xi ).h=a(\xi )hb(\xi )$ for each $h\in {\bf H}=\bf
K$ or $f'(\xi ).h=c(\xi )(Y(\xi )h)$ for each $h\in {\bf O}=\bf K$
and each $\xi \in U$, where $Y$ is an ordered product of operators
$Y_{k,m}$ from Theorem 5, $a(\xi )$, $b(\xi )$, $c(\xi )$ and $Y(\xi
)$ are the holomorphic functions and the operator valued function
over $\bf H$ and $\bf O$ respectively such that $a(\xi )\ne 0$,
$b(\xi )\ne 0$, $c(\xi )\ne 0$ for each $\xi \in U$. Then the
function $f$ is pseudoconformal.}
\par {\bf Proof.} In view of the condition of this corollary
$\partial f(z)/\partial {\tilde z}=0$ for each $z\in U$. In
accordance with Theorem 18 $f(z)=\int_{z_0}^z(\partial f(\xi
)/\partial \xi ).1d\xi $ and this integral is independent of a
rectifiable path, but depends only  an initial $z_0$ and final $z$
points in $U$. Thus $f$ satisfies all conditions of Definition 1.
\par {\bf 19. Remark.} In the particular case of ${\bf C}={\cal A}_1$ we
simply have ${\hat f}=f$, but here the case of arbitrary $r\in \bf
N$ or $r=\Lambda $ with $card (\Lambda )\ge \aleph _0$ is
considered. For example, the pair of the function and its phrase
$((ab)(c{\tilde z}d),(aeb)(c{\tilde z}d))$ belong to $ C^0_{\cal
P}(U,{\cal A}_r)$ and there exists $(g,\mu )=((azb)(c{\tilde z}d),
(azb)(c{\tilde z}d))\in C^0_{\cal P}(U,{\cal A}_r)$ such that
$(\partial (g,\mu )/\partial {\tilde z}).1=(f,\nu )$, but $(f,\nu )$
is not integral holomorphic for each $r\ge 1$, where $a, b, c, d\in
{\cal A}_r$ are constants. In the commutative case of $\bf C$
conditions of section 14 simply reduce to the power series of
complex holomorphic functions and symbols $e$ and $\tilde e$ can be
omitted, but in the noncommutative case, when $r\ge 2$, the
conditions of section 14 are generally essential for Theorems 17 and
18.
\par {\bf 20. Theorem.} {\it Let $\cal F$ be a family of functions
quaternion or octonion holomorphic on a domain $U$ in ${\bf K}^n$
such that each $f\in \cal F$ for each $j=1,...,n$ at each $z\in U$
is either pseudoconformal by $\mbox{ }_jz$ or $\partial
f(z)/\partial \mbox{ }_jz=0$. If $\cal F$ is uniformly
$C^1$-continuous, then $\cal F$ is the normal family.}
\par {\bf Proof.} It is sufficient to prove, that the family
$\cal F$ is normal at each point $P$ of a ball $B({\bf K}^n,P_0,R)$
contained in $Int (U)$, where $0<R<\infty $. The uniform $C^1$
continuity means that, for each $\epsilon >0$ there exists $\delta
>0$ such that for each $z_1, z_2\in U$ with $|z_1-z_2|<\delta $
there is satisfied the inequality $|f(z_1)-f(z_2)|+\|
f'(z_1)-f'(z_2) \| <\epsilon $ for each $f\in \cal F$. The segment
$P_0P$ divide into $(k+1)$ segments of lengths not larger, than
$\delta $ such that $(k+1)\delta \ge R$. This gives points
$P_1,..,P_k$ in $[P_0,P]$.
\par Let $T$ be a topological space
and $C(T,c)$ denotes the space of all continuous mappings $f: T\to
F$ endowed with the compact open topology, where $F=\bf R^m$, $m\in
\bf N$. In view of the Ascoli Theorem (6.6.5) \cite{nari} if $H$ is
an equicontinuous subset of $C(T,c)$ such that $t'(H):= \{ h(t):
h\in H \} $ is a bounded subset in $F$ for each $t\in T$, then $H$
is relatively compact, that is, its closure is compact. On the other
hand, $C^n(B,{\bf K})$ is with the first axiom of countability, that
is, the topological character at each point is not greater, than
$\aleph _0$. In accordance with the Theorem 3.10.31 \cite{eng} the
sequential compactness and countable compactness are equivalent in
the class of sequential $T_1$-spaces, in particular, $T_1$-spaces
with the first countability axiom. The space $C^n(B,{\bf K})$ is
Lindel\"of, hence $H$ in it is relatively sequential compact by
Theorem 3.10.1 \cite{eng}.
\par If there exists $\sup_{f\in \cal F} |f(P_0)|<\infty $, then the
family $\cal F$ is normal on $B({\bf K}^n,P_0,R)=:B$, since ${\cal
F}\subset C^{\infty }$ and each embedding of $C^{n+1}$ into $C^n$ is
the compact operator, where $C^n=C^n(B,{\bf K})$ denotes the family
of all $n$ times continuously differentiable functions $g$ on $U$
with values in $\bf K$. Each $g\in C^n$ has the bounded norm $\| g
\| _{C^n(B,{\bf K})}<\infty $.  \par Thus, if there exists a bounded
at $P_0$ sequence $f_n(P_0)$ in $\cal F$, then $\{ f_n: n \} $ is
the normal sequence. If no any sequence $\{ f_n: n \} $ is bounded
at $P_0$, then there exists a subsequence $f_{n_k}$ with $\lim_{k\to
\infty }|f_{n_k}(P_0)|=\infty $, where $n_1<n_2<n_3<...$. In view of
$|f_{n_k}(P)|>|f_{n_k}(P_0)|-(k+1)\epsilon $ there exists
$\lim_{k\to \infty } f_{n_k}=\infty $ uniformly in $B({\bf
K^n},P_0,R)$. In any case $\cal F$ is normal in $B({\bf K^n},P_0,R)$
and hence in $U$.
\par {\bf 21. Theorem.} {\it Let $U$ be an open domain in $\bf K$.
If $f$ is a pseudoconformal function on $U$ and $f$ has a set $E$ of
zeros in $U$, then each zero is an isolated point of a set $E$.}
\par {\bf Proof.} Suppose on the contrary, that either
$E$ has a topological submanifold $Z$ of a dimension in the sense of
coverings greater than zero, or a zero $z_0\in U$ is a limit point
of $E$, where $E:=\{ z\in U: f(z)=0 \} $. In both cases there exists
a sequence $\{ z_n: z_n\in E\setminus \{ z_0 \}, n\in {\bf N} \} $
of pairwise distinct $z_n\ne z_m$ zeros for each $n\ne m\in \bf N$
such that $\lim_{n\to \infty }z_n=z_0\in E$, since $f$ is
continuous. The pseudoconformal function $f$ satisfies conditions
$\partial f(z)/\partial {\tilde z}=0$ and $(\partial f(z)/\partial
z)\ne 0$ for each $z\in U$, then
\par $(i)$ $df(z).h=(\partial f(z)/\partial z).h $ and
\par $(ii)$ $f(z+h)-f(z)=df(z).h+\alpha (z,h)$ for each $z\in U$
and $h\in \bf K$, where $\lim_{h\to 0}\alpha (z,h)/|h|=0$. Take
$h=h_n:=z_n-z_0$, then $h_n\ne 0$ and $|(\partial f(z)/\partial
z)|_{z=z_0}.h_n|>0$ for each $0<n\in \bf N$. On the other hand, for
each $v\in \bf K$ there is satisfied the equality \par $(iii)$ $Re
\{ [(\partial f(z)/\partial z).h_n] [(\partial f(z))/ \partial
z).v]^{*} \} |h_n| |v| $  \par $= |(\partial f(z)/\partial z).h_n|
|(\partial f(z)/\partial z).v| Re (h_n{\tilde v})$ \\
for $z=z_0$ in accordance with Definition 1. But $f(z_n)-f(z_0)=0$
for each $n$ and together with $(i,ii)$ this contradicts $(iii)$.
Therefore, $E$ consists only of isolated points.
\par {\bf 22. Definition.} A function $f$ on an open domain
$U$ in $\bf K$ we call $p$-pseudoconformal at a point $z_0$ in $U$,
where $p\in \bf N$, if there exists a pseudoconformal mapping $g$ at
$0$ such that $f(z_0+\zeta ^p)=g(\zeta )$ for each $\zeta $ in a
neighborhood of $0$. (For $p=1$ it was called pseudoconformal.)
\par {\bf 23. Theorem.} (Argument principle). {\it Let $f$ be a
holomorphic function on an open region $U$ in $\bf K$ satisfying
conditions of Remark 12 and let $f$ be $p$-pseudoconformal in a
neighborhood of each its zero in $U$, and let $\gamma $ be a closed
rectifiable curve contained in $U$. Then ${\hat I}n (0; f\circ
\gamma )=\sum_{\partial _f(a)\ne 0} {\hat I}n (a; \gamma )\partial
_f(a)$.}
\par {\bf Proof.} The function $Ln (z)$ is pseudoconformal
and $g$ is pseudoconformal on a neighborhood of each its zero. In
accordance with Theorem 21 each zero $z_0$ of $f$ is an isolated
point. Thus for $z_0$ there exists a neighborhood $V$ in $U$ such
that $Ln (f(z))$ is $p$-pseudoconformal on $V\setminus \{ 0 \} $ by
Theorem 6, since $Ln (g(z))$ is pseudoconformal. Each zero of a
pseudoconformal function is simple. In view of Conditions $(i-iii)$
of Definition 1 and Theorems 4 and 5 above, ${\hat I}n(0;f\circ
\gamma ) =\int_{\zeta \in \gamma }d Ln (f(\zeta ))=$ $\int_0^1d Ln
(f\circ \gamma (s))=$ $\int_{\gamma }f^{-1}(\zeta )df(\zeta )={\hat
I}n (a; \gamma )p$, where $a=z_0$ (see also the proof of Theorem
3.28 \cite{ludfov}). Using the homotopy Theorem 2.15 and Theorem
3.8.3 \cite{ludfov} we get the statement of this theorem.
\par {\bf 24. Theorem.} {\it Let $U$ be a proper open subset
of $\bf K$, let also $f$ and $g$ be two continuous functions from
${\bar U}:=cl (U)$ into $\bf K$ such that on a topological boundary
$Fr (U)$ of $U$ they satisfy the inequality $|f(z)|<|g(z)|$ for each
$z\in Fr (U)$. Suppose $f$ and $g$ are $\bf K$-meromorphic functions
in $U$ and $h$ be a unique continuous map from $\bar U$ into $\bf K$
such that $h|_E=f|_E +g|_E$, where $E:= \{ z: f(z)\ne \infty ,
g(z)\ne \infty \} $, $g$ is $p$-pseudoconformal in a neighborhood
$U_{z_0}$ in $U$ of each its zero $z_0$ and $1/g(z)$ is
$p$-pseudoconformal in $U_{z_0}\setminus \{ z_0 \} $ for each pole
$z_0$, where $p\in \bf N$ may depend on $z_0$. Then $\nu
_{g|_U}(0)-\nu _{g|_U}(\infty )=\nu _{h|_U}(0)-\nu _{h|_U}(\infty
)$.}
\par {\bf Proof.} If $z_0$ is a pole of $f$ at $z_0$, then
$1/f(z)$ has a zero at $z_0$. At each $z\in Fr (U)$ we have $g(z)\ne
0$ and $h(z)\ne 0$. Thus the argument principle 23 can be applied,
since $Ln (1/f)=-Ln (f)$. Moreover, $h=f+g=g(1+(1/g)f)$ and
$|(1/g)f|<1$ on $Fr (U)$. Thus the point $w=(1/g)f(z)$ is within the
unit ball $ \{ \zeta \in {\bf K}: |\zeta |<1 \} $. Therefore, the
vector $v=1+w$ can not rotate on $2\pi $ around zero. Thus the
winding numbers of $h$ and $g$ around the zero are the same. From
this the statement of this theorem follows.
\par {\bf 25. Lemma.} {\it Let $\{ f^k: k\in {\bf N} \} $ be a sequence
of functions quaternion or octonion holomorphic on a domain $U$ in
${\bf K}^n$ such that each $f^k$ for each $j=1,...,n$ at each $z\in
U$ is either pseudoconformal by $\mbox{ }_jz$ or $\partial
f(z)/\partial \mbox{ }_jz=0$. If $f^k$ converges to some function
$f$ relative to the metric $\rho =\rho _V$ for a closed bounded
subset $V=cl (Int (V))$ in $U$ and $f$ is not identically zero on
$V$, then each point $z\in Int (V)$ such that $f(z)=0$ is the limit
point of a family of all points $\zeta $ such that $f^k(\zeta )=0$
for some $k$.}
\par {\bf Proof.} Without restriction of the generality we can
suppose, that $0\in V$ and $f(0)=0$. In view of Lemma 11 a function
$f$ is holomorphic and at for each $j=1,...,n$ at each $z\in Int
(V)$ is either pseudoconformal by $\mbox{ }_jz$ or $\partial
f(z)/\partial \mbox{ }_jz=0$. Then $f$ has a power series
decomposition given by Formulas 13(ii,iii). Moreover,
\par $(i)$ $f(z) = \sum_{|m|=0}^{\infty } \phi _{f,|m|} (z-z_0)$, where
\par $(ii)$ $\phi _{f,c}(z-z_0) =
\sum_{k+|m'|=c} \{ (\phi _{k,m'}(\mbox{ }_2z,...,\mbox{ }_nz),
\mbox{ }_1z^k) \} _{q(2v)}$, \\
$c=0,1,2,...$, $\phi _{k,m'}$ are holomorphic functions by $\mbox{
}_2z,...,\mbox{ }_nz$, $m = (m_1,...,m_v)=(k,m')$, $k=m_1$,
$m'=(m_2,...,m_v)$, $m_j=0,1,2,...$, $|m'| := m_2+...+m_v$, $v\in
\bf N$. The series $(i,ii)$ converges uniformly in the domain $\{ z:
|\mbox{ }_1z|<R,$  $|(\mbox{ }_2z,...,\mbox{ }_nz)|<R' \} $, where
$0<R$ and $0<R'$. Evidently, all functions $\phi _{k,m'}$ can not be
identically zero on $Int (V)$, since $f$ is not identically zero on
$Int (V)$. Let $\phi _{k,m'}$ be the first such function, where they
are ordered by increasing $|m|=k+|m'|$ and for equal $|m|=c$
lexicographic ordering of vectors $m$ is used. This process of
decomposition can be extended by $\mbox{ }_2z,...,\mbox{ }_nz$.
\par In accordance with Theorem 3.17 \cite{ludfov}
each polynomial over ${\cal A}_r$ has roots in ${\cal A}_r$. In view
of Theorem 24 zeros of $f(\mbox{ }_1z,z')$ by $\mbox{ }_1z$ in the
ball $\{ |z|<\epsilon \} $ are limits of zeros of $f^k$, when $k$
tends to the infinity. A zero $\mbox{ }_1z_0$ of order $\partial
_f(\mbox{ }_1z_0,z')=p>0$ is the limit of zeros of $f^k(\mbox{
}_1z,z')$ by $\mbox{ }_1z$ of order $p$. Since $\epsilon
>0$, $R>0$ and $R'>0$ are arbitrary small, then the origin $0$ of
the coordinate system is the limit point  of points $z$ with
$f^k(z)=0$ and this lemma is proven in this case.
\par Let now $p=0$. For $|\mbox{ }_1z|<R'$ the polynomial
$\phi _{0,m'}(\mbox{ }_2z,...,\mbox{ }_nz)$ has $|m'|$ zeros
${z'}_0$ by $z'$ in $\bf K^{n-1}$, where $|z'|<R'$. These zeros are
the unique zeros of $f(\mbox{ }_1z_0,z')$ for $|z'|<R'$. In view of
Theorem 24 there exists $N\in \bf N$ such that $f^k(z)$ also has
$|m'|$ zeros in $ \{ |\mbox{ }_1z|<R, |z|<R' \} $ for each $k>N$.
Since $R>0$ and $R'>0$ are arbitrary small, then the lemma is proven
in this case also.
\par {\bf 26. Corollary.} {\it If a sequence of pseudoconformal functions
$f^k$ converges uniformly relative to the metric $\rho $ on the
domain $W:=\{ |\mbox{ }_1z|<R, |z'|<R' \} $ and if no any function
$f^k$ for each $k>N$ for some $N\in \bf N$ is not zero in $W$, then
the limit function $f$ is not zero in $W$ besides the case $f(z)=0$
for each $z\in W$.}
\par {\bf 27. Theorem.} {\it Let a domain $U$ in ${\bf K}^n$ satisfies
the condition of Remark 12 and a family $\cal F$ consists of
functions $\bf K$ holomorphic on $U$ such that each $f\in \cal F$
for each $\xi \in Int (U)$ by each variable $\mbox{ }_jz$ is either
pseudoconformal or $(\partial f(z)/\partial \mbox{ }_jz)_{\xi
=z}=0$. If $Int ({\bf K}\setminus \bigcup_{f\in \cal F} f(U))\ne
\emptyset $, then $\cal F$ is normal.}
\par {\bf Proof.} Denote by $a$ a centre of a ball $B:=B({\bf
K^n},a,\epsilon )$ such that $B\cap (\bigcup_{f\in \cal F}
f(U))=\emptyset $, then  $|f(z)-a|>\epsilon $ for each $f\in \cal
F$. Take an infinite sequence $f^k$ from $\cal F$. Compose new
sequence $\phi ^k(z) := [f^k(z)-a]^{-1}$, then $|\phi
^k(z)|<1/\epsilon $. Functions $f^k$ are bounded on $V$, then due to
Theorem 13 there exists a subsequence $\phi ^{n_k}$ which converges
uniformly relative to $\rho $ to a function $\Phi $ in a domain
$V=cl (Int (V))\subset U$, where $V$ is bounded. This function $\Phi
$ is holomorphic in $V$. No any $\phi ^k$ is zero in $V$. Thus $\Phi
$ can not be zero in $V$ besides the case given by Corollary 26,
when $\Phi =const $ in $V$. Therefore, consider two cases.
\par Let $\Phi (z)\ne 0$ in $V$. Therefore, there exists the
function $\Phi (z)=:[F(z)-a] ^{-1}$, hence $F(z)=a+(\Phi (z))^{-1}$
by each variable $\mbox{ }_jz$ is either pseudoconformal or
$(\partial F(z)/\partial \mbox{ }_jz)_{\xi =z}=0$. The subsequence
$f^{n_k}$ converges to $F(z)$ uniformly relative to the metric $\rho
$ on a canonical closed bounded subset $W$ in $Int (V)$. In fact
$F(z)-f^{n_k}(z)=[\Phi (z)]^{-1}-[\phi ^{n_k}(z)]^{-1}$ $=\phi
^{-1}[\phi ^{n_k}-\Phi ](\phi ^{n_k})^{-1}$, where $\Phi $ and $\phi
^{n_k}$ are not zero in $W$, when $|\Phi |>\lambda $ on $W$. Since
$\phi ^{n_k}$ converges to $\Phi $ relative to $\rho $, then there
exists $k_0\in \bf N$ such that $\rho _W(\phi ^{n_k},0)>\lambda
-\epsilon >0$ for each $k>k_0$. Thus $\rho _W (F,f^{n_k})<\epsilon
/[\lambda (\lambda -\epsilon )]$ for each $k>k_0$ and this implies
the uniform convergence of $f^{n_k}$ to $F$ on $W$.
\par In the case $\Phi (z)=0$ for each $z\in V$, for each $\epsilon
>0$ there exists $k_0\in \bf N$ such that $|\phi
^{n_k}|<\epsilon $ for each $k>k_0$. Since $f^{n_k}=a+1/\phi
^{n_k}$, then $|f^{n_k}|>1/\epsilon -|a|$. Then $f^{n_k}$ converges
uniformly to $\infty $ on $W$, where $W=cl (Int (W))\subset Int
(V)$. Thus $\cal F$ is the normal family.
\par {\bf 28. Theorem.} {\it Let $U$ be a domain in $\bf K$ satisfying
conditions of Remark 12 and $f(z)=(\partial g(z)/\partial z).1$ for
each $z\in U$, where $g$ is a holomorphic either pseudoconformal
function at $\xi $ or $g'(\xi )=0$ for each $\xi \in U$. If $\gamma
$ is a rectifiable loop in $U$, then for each $z$ encompassed by
$\gamma $ there is satisfied the inequality:
\par $(i)$ $|f(z)|\le \sup_{\xi \in \gamma }|f(\xi )|$. \par
 If $f$ is bounded on
$W$, where $W=cl (Int (W))$, $U\subset Int (W)$, then \par $(ii)$
$|f(z)|\le \sup_{\xi \in \partial U}|f(\xi )|$.}
\par {\bf Proof.} In view of Theorems 3.8.3 and 3.9 \cite{ludfov}
$f(z)=(2\pi )^{-1}(\int_{\gamma } f(\xi )dLn \xi )M^*$, where $M\in
\bf K$ is characterised by $\gamma $, $Re (M)=0$, $|M|=1$. On the
other hand, ${\hat f}(\xi )=(\partial g(\xi )/\partial \xi )$ for
each $\xi \in U$. If $g$ is pseudoconformal at $\xi $, then by
Theorems 4 and 5 $|(\partial g(\xi )/\partial \xi ).h|=|(\partial
g(\xi )/\partial \xi ).1| |h|$ for each $h\in \bf K$. If $g'(\xi
)=0$, then $|g'(\xi ).h|=0=|g'(\xi ).1| |h|$. Then from Formulas
17.(vi,vii) Inequality $(i)$ follows. Considering arbitrary
rectifiable loops in $U$ close to the boundary $\partial U$ of the
domain $U$ we get Inequality $(ii)$, since $f$ is continuous and
bounded on $U$.
\par {\bf 29. Remark.} The latter theorem is the quaternion and
octonion analog of the maximum principle known in complex analysis.
\par Henceforth, by a domain $U$ in $\bf K^n$ is undermined an open
connected subset in $\bf K^n$.
\par {\bf 30. Theorem.} {\it If a function $f$ on a domain
$U$ in $\bf K$ is pseudoconformal, then its image $V:=f(U)$ is also
a domain in $\bf K$.}
\par {\bf Proof.} Since $f$ is pseudoconformal on $U$, then
it is not constant, since $f'(z)\ne 0$ on $U$. \par Let $w_1$ and
$w_2$ be two arbitrary points in $V$, take $z_1\in f^{-1}(w_1)$ and
$z_2\in f^{-1}(w_2)$. Since $U$ is connected, then it is acrwise
connected and there exists a path $\eta : [\alpha ,\beta ]\to U$
such that $\eta (\alpha )=z_1$ and $\eta (\beta )=z_2$. In view of
the continuity of $\eta $ its image $\zeta :=f\circ \eta $ is the
path connecting points $w_1$ and $w_2$. Since $f(U)=V$, then $f(\eta
)\subset V$, hence $V$ is connected.
\par The set $U$ is open in $\bf K$, hence for each $z_0\in U$ there
exists $R>0$ such that $B({\bf K},z_0,R)\subset U$ and $B({\bf
K},z_0,R)\cap \partial U=\emptyset $. For each $w_0\in V$ there
exists $z_0\in U$ and $R>0$ such that $f(z_0)=w_0$ and $f(z)\ne w_0$
for each $z\in B({\bf K},z_0,R)\setminus \{ z_0 \} $ in accordance
with Theorem 21, since $f$ is the pseudoconformal function. Let
$\gamma := \partial B({\bf K},z_0,R) = \{ z\in {\bf K}: |z-z_0|=R \}
$. Put \par $(i)$ $\mu := \min_{z\in \gamma } |f(z)-w_0|$, \\
then $\mu >0$, since the continuous function $|f(z)-w_0|$ attains it
minimal value on compact $\gamma $, but for $\mu =0$ there would be
a point $w\ne w_0$, $w\in \gamma $ such that $f(z)=w_0$, but this
would be in the contradiction with the defintion of $B({\bf
K},z_0,R)$. \par Prove, that $ \{ w\in {\bf K}: |w-w_0|<\mu \}
\subset V$. Let $|w_1-w_0|<\mu $, then $f(z)-w_1=f(z)-w_0+(w_0-w_1)$
and $|f(z)-w_0|\ge \mu $ on $\gamma $ due to $(i)$. Since
$|w_0-w_1|<\mu $, then the function $f(z)-w_1$ has in $B({\bf
K},z_0,R)$ the same number of zeros as $f(z)-w_0$ in accordance with
Theorem 24, that is, at least one. Thus $w_1\in V$, hence $\{ w:
|w-w_0|< \mu \} \subset V$ and inevitably $V$ is open in $\bf K$.
\par {\bf 31. Definition.} Let $\hat {\bf K}$ be a one-point (Alexandroff)
compactification of $\bf K$ such that ${\hat {\bf K}}\setminus {\bf
K}=\infty $. A function $f: {\hat {\bf K}}\to {\hat {\bf K}}$ is
called pseudoconformal at infinity, if $g(z):=f(1/z)$ is
pseudoconformal at zero. A function $f: {\hat {\bf K}}\to {\hat {\bf
K}}$ is called a fractional $\bf R$-linear pseudoconformal mapping,
if it is a composition of mappings of the types: $(i)$ $z\mapsto
z+c=:S_c(z)$ for each $z\in {\hat {\bf K}}$ and a fixed $c\in \bf
K$; $(ii)$ $z\mapsto z^{-1}=:Inv (z)$ for each $z\in {\hat {\bf
K}}$; $(iii)$ $z\mapsto azb=:M_{a,b}(z)$ for each $z\in {\hat {\bf
H}}$ and fixed $a, b\in \bf H$ with $ab\ne 0$, or $z\mapsto Yz$ for
each $z\in {\hat {\bf O}}$, where $Y$ is the ordered product of
operators $Y_{k,m}$ from Theorem 5.
\par {\bf 32. Lemma.} {\it The family of all fractional $\bf
R$-linear pseudoconformal mappings form the group $F_L$ and each its
element is the pseudoconformal mapping.}
\par {\bf Proof.} Mention that each mapping of the types $(i-iii)$
is pseudoconformal due to Theorems 4 and 5, since each composition
of pseudoconformal mappings is pseudoconformal. The inverses of the
mappings $(i-iii)$ are: $S_{-a}\circ S_a(z) = id(z)$, $Inv \circ Inv
(z) = id(z)$, $M_{(a^{-1},b^{-1})}\circ M_{a,b}(z)=id (z)$,
$Y_{k,m}^{-1}=\exp (-t_{k,m}X_{k,m})$ for each $0\le k<m\le 7$ and
each given parameter $t_{k,m}\in \bf R$, where $id (z)=z$ for each
$z\in {\hat {\bf K}}$.
\par Compositions of diffeomorphisms are diffeomorphisms and compositions
of functions are associative. The family $F_L$ forms the group,
since each fractional $\bf R$-linear pseudoconformal mapping is
bijective and surjective and inverses of them are of the same type
due to Theorem 6. It remains only to prove, that the mapping
$z\mapsto z^{-1} =: Inv (z)$ is pseudoconformal. Since $z^{-1}z=1$,
then $((z^{-1})'.h)z +z^{-1}h=0$ for each $z\in \bf K$ and $h\in \bf
K$, hence $(z^{-1})'.h=-(z^{-1}h)z^{-1}$ due to the alternativity of
$\bf K$. Thus $Inv (z)$ is pseudoconformal at each $z\in \bf K$.
Since $Inv (z^{-1})=z$, then $Inv $ is pseudoconformal at the
infinity also.
\par {\bf 33. Theorem.} {\it Each fractional $\bf R$-linear
pseudoconformal mapping transforms each hypersphere in ${\hat {\bf
K}}$ into a hypersphere in ${\hat {\bf K}}$.}
\par {\bf Proof.} Each hypersphere in ${\hat {\bf K}}$ can be
written in the form \par $(i)$ $Ez{\tilde z} + J {\tilde z} +
z{\tilde J} +D=0$, where $E, D\in \bf R$, $J\in \bf K$. Indeed, the
usual hypersphere with a centre $z_0$ and radius $R>0$ is: $ \{ z\in
{\bf K}: |z-z_0|=R \} $, where $|z-z_0|^2=z{\tilde z}-z_0{\tilde
z}-z{\tilde z}_0+z_0{\tilde z}_0=R^2$ such that $E=1$, $J=-z_0$,
$D=|z_0|^2-R^2$. The inversion $Inv (z)$ gives instead of the
hypersphere with centre zero the hyperspehere with centre infinity.
If $E=0$, then this equation gives a hyperplane, which is a
hypersphere of the infinite radius. Each mapping $S_a$ simply shifts
the centre of the hypersphere. Each mapping $M_{a,b}$ or $Y$ rotates
${\hat {\bf K}} = {\hat {\bf H}}$ or ${\hat {\bf K}} = {\hat {\bf
O}}$ respectively, hence transforms a hypersphere into a
hyperspehere. It remains to verify, that the equation is of the same
type after application of the mapping $Inv (z)$. Let $w=z^{-1}$,
then $(i)$ gives $E|w|^{-2}+J{\tilde w}^{-1}+w^{-1}{\tilde J}+D=0$,
where $w^{-1}={\tilde w}|w|^{-2}$, since $w{\tilde w}=|w|^2$. After
multiplication of both its sides on $|w|^2$ we get: $E+Jw+{\tilde w}
{\tilde J}+Dw{\tilde w}=0$ which is of the same type as Equation
$(i)$.
\par {\bf 34. Definition.} Two points $z_1$ and $z_2$ are called
symmetrical relative to a hypersphere $S$ with the centre $z_0$ and
radius $R>0$ in $\bf K$, if $(i)$ $|z_1-z_0| |z_2-z_0|=R^2$ and
$(ii)$ $Arg (z_2-z_0)= Arg (z_1-z_0)$, where $z=|z|\exp (Arg (z))$
for each $z\in \bf K$, $Arg (z)\in \bf K$, $Re (Arg (z))=0$.
\par {\bf 35. Theorem.} {\it Each fractional $\bf
R$-linear pseudoconformal mapping $f$ transforms each two
symmetrical points $z_1$ and $z_2$ relative to a hypersphere $S$ in
$\bf K$ into two points $w_1=f(z_1)$ and $w_2=f(z_2)$ symmetrical
relative to the image hypersphere $f(S)$.}
\par {\bf Proof.} In view of Proposition 3.2 and Corollary 3.6
\cite{ludfov} each $z\in \bf K$ has the polar decomposition
$z=|z|\exp (Arg (z))$. In accordance with Theorem 33 the image
$f(S)$ is the hypersphere. Each mapping $S_a$, $M_{a,b}$ and $Y$
transforms symmetrical points into symmetrical, since $S_a$ only
shifts all points on $a$, while $M_{a,b}$ and $Y$ induce rotations
of the real shadow of $\bf K$ or rotates $\bf K$. But equations
$(i,ii)$ are invariant relative to such transformations. The mapping
$Inv \circ S_{z_0}(z)$ transforms $Arg (z-z_0)$ into $Arg
((z-z_0)^{-1})= - Arg (z-z_0)$. Equations $(i,ii)$ are equivalent to
the condition \par $(iii)$ $z_2-z_0 = R^2 ({\tilde z}_1-{\tilde
z}_0)^{-1}$. \\  Then $Inv \circ S_{z_0}(S)$ is the hypersphere with
centre infinity and radius $1/R$. If $w_1=Inv \circ S_{z_0}(z_1)$
and $w_2=Inv \circ S_{z_0}(z_2)$, then $(iii)$ takes the form $w_2 =
R^{-2}{\tilde w}_1^{-1}.$ Thus each $f\in F_L$ being the composition
of mappings of types $(i-iii)$ from Definition 31 transforms
symmetrical points $z_1$ and $z_2$ relative to $S$ into symmetrical
points $w_1$ and $w_2$ relative to $f(S)$.
\par {\bf 36. Lemma.} {\it Let two nonintersecting domains
$U_1$ and $U_2$ in $\bf K$ have a common piece $\gamma $ of their
borders which is a hyperplane over $\bf R$ in $\bf K$, let also
functions $f_1$ and $f_2$ be holomorphic in $U_1$ and $U_2$ and they
are continuous on the sets $U_1\cup \gamma $ and $U_2\cup \gamma $
respectively. Then, if $f_1(z)=f_2(z)$ for each $z\in \gamma $, then
the function $f(z):=f_1(z)$ for each $z\in U_1\cup \gamma $ and
$f(z)=f_2(z)$ for each $z\in U_2$ is holomorphic in $U:=U_1\cup
\gamma \cup U_2$.}
\par {\bf Proof.} From the condition of the lemma it follows, that
$f$ is continuous in $U$. The function $f$ is holomorphic if and
only if it is integral holomorphic in accordance with Theorems 2.11,
2.16 and 3.10 \cite{ludfov}. Let $\eta $ be a rectifiable loop in
$U$. If $\eta \subset Int (U_1)$ or $\eta \subset Int (U_2)$, then
$\int_{\eta } f(z)dz=0$. Consider the case $\eta \cap \gamma \ne
\emptyset $. \par Let $\gamma $ divides $\eta $ into two pieces
$\eta _1\subset U_1$ and $\eta _2\subset U_2$. Then $\int_{\eta
}f(z)dz= \int_{\eta _1}f(z)dz + \int_{\eta _2}f(z)dz +\int_{\eta
_3}f(z)dz - \int_{\eta _3}f(z)dz$, where $\eta _3\subset \gamma $,
$\eta _3$ has two ends as points of intersections of $\eta $ with
$\gamma $, $\eta _3$ is the rectifiable path parametrized by $t\in
[0,1]$. Let $\eta _{3,\delta }$ be a family of rectifiable paths
$\eta _{3,\delta } \subset U_1$ for each $0<\delta <1$ such that
$\eta _{3,\delta }$ intersects with $\eta _1$ in two points and
$\lim_{\delta \to 0}\eta _{3,\delta }(t)=\eta _3(t)$ uniformly by
$t\in [0,1]$.  For $-1<\delta <0$ consider also the family $\eta
_{3,\delta }(t)$ of rectifiable paths parametrized by $t\in [0,1]$
having two common points with $\eta _2$, such that $\eta _{3,\delta
}$ converges to $\eta _3$ uniformly by $t\in [0,1]$ and $\eta
_{3,\delta }\subset U_2$ for each $-1<\delta <0$, where $\eta
_{3,\delta }$ with $\delta >0$ and $\eta _{3,\delta }$ with $\delta
<0$ have opposite directions. Therefore, $\lim_{\delta \to 0}
\int_{\eta _{3,\delta }} f(z)dz=\int_{\eta _3}f(z)dz$, since $f$ is
continuous on $U$. \par Denote by $\eta _{\delta }$ the loop which
is the join of the path $\eta _{3,\delta }$ and a piece $\eta
_{1,\delta }$ of $\eta _1$ in $U_1$ for $\delta >0$ or a piece $\eta
_{2,\delta }$ of $\eta _2$ in $U_2$ for $\delta <0$ having common
ends with $\eta _{3,\delta }$. Then $\int_{\eta _{\delta }}
f(z)dz=0$ for each $-1<\delta <1$. Since $\lim_{\delta \to
0}\int_{\eta _{\delta }\cup \eta _{-\delta }}f(z)dz=\lim_{0<\delta
\to 0}\int_{\eta _{1,\delta }}f(z)dz + \lim_{0<\delta \to
0}\int_{\eta _{2,-\delta }} f(z)dz=\int_{\eta }f(z)dz$, then
$\int_{\eta }f(z)dz=0$ for each rectifiable loop $\eta $ in $U$,
consequently, $f$ is holomorphic on $U$.
\par {\bf 37. Definition.} Let $V$ be a surface of real dimension
$m$ in $\bf K^n$ such that $V$ is parametrized by $m$ real
parameters $\theta _1,...,\theta _m \in [0,1]$, $\gamma : [0,1]^m\to
\bf K^n$, $V = \gamma ([0,1]^m)$. We call $V$ the Jordan surface, if
the mapping $\gamma $ is bijective from $[0,1]^m$ onto $V$. A domain
$U$ in $\bf K^n$ we call a Jordan domain, if $U$ is an open domain
in $\bf K^n$ and its border $\partial U$ consists of a finite number
$k$ of pieces $V_j$ each of which is a Jordan hypersurface, that is,
of real dimension $n2^r-1$, where $r=2$ for ${\bf K}=\bf H$, $r=3$
for ${\bf K}=\bf O$, $n\in \bf N$, $dim_{\bf R}V_j\cap V_l\le
n2^r-2$ for each $j\ne l$, $\partial U=\bigcup_{j=1}^kV_j$.
\par {\bf 38. Theorem.} {\it Let $U$ and $W$ be two Jordan bounded
domains in $\bf K$, $\partial U=V$ and $\partial W=S$, where $U$ and
$W$ satisfy conditions of Remark 12. If a function $f: cl (U)\to cl
(W)$ is continuous, pseudoconformal on $U$ and $f: V\to S$ is
bijective, then $f$ is the pseudoconformal diffeomorphism of $U$ on
$W$.}
\par {\bf Proof.} Due to Definition 37 $U$ is open in $\bf K$.
Let $w_0\in W$ be an arbitrary point. Since $f(V)=S$, then $f\ne
w_0$ on $V$ and due to continuity of $f$ there exists a $\delta
$-band $V_{\delta } := \{ z\in cl (U): dist (z,V)<\delta \} $ for
some $\delta >0$ such that $f(V_{\delta })$ does not contain $w_0$,
where $dist (z,A) := \inf_{a\in A}|z-a|$ for a point $z\in \bf K$
and a subset $A\subset \bf K$. Take a branch of the function $Arg$
and a rectifiable loop $\eta $ in $V$. If $\eta $ is characterized
by $M\in \bf K$, $|M|=1$, $Re (M)=0$, where $\int_{\eta } dLn (z) =
2\pi M$, $M=M(\eta )$, then $N := (2\pi )^{-1} {\hat I}n(w_0;f\circ
\eta ) M^*$ is the integer number, hence $N$ is constant under
continuous deformation of the rectifiable loop $\eta $. \par The
function $f|_V: V\to S$ is the homeomorphism, then $N=1$ for a loop
$\eta $ winding once around $z_0$, since $f\circ \eta $ winds once
around $w_0$, where $f(z_0)=w_0$. For each rectifiable loop $\gamma
$ homotopic to $\eta $ relative to $V_{\delta }$ as in Theorem 3.9
\cite{ludfov} piecewise $\eta _j$ relative to the corresponding
subsets $V_{\delta ,j}$ in $V_{\delta }$ there is the equality
$N=(2\pi )^{-1} {\hat I}n (w_0; f\circ \gamma ) M^*$, where
$M=M(\gamma )$, $\bigcup_j\eta _j=\eta $, $\bigcup_jV_{\delta
,j}=V_{\delta }$. In view of Theorem 23 each equation $f(z_0)=w_0$
has precisely one solution.
\par The proof above shows, that $N := (2\pi
)^{-1} {\hat I}n(w_0;f\circ \eta ) M^*=0$ for each $w_0\notin cl
(W)$. In view of Theorem 30 $f|_U$ does not attain values from $S$,
since in this case it would take values from ${\bf K}\setminus cl
(W)$. Thus $f(cl (U))=cl (W)$, moreover, $f$ has each value $w_0\in
W$ once for some $z_0\in U$, $f(z_0)=w_0$. Therefore, $f: U\to W$ is
the homeomorphism. Since $f$ is pseudoconformal, then $f^{-1}$ is
pseudoconformal by Theorem 6, hence $f$ is the pseudoconformal
diffeomorphism of $U$ on $W$.
\par {\bf 39. Theorem.} {\it Let $U_1$ and $U_2$ be two Jordan
domains in $\bf K$, where $\partial U_j$ contains a domain $V_j$ in
a hyperplane or in a hypersphere over $\bf R$ in $\bf K$ for $j=1$
and $j=2$. Suppose that $W_1$ and $W_2$ are Jordan domains
symmetrical with $U_1$ and $U_2$ relative to $V_1$ and $V_2$
respectively, where $U_j\cap W_j=\emptyset $ for $j=1$ and $j=2$. If
a function $f$ maps $U_1$ on $U_2$ pseudoconformally, where
$f_1(V_1)=V_2$ and there exists the limit \par $(i)$ $\lim_{z\to
\zeta \in V_1}\partial f_1(z)/\partial z\ne 0$ for each $\zeta \in
\gamma $, then $f_1$ has a pseudoconformal extension $f$ on $W_1$
and $f$ maps $U_1\cup V_1\cup W_1$ on $U_2\cup V_2\cup W_2$.}
\par {\bf Proof.} Consider at first the case, when $V_1$ and $V_2$
are domains in hyperplanes contained in ${\bf R}i_0\oplus {\bf R}i_1
\oplus ... \oplus {\bf R}i_{2^r-2}$. In the domain $W_1$ symmetrical
with $U_1$ relative to $V_1$ define the function $f_2(z):= \theta
_{2^r-1}[f_1(\theta _{2^r-1}(z))]$, where $\theta
_{2^r-1}(z)=z_0i_0+...+z_{2^r-2}i_{2^r-2}-z_{2^r-1}i_{2^r-1}$ for
each $z=z_0i_0+...+z_{2^r-1}i_{2^r-1}$. For $z\in W_1$ the point
$\theta _{2^r-1}(z)$ belongs to $U_1$. Since the function $f_1$ is
pseudoconformal by $z$ in $U_1$, then the function $\theta
_{2^r-1}\circ f_2(z) = f_1(\theta _{2^r-1}(z))$ is pseudoconformal
by $\theta _{2^r-1}(z)$. Use the projection operators $\pi
_j(h):=h_j=(-hi_j+ i_j(2^r-2)^{-1} \{ -h
+\sum_{l=1}^{2^r-1}i_l(hi_l^*) \} )/2$ for each $j=1,2,...,2^r-1$, \\
$\pi _0(h):=h_0=(h+ (2^r-2)^{-1} \{ -h
+\sum_{l=1}^{2^r-1}i_l(hi_l^*) \} )/2$ on ${\cal A}_r$, where $r=2$
for quaternions, $r=3$ for octonions. Therefore, there is the
expression of $\theta _{2^r-1}(z)$ with the help of generators of
$\bf K$ and $z$ such that in this representation $\partial \theta
_{2^r-1}(z)/\partial z=0$ and inevitably $f_2(z)$ is holomorphic by
$z$ in $W_1$. \par By Theorem 38 $f_1$ is continuous on $cl (U_1)$.
By the conditions of this theorem $f_1(V_1)=V_2$, where it was
supposed above, that $V_2$ is contained in the hyperplane $\{ z\in
{\bf K}: z_{2^r-1}=0 \} $. Thus if $z\in W_1$ tends to $x\in V_1$,
then $\theta _{2^r-1}(z)$ tends to $x$ and $f_2(z)$ converges to
$\theta _{2^r-1}f_1(x)=f_1(x)$, hence $f_1(x)=f_2(x)$ for each $x\in
V_1$.
\par Since $U_1\cap W_2=\emptyset $, then $f_1$ and $f_2$ satisfy
the conditions of Lemma 36, consequently, the function $f(z)$ such
that $f(z)=f_1(z)$ in $U_1\cup V_1$ and $f(z)=f_2(z)$ in $W_1$ is
holomorphic in $U_1\cup V_1\cup W_1$. If $f'|_{V_1}\ne 0$, then from
Condition 2 of Definition 1 for each $z\in U_1$ it follows, that
this condition is also true for each $z\in V_1$ due to Condition
$(i)$ of this theorem. By the construction $f$ is the
pseudoconformal function of $U_1\cup V_1\cup W_1$ on $U_2\cup
V_2\cup W_2$ (see Theorems 4 and 5), since if $A, B\in O(m)$ and
$det (A)<0$, $det (B)<0$, then $AB\in SO(m)$, particularly, for
$m=2^r$.
\par The general case reduces to the one considered above with the help of
fractional $\bf R$-linear pseudoconformal mappings $p_1$, $p_2$
transforming domains $V_1$, $V_2$ in hyperspheres into domains in
hyperplanes. Then the function $g_1:=p_2\circ f_1\circ p_1^{-1}$ due
to the proof given above has the pseudoconformal extension on the
domain $p_1(W_1)$ symmetrical with $p_1(U_1)$ relative to $p_1(V_1)$
such that the extension $g_2$ maps $p_1(W_1)$ on $p_2(W_2)$. Thus
the function $f_2:=p_2^{-1}\circ g_2\circ p_1$ is the
pseudoconformal mapping of $W_1$ on $W_2$ extending $f_1$.
\par {\bf 40. Remark.} A domain $V_2$ in a hypersphere can contain
the infinity point, then the extension $f$ will be meromorphic and
pseudoconfromal everywhere on $U_1\cup V_1\cup W_1$ besides a point
$z_0\in V_1$ corresponding to $\infty \in V_2$, where $f$ has a
simple pole, since $f$ is pseudoconformal and univalent on $U_1$.
\par {\bf 41. Definition.} Let $V$ be an open domain in $\bf K$ and
$g: V\to \bf K$ be a holomorphic (or pseudoconformal) diffeomorphism
on $g(V)$. Suppose that $\pi $ is a hyperplane or a hypersphere over
$\bf R$ embedded into $\bf K$ such that $V\cap \pi =: p$ is open in
$\pi $. Then $g(p)=:S$ we call a holomorphic (or pseudoconformal)
hypersurface in $\bf K$. If $\pi $ is a real hypersphere embedded
into $\bf K$ and $p=V\cap \pi $, then $g(p)=:S$ we call the closed
holomorphic (or pseudoconformal respectively) hypersurface.
\par {\bf 42. Theorem.} {\it  If a boundary $\partial U$ of an open
domain $U$ in $\bf K$ is such that it contains a pseudoconformal
hypersurface $S$ and if $f: U\to \bf K$ is a pseudoconformal
function such that $\lim_{z\to \zeta \in S} f'(z)\ne 0$ for each
$\zeta \in S$, then $f$ has a pseudoconformal extension through
$S$.}
\par {\bf Proof.} For each $t_0\in p$ there exists a neighborhood
$B := B({\bf K},t_0,\rho )$ contained in $V$ such that $g(t)$ on it
is pseudoconformal in accordance with Definition 41. The function
$g$ maps $p$ on $S$. Let $B^+$ denotes one half of $B$ contained in
$B\setminus p$ such that $g$ maps $B^+$ into $U$. Then the function
$w=f\circ g$ in $B^+$ satisfies the conditions of the symmetry
principle, consequently, has a pseudoconformal extension in $B$.
Thus $f(z)$ has the pseudoconformal extension through $S$.
\par {\bf 43. Corollary.} {\it Let $U$ and $W$ be two domains in $\bf K$
with pseudoconformal hypersurfaces $\partial U$ and $\partial W$
(may be both closed) and if $f: U\to \bf K$ is a pseudoconformal
function such that $\lim_{z\to \zeta \in \partial U} f'(z)\ne 0$ for
each $\zeta \in \partial U$, then $f$ has a pseudoconformal
extension in $cl (U)$.}
\par {\bf 44. Definitions.} Let $f_0$
be a locally $z$-analytic function on a subset $U$ in $\bf O$. If
there exists a subset $W$ in $\bf O$ such that $U\subset W$ and
$U\ne W$ and a locally $z$-analytic function $f$ on $W$ such that
$f|_U=f_0$, then we say, that $f_0$ has an analytic extension.
\par By a canonical element we mean a pair $F := (B,f)$, where
$B=B({\bf O},z_0,r^{-}) := \{ z\in {\bf O}: |z-z_0|<r \} $ is an
open maximal ball for which $f$ has a power series expansion with
centre in $z_0$ converging on $B$, $f(z)=\sum_k \{ (b_k, (z-z_0)^k)
\} _{q(2|k|)} $, $k=(k_1,...,k_m)$, $0\le k_j\in \bf Z$ for each
$j$, $|k|=k_1+...+k_m$, $b_k=(b_{k,1}, ...b_{k,m(k)})$, $b_{k,j}\in
\bf O$ for each $j$, $m=m(k)\in \bf N$, $\{ (b_k,(z-z_0)^k) \}
_{q(2|k|)} := \{ b_{k,1}(z-z_0)^{k_1}...b_{k,m(k)}(z-z_0)^{k_{m(k)}}
\} _{q(2|k|)}$. The point $z_0$ is called the centre of $F$ and $B$
is called its ball, $f(z)$ is called the value of $F$ at $z\in B$,
$m(k) := \max \{ j: k_j\ne 0, k_i=0 \mbox{ for each } i>j \} $.
\par Two elements $F=(B,f)$ and $G=(C,g)$ are immediate analytic
extensions of each other, if $W:=B\cap C\ne \emptyset $ and
$f|_W=g|_W$. Two elements are analytic extensions of each other, if
there exists a finite chain of elements $F_j=(B_j,f_j)$,
$j=1,...,n\in \bf N$, such that $F_1$ is the immediate extension of
$F$,..., $F_{j+1}$ is the immediate extension of $F_j$,...,$F_n$ is
the immediate extension of $G$.
\par We say, that a canonical element $F_0=(B_0,f_0)$ is extendible
along a path $\psi : [0,1]\to \bf O$, $z_0=\psi (0)$, if there
exists a family of elements $F_t=(B_t,f_t)$ with centres $z_t:=\psi
(t)$ and nonzero radii $r_t$, such that if there exists an open
segment $(\alpha ,\beta )$ such that $t_0\in (\alpha ,\beta )\subset
[0,1]$ and $\psi (\alpha ,\beta )\subset B_{t_0}$, then for each
$t\in (\alpha ,\beta )$ the element $F_t$ is the immediate extension
of $F_{t_0}$.
\par {\bf 45. Theorem.} {\it Let $\gamma _0$ and $\gamma _1$
be two homotopic paths in a $7$-connected domain $U$ in $\bf O$ and
an analytic element $F$ has an analytic extension along each path
$\gamma _p$, $p\in [0,1]$, where $\{ \gamma _p: p\in [0,1] \} $ is a
family of paths realizing the homotopy. Then all results of
extensions of $F$ along $\gamma _0$ and $\gamma _1$ coincide.}
\par {\bf Proof.} Let $B$ be a ball, $B=B({\bf O},z_0,r^{-})$,
on which $f(z)=\sum_k(b_k,(z-z_0)^k)$ is the convergent series. Then
write it in the form $f(z)=\sum_nP_n(z-z_0)$, where $P_n(\lambda
(z-z_0))=\lambda ^nP_n(z-z_0)$ for each $\lambda \in \bf R$, that
is, $P_n(z-z_0)=\sum_{\eta (k)=n}\{ (b_k,(z-z_0)^k) \} _{q(2|k|)}$,
$\eta (k):=k_1 +...+ k_{m(k)}$. Since $f(z)$ is infinite
differentiable inside $B$, then there exists $\partial
_z^nf(z)|_{z=z_0}.(h_1,...,h_n)=
\partial _z^nP_n(z)|_{z=z_0}.(h_1,...,h_n)$ for each $h_1,...,h_n\in
\bf O$. Since $(z-z_0) = (z-z_1) + (z_1-z_0)$ and $(z-z_0)^n=
(z-z_1)^n + (z-z_1)^{n-1}(z-z_0) + ((z-z_1)^{n-2}(z-z_0))(z-z_1)$
$+...+ (z-z_0)(z-z_1)^{n-1}$ $+(z-z_1)^{n-2}(z-z_0)^2 +
(((z-z_1)^{n-3}(z-z_0))(z-z_1))(z-z_0)$  $+... +
(z-z_0)^{n-1}(z-z_1)+ (z-z_0)^n$, then a series in $(z-z_0)$ can be
decomposed in the $(z-z_1)$-variable in a corresponding ball for
$z_1\in B$. Suppose that $f(z)$ is equal to zero in the ball $B({\bf
O},z_1,c^{-})\subset B({\bf O},z_0,r^{-})$ for some $0<c\le r$ and
$z_1\in B({\bf O},z_0,r^{-})$. In view of arbitrariness of
$h_1,...,h_n$ we get from $\partial _z^nf(z)|_{z=z_0}=0$, that
$P_n=0$ for each $n$. Thus $f(z)=0$ on the entire $B$.
\par $(ii)$. Let two elements $F_1$ and $G_1$ be obtained by extensions
of $F$ along $\gamma $ with the help of families $F_t$ and $G_t$
respectively, where $F_0=G_0$, $t\in [0,1]$. The set $\Psi := \{
t\in [0,1]: F_t=G_t \} $ is nonvoid, since $0\in \Psi $. Moreover,
$\Psi $ is the open subset in $[0,1]$, since for each $t_0\in \Psi $
there exists a neighborhood $U$ of $t_0$ such that $U\subset [0,1]$
and for each $t\in U$ a point $\gamma (t)$ belongs to the common
ball of convergence of elements $F_{t_0}$ and $G_{t_0}$. In view of
Definitions 44 $F_t$ and $G_t$ are immediate analytic extensions of
$F_{t_0}=G_{t_0}$ for each $t\in U$. Thus $U\subset \Psi $.
\par Suppose, that $t_0$ is a limit point of $\Psi $, and let
$U$ be a neighborhood of $t_0$ such that for each $t\in U$ points
$\gamma (t)$ belong to a less ball $B$ of convergence of $F_{t_0}$
and $G_{t_0}$. There exists $t_1\in U\cap \Psi $ such that
$F_{t_1}=G_{t_1}$ in it. Since $\gamma (t_1)\in B$, then $F_{t_0}$
and $G_{t_0}$ are immediate analytic extensions of equal elements
$F_{t_1}$ and $G_{t_1}$. Let $B_1$ the ball of convergence of
$F_{t_1}=G_{t_1}$, then $f_{t_0}=g_{t_0}$ on $B\cap B_1$. There
exists a ball contained in $B\cap B_1$. But from the first part
$(i)$ of the proof it follows, that $f_{t_0}=g_{t_0}$ on $B$, hence
$F_{t_0}=G_{t_0}$ and inevitably $t_0\in \Psi $. Thus $\Psi $ is
closed. Therefore, if an element $F_0$ has an analytic extension
along a path $\gamma $, then one gets a definite element, which is
independent of a choice of a family $\{ F_t: t\in (0,1] \} $
realizing an analytic extension. \par The set $\Psi $ is closed and
open in $[0,1]$, $\Psi \ne \emptyset $, hence $\Psi =[0,1]$ and
inevitably $F_1=G_1$.
\par {\bf 46. Theorem.} {\it Let $U$ be an open subset in $\bf O$
satisfying conditions of Remark 12. If an analytic element $F_0$ has
an analytic extension along each path in $U$ with the initial point
at the centre of $F_0$, then a result of an analytic extension along
paths does not depend on paths and is uniquely characterized by
their ends.}
\par {\bf Proof.} In the domain $U$ each two
paths with common ends are homotopic. Therefore, this theorem
follows from Theorem 45.
\par {\bf 47. Theorem.} {\it Let $U$ be an open
subset in $\bf O$ satisfying conditions of Remark 12 and with a
boundary $\partial U$ consisting more, than one point. Then $U$ is
pseudoconformally equivalent with the open unit ball in $\bf O$.}
\par {\bf Proof.} A family $\{ f_j: j\in \Lambda \} $ of functions
$f_j: U\to \bf O$ for a domain $U$ in $\bf O$ we call locally
uniformly bounded, if for each proper compact subset $V$ in $U$
there exists a constant $C>0$ such that $|f_j(z)|\le C$ for each
$j\in \Lambda $ and $z\in V$, where $\Lambda $ is a set. A family
$\{ f_j: j\in \Lambda \} $ we call locally uniformly continuous, if
for each $\epsilon >0$ and each compact proper subset $V$ in $U$
there exists $\delta =\delta (\epsilon ,V)$ such that
$|f_j(z)-f_j(y)|<\epsilon $ for each $j\in \Lambda $ and each $z,
y\in V$ with $|z-y|<\delta $. Mention, that each compact proper
subset $K$ in $U$ can be encompassed by a rectifiable loop contained
in $U$. Then from Formulas $(2.12)$ and Theorem $3.9$
\cite{ludoyst,ludfov} it follows, that if a family $\{ f_j: j\in
\Lambda \} $ of functions octonion holomorphic in a region $U$ is
locally uniformly bounded, then it is locally uniformly continuous,
where $\Lambda $ is a set.
\par Take two different points $\alpha $ and $\beta $ in the boundary
$\partial U$. Apart from the complex case the function $((z - \alpha
) (z - \beta )^{-1})^{1/2}$ while octonion holomorphic extension in
a domain $U$ has an infinite family of branches, but it is possible
to choose two of them which are different in their signs. We denote
them by $\phi _1$ and $\phi _2$. Since $(z - \alpha ) (z - \beta
)^{-1} = 1 + (\beta - \alpha ) (z - \beta )^{-1}$, then $(z - \alpha
) (z - \beta )^{-1}$ is univalent (that is, bijective), hence each
branch $\phi _1$ and $\phi _2$ is univalent. Suppose that $z_1,
z_2\in U$ and $\phi _1(z_1) = \phi _2(z_2)$, hence $\phi _1^2(z_1) =
\phi _2^2(z_2)$, consequently, $z_1=z_2$, then $\phi _1(z_1)=-\phi
_2(z_2)$ leads to the contardiction. Thus, $\phi _1(U)\cap \phi
_2(U) = \emptyset $.
\par In view of Theorem 6 a composition of pseudoconformal mappings is
pseudoconformal. Functions $\phi _1$ and $\phi _2$ are
pseudoconformal due to Corollary 3, Theorem 5 and Corollary 7. Then
the function $f_1(z) := r [\phi _1(z) - z_0]^{-1}$ is univalent and
conformal in $U$. Moreover, $|f_1 (z)|\le 1$ for each $z\in U$. The
argument principle over $\bf O$ (see Theorem 23 above) is satisfied
for pseudoconformal functions, hence if a sequence of functions
$f_n$ pseudoconformal in the domain $U$ is uniformly in each compact
subset $K$ in $U$ converging to a function $f\ne const $, then from
$f(z_0)=0$ for $z_0\in U$ it follows, that for each ball $B({\bf
O},z_0,\rho ^{-})\subset U$ there exists $n_0\in \bf N$, such that
for each $n>n_0$ there exists $z_n\in B({\bf O},z_0,\rho ^{-})$ with
$f_n(z_n)=0$.
\par Consider the family $\cal S$ of all univalent conformal functions
$f$ on $U$ with $\sup_{z\in U}|f(z)|\le 1$. Since $f_1\in \cal S$,
then ${\cal S}\ne \emptyset $. Then the family ${\cal S}_1$ of
functions from $\cal S$ with $\| f'(z_1) \| \ge \| {f_1}' (z_1)\|
>0$ at some point $z_1\in U$ is compact, since $\cal S$ is compact,
hence ${\cal S}_1$ is complete. Therefore, there exists a function
$f_0\in {\cal S}_1$ such that $\| f'(z_1) \| \le \| {f_0}'(z_1) \| $
for each $f\in {\cal S}_1$. Thus, $f_0(U)\subset B({\bf O},0,1)$.
\par Show that $f_0(z_1)=0$. Otherwise, in ${\cal S}_1$ would be
a function \par $g(z)=[f_0(z)-f_0(z_1)] [1-(f_0(z_1))^*f_0(z)]^{-1}$
with \par $g'(z).h=f_0'(z).h[1-(f_0(z_1))^*f_0(z)]^{-1} -
[f_0(z)-f_0(z_1)] (([1-(f_0(z_1))^*f_0(z)]^{-1}$
\par $(f_0(z_1))^*f_0'(z).h)
[1-(f_0(z_1))^*f_0(z)]^{-1})$ \\
 for each $h\in \bf O$. Therefore,
$|g'(z_1)| = |f_0'(z_1)| [1-|f_0(z_1)|^2]^{-1}> |f_0'(z_1)|$ despite
the extremal property of $f_0$.
\par Prove, that $f_0(U)=B$. Suppose contrary, that there exists
$b\in B\setminus f(U)$. Since $f_0(z_1)=0$, then $b\ne 0$. Then
$|{\tilde b}^{-1}|>1$, hence ${\tilde b}^{-1}\notin f(U)$. In
accordance with Theorem 46 there exists a univalent branch of the
root $w(z) = [(f_0(z)-b) (1-{\tilde b}f_0(z))^{-1}]^{1/2}$ which
belongs to $\cal S$. Therefore, $v(z):=[w(z)-w(z_1)] [1-
(w(z_1))^*w(z)]^{-1}$ also belongs to $\cal S$ by Theorem 6. Choose
a branch of $Ln (z)$ corresponding to a branch of $z^{1/2}$ such
that $z^{1/2}=\exp (2^{-1}Ln (z))$. Then $(\partial z^{1/2}/\partial
z).h=2^{-1}z^{-1/2}h$ for each $h\in \bf R$. The choosing $h_1\in
\bf O$, $h_1\ne 0$, such that $f_0'(z_1).h_1\in \bf R$ and the using
$|f'(z).t|=|f'(z).1| |t|$ for each pseudoconformal mapping $f$ at
$z$ and each $t\in \bf O$ due to Theorem 5 gives that
$|v'(z_1)|=(1+|b|)|f_0'(z_1)||b|^{-1/2}2^{-1}$. On the other hand,
$1+|b|>2|b|^{1/2}$, since $|b|<1$, consequently, $h\in {\cal S}_1$
and $|v'(z_1)|>|f_0'(z_1)|$ despite the extremal property of $f_0$.
\par {\bf 48. Example. The modular function.} Consider an open
domain $Q_0$ contained in the unit ball $B({\bf K},0,1) := \{ z\in
{\bf K}: |z|\le 1 \} $ such that $Q_0$ has a boundary $\partial Q_0$
consisting of pseudoconformal hypersurfaces joined together in a way
that $Q_0$ has $5$ vertices $A_1,...,A_5$ if ${\bf K}=\bf H$ or $9$
vertices $A_1,...,A_9$ if ${\bf K}=\bf O$, that is, $Q_0$ is the
twisted simplex. Suppose that each hypersurface of $\partial Q_0$ is
a piece of a hypersurface $\partial B({\bf K},z_j,1)$ such that
hyperspheres $\partial B({\bf K},0,1)$ and $\partial B({\bf
K},z_j,1)$ are orthogonal relative to the scalar product $(x,y) :=
Re (x{\tilde y})$. This means that in each plane $\pi $ of real
dimension two containing $z_0$ and $z_j$ and an intersection point
$w\in \partial B({\bf K},z_j,1)\cap \partial B({\bf K},0,1)$ the
corresponding circles $\pi \cap \partial B({\bf K},z_j,1)$ and $\pi
\cap \partial B({\bf K},z_0,1)$ are orthogonal. In view of Theorem
2.1.5.7 \cite{lusmldg05} and Theorem 47 above there exists a
pseudoconformal mapping $\mu $ of $Q_0$ on a semispace ${\bf K}^+ :=
\{ z\in {\bf K}: Re (i_{2^r-1}{\tilde z})>0 \} $ such that $\mu $
has a continuous extension on $\partial Q_0$ with $\mu (\partial
Q_0)=P$, $\mu (A_1)=0$, $\mu (A_2)=1$, $\mu (A_3) = \infty $, $\mu
(A_4)=i_1$, ...,$\mu (A_{2^r+1})=i_{2^r-2}$, where $P = \partial
{\bf K}^+$. \par Consider open subdomains $Q_1^k$ in $B({\bf
K},0,1)$ symmetrical with $Q_0$ relative to its faces,
$k=1,...,2^r+1$, since each face $\partial Q_{0,k}$ is characterized
by $2^r$ vertices. In accordance with Theorem 42 there exists a
pseudoconformal extension of $\mu $ into $Q_1^k$ for each $k$. The
images $A_l^k$ of $A_l$ symmetrical relative to the face $\partial
Q_{0,k}$ belong to $\partial B({\bf K},0,1)$ also for each $k$.
Since the inversion mapping or more generally a $\bf R$-fractional
pseudoconformal function gives that each face of $Q_1^k$ is a domain
of a hypersphere also (see Theorem 35). The extended function $\mu $
maps each $Q_1^k$ onto ${\bf K}^- := \{ z\in {\bf K}: Re
(i_{2^r-1}{\tilde z}) <0 \} $. Then consider subdomains $Q_2^{k,p}$
symmetrical with $Q_1^k$ relative to its faces, $p=1,...,2^r$. Then
extend $\mu $ from $Q_1^k$ onto $Q_2^{k,p}$ for each $p$ and each
$k$. \par Then by the mathematical induction continue this process.
As the result we construct the function $\mu $ pseudoconformal in
$Int (B({\bf K},0,1))$. Call this $\mu $ the modular function over
$\bf K$. The modular function is pseudoconformal in $Int (B({\bf
K},0,1))$ due to its construction. The interior of the unit ball
$Int (B({\bf K},0,1))$ is its domain of holomorphy, since on
$\partial B({\bf K},0,1)$ the set $S_t$ of all points obtained from
$A_t$ by consequtive symmetrical reflections relative to faces and
which are vertices of domains $Q_s^{k_1,...,k_s}$ is dense, where
$k_1,...,k_s = 1,2,..., 2^r$, $t=1,...,2^r+1$. If a sequence $ \{
z_n: n\in {\bf N} \} $ in $Int (B({\bf K},0,1))$ tends to $A_t$
while $n$ tends to the infinity, then $\mu (z_n)$ tends to $a_t$,
where $a_1=0$, $a_2=1$, $a_3=\infty $,
$a_4=i_1$,...,$a_{2^r+1}=i_{2^r-2}$. Thus $\mu $ can not be extended
onto $B({\bf K},0,1)$ even continuously.
\par Consider now the inverse function $\mu ^{-1}$.
Choose its branch in ${\bf K}^+$ and mapping it into $Q_0$. By the
symmetry principle (see Theorem 39 above) this branch has the
extension into ${\bf K}^-$. Extensions may be through different
pseudoconformal hypersurfaces, hence branches may be different. This
process of extension may be done infinitely. Thus $\mu ^{-1}$ has
infinite countable number of values such that points $0, 1, \infty ,
i_1,...,i_{2^r-2}$ are logarithmic branching points, in another
words branching points of infinite order. Moreover, all values of
$\mu ^{-1}$ are in the unit ball.
\par {\bf 49. Theorem.} {\it Let $U$ be a domain in $\bf K^n$
satisfying conditions of Remark 12. Suppose that $\cal F$ is a
family of functions $f(z)=g'(z).1$ holomorphic on $U$ such that both
each $f$ and the corresponding to it $g$ satisfy the condition:
$(i)$ for each $j=1,...,n$ and at each $z\in U$ a function is either
pseudoconformal by the variable $\mbox{ }_jz$ or its partial
derivative by $\mbox{ }_jz$ is zero. If each $f$ does not take no
any value from the set $ \{ 0, 1, i_1, i_2,..., i_{2^r-2} \} $, then
they form a normal family $\cal F$ in $U$, where $r=2$ for ${\bf
K}=\bf H$ and $r=3$ for ${\bf K}={\bf O}$.}
\par {\bf Proof.} It is sufficient to prove this for each point
$P_0$ in $Int (B)$, where $B$ is a closed ball $B=B({\bf
K^n},z_0,R)$ contained in $U$. Let $\nu (z)=\mu ^{-1}(z)$ be a
branch of the inverse function to the modular function (see \S 48).
To each $f\in \cal F$ there corresponds $\phi (z) := \nu (f(z))$,
where $\nu (f(z))$ and $\nu (g(z))$ are pseudoconformal for each
$f\in \cal F$ and the corresponding to it $g(z)$ due to Theorem 6,
since the only critical points of $\nu $ are $\{ 0, 1, \infty ,
i_1,..., i_{2^r-2} \} $. Let $\pi _s$ denotes the projection on
${\bf R}i_s$ such that $\pi _s(z)=z_s$, where $z=z_0+z_1i_1+...
z_{2^r-1}i_{2^r-1}$, $z\in \bf K$, $z_0,...,z_{2^r-1}\in \bf R$.
Thus the values of $\phi $ do not cover the semispace $\{ z\in {\bf
K}: \pi _s(z)>0 \} $ for each $s=0,...,2^r-1$. Therefore, by Theorem
27 the family $\{ \phi : \phi =\nu (f), f\in {\cal F} \} $ is
normal.
\par Consider a sequence $ \{ f_n: n\in {\bf N} \} $ and its values
at the point $P_0$. Let at first $ \{ f_n : n \} $ has at least one
limit point $\alpha $ different from $ \{ 0, 1, \infty , i_1,...,
i_{2^r-2} \} $. Thus there exists a subsequence $f_{n_1} (P_0)$,
$f_{n_2} (P_0)$,... converging to $\alpha $. Therefore, the sequence
$\phi _{n_1}(P_0)$, $\phi _{n_2}(P_0)$,... converges to $\nu (\alpha
)$, since $\alpha \notin \{ 0, 1, \infty , i_1,..., i_{2^r-2} \} $
and $\nu $ is regular at $\alpha $. The limit function $\Phi $ is
different from $\infty $ and satisifies Condition $(i)$ in
accordance with Inequality 2.7(4) and Theorem 3.10 and Remarks 3.34
\cite{ludfov}, Theorems 17 and 18 and Lemma 11 above. At each point
of the hypersphere $\partial B$ they satisfy the inequality $\pi _0
(\phi _{n_k}(z))>0$. But $\phi _{n_k}$ does not take any real value
on $\partial B$. Therefore, $\Phi $ could not take any real value at
any point $z\in \partial B$ besides the case, when $\Phi $ is a real
constant, but this is impossible, consequently, $\pi _0 (\Phi
(z))>0$ on $\partial B$.
\par When $P_0$ is arbitrary in $Int (B)$, then the point $\Phi
(P_0)$ is outside the real axis on definite distance as well as from
values $ \{ i_1,..., i_{2^r-2} \} $. Thus the subsequence
$f_{n_k}=\mu (\phi _{n_k})$ converges to the function $F=\mu (\Phi
)$ satisfying Condition $(i)$ uniformly on $Int (B)$.
\par Suppose now that the unique limit points of $f_n(P_0)$ may be
$\alpha \in  \{ 0, 1, \infty , i_1,..., i_{2^r-2} \} $. Then there
exists a subsequence $f_{n_1}$,...,$f_{n_k}$,... converging to
$\alpha $, for example, if $\alpha \ne 0$ and $\alpha \ne \infty $.
Therefore, it can be supposed that no any $f_{n_k}$ takes value zero
on $Int (B)$. In this case $Ln f_{n_k}$ satisfies Condition $(i)$.
Choose such branch of $Ln$ that $Ln f_{n_k}(P_0)=:b$ satisfies:
$|b+{\tilde b}|<2\pi $, hence $g_{n_k}(z):=[2\pi i_s+Ln
f_{n_k}(z)](4\pi i_s)^{-1}$ satisfies Condition $(i)$, where $s=1$
for $\alpha =1$, or $s=p$ for $\alpha =i_p$. Then $g_{n_k}(P_0)$
does not take values $\{ 0, 1, \infty , i_1,..., i_{2^r-2} \} $,
since $Ln i_s=\pi i_s/2$. Therefore, there exists a subsequence
$g_{n_1},...,g_{n_k},...$ converging uniformly relative to the
metric $\rho $ from Definitions 8 on $Int (B)$ to the function $G$
satisfying Condition $(i)$ due to Inequality 2.7(4) and Theorem 3.10
and Remarks 3.34 \cite{ludfov}, Theorems 17 and 18 and Lemma 11
above, since in the case of a necessity $B$ can be shrinked a little
and $B$ is compact. We have $G(P_0)=[2\pi i_1+ \pi i_s/2][4\pi
i_1]^{-1}=:c$, but no any $g_{n_k}(P_0)$ is equal to this value.
Therefore, $G(z)=c$ on $Int (B)$. From $f_{n_k}=\exp ((2g_{n_k} -1)
2\pi i_s)=\exp (g_{n_k}4\pi i_s)$ it follows, that the subsequence
$f_{n_k}$ converges uniformly on $Int (B)$ to $1$ if $\alpha =1$, or
to $i_p$ for $p>0$. If $\alpha =0$, then consider $1-f_n$ instead of
$f_n$. If $\alpha = \infty $, then consider $1/f_n$ instead of
$f_n$. Hence there exists a subsequence $f_{n_k}$ converging
uniformly to $\infty $ on $Int (B)$ in the latter case.
\par {\bf 50. Theorem.} {\it If there is
given a family $\cal F$ of functions holomorphic in a neighborhood
$U$ of a point $z=0$ in $\bf K^n$, $2\le n\in \bf N$, such that
$\cal F$ is normal at each point $z=(z',\mbox{ }_{n-k+1}z,...,\mbox{
}_nz)$ with $z'=0\in \bf K^{n-k}$, $z':=(\mbox{ }_1z,...,\mbox{
}_{n-k}z)$, $\mbox{ }_1z,...,\mbox{ }_nz\in \bf K$, $0<|(\mbox{
}_{n-k+1}z,...,\mbox{ }_nz)|<R$, $1\le k\le n-1$, but $\cal F$ is
not normal at $z=0$, then for each $\epsilon >0$ there exists
$\delta >0$ such that for each $z'$ with $0<|z'|<\delta $ there
exists $(\mbox{ }_{n-k+1}z,...,\mbox{ }_nz)$ with $|(\mbox{
}_{n-k+1}z,...,\mbox{ }_nz)|<\epsilon $ such that the family $\cal
F$ is not normal at $z=(z',\mbox{ }_{n-k+1}z,...,\mbox{ }_nz)\in
U$.}
\par {\bf Proof.} There exists $R'$ such that $0<R'<R$ and
$\cal F$ is normal at each point \\  $(0,\mbox{
}_{n-k+1}z,...,\mbox{ }_nz)\in U$ with $0<|(\mbox{
}_{n-k+1}z,...,\mbox{ }_nz)|<R'$. The boundary $\partial B({\bf
K^k},0,R')$ can be covered by a finite number of balls $B({\bf
K^k},y,R_1)$, where $y\in
\partial B({\bf K^k},0,R')$, $0<R_1\le R'$. The family $\cal F$ is
normal at each point $z=(x,y)$, $x=z'$, $y=(\mbox{
}_{n-k+1}z,...,\mbox{ }_nz)$, such that $|x|\le b$, $\eta _1\le
|y|\le \eta _2$ for some constants $0<b<\infty $, $0<\eta _1<\eta
_2<\infty $. Therefore, there exists a subsequence $ \{ f_n(z): n \}
\subset \cal F$ which converges relative to the metric $\rho $ to
the limit function $f$ or to $\infty $ for $|x|\le b$ and $\eta
_1\le |y|\le \eta _2$, in particular, for $|x|=b$ and $|y|=\eta _2$.
In view of Inequality 2.7(4) and Theorem 3.10 and Remarks 3.34
\cite{ludfov} this convergence is uniform on some ring $0<c\le
|x|\le b$, $0<\eta _1\le |y_2|\le \eta _2$. In view of Theorem 2.12
\cite{luladfcdv} and Formula 3.34(i) \cite{ludfov} the sequence
$f_n$ converges relative to the metric $\rho $ to a holomorphic
function $F$ satisfying Condition $49(i)$ on $\{ (x,y): |x|\le b,
|y|\le \eta _2 \} $, since $c>0$ and $\eta _1>0$ are arbitrary
small.
\par If this happens for each infinite convergent sequence from
$\cal F$, then $\cal F$ would be normal at $z=0$ also. Therefore,
there exists a sequence $f_1, f_2,...$ bounded at $0$, but
converging to $\infty $ in $ \{ z=(x,y): |x|\le b, |y|\le \eta _2 \}
$. That is, there exists $A=const >0$ such that $|f_n(0)|<A$ for
each $n$, while for chosen $0<\eta _1 <\eta _2$ and for each $C>0$
there exists $N=N(C)>0$ such that $|f_n(x,y)|>C$ for each $n>N$ and
each $z=(x,y)$ with $|x|\le b$ and $\eta _1\le |y|\le \eta _2$. In
particular, $|f_n(0,y)|>A$ for each $n>N(A)$ and each $z=(x,y)$ with
$|x|\le b$ and $\eta _1\le |y|\le \eta _2$. The function
$|f_n(0,y)|$ is continuous, $|f_n(0,0)|<A$ and $|f_n(0,y)|>A$ for
$\eta _1<|y|<\eta _2$. Hence the minimum of $|f_n(0,y)|$ is attained
for some $y_0$ with $|y_0|<\eta _1$. In view of Theorem 24 there are
zeros $z=(0,y_1)$ of $f_n(0,y)$ in $\{ y: |y|<\eta _1 \} $ and
inevitably there is $y_1=y_0$. Choose $y_1=:y_{n,0}$ with the
minimum $|y_{n,0}|$. \par Consider a limit point $Y_0$ of the
bounded sequence $ \{ y_{n,0}: n \} $. Suppose that $Y_0\ne 0$, then
without restriction of the generality $0<\eta _1\le |Y_0|\le \eta
_2$, since $\eta _1<R'<\eta _2<R$. In this ring $f_n(0,y)$ converges
to the infinity uniformly, on the other hand, $f_n(0,y_{n,0})=0$,
that leads to the contradiction. Thus there exists $\lim_{n\to
\infty }y_{n,0}=0$, hence $Y_0=0$. \par Without restriction of the
generality we can suppose that for some $\eta _1<t<\eta _2$ each
$f_n(0,y)$ has zeros in $B({\bf K},0,t)$. Consider an impilic
holomorphic function $y_n=y_n(x)$ characterized by the equation
$f_n(x,y(x))=0$ such that $y_n(0)=y_{n,0}$. There are possible only
two cases. 1. For each $|x|\le b$ and each $n$ with $n>n_0$ there is
satisfied the inequality $|y_{n,0}(x)|\le t$. We can choose a
subsequence and suppose that $n_0=1$. 2. For each $b_p>0$ with
$0<b_p\le b$ there exist $n=n_p$ and $x=x_p$ with $|x|\le b_p$ such
that $|y_{n_p}(x)|>t$. Prove that the second option is impossible.
For this consider a sequence $0<...<b_{p+1}<b_p<...$ with
$\lim_{p\to \infty }b_p=0$. There is an accumulation point $Y$ of a
sequence $y_{n_p}(x_p)$, hence $\lim_{p\to \infty
}(x_p,y_{n_p})=(0,Y)$, where $t\le |Y|\le \eta _2$. At $(0,Y)$ the
family $\cal F$ is normal and $|f_n|$ tends to the infinity while
$n$ tends to the infinity. This contradicts the equation
$f_{n_p}(x_p,y_{n_p})=0$. Thus the second case is impossible. It
remains the first case. \par There exists $b>0$ such that the
implicit function $y_n(x)$ defined by $f_n(x,y_n(x))=0$ with
$y_n(0)=y_{n,0}$ satisfies the inequality $|y_n(x)|\le t$ for each
$n$. Choose a subsequence $y_{n_p}$ such that $\lim_{p\to \infty
}y_{n_p}(x_0)=Y_0$ and $|Y_0|\le t<\eta _2$ for some $x_0$ with
$|x_0|<b$. Join $Y_0$ with $Y_1$ such that $t<|Y_1|<\eta _2$ with
the segment $[Y_0,Y_1]$. The sequence $f_n$ is normal at $(x_0,Y_1)$
and $\lim_{n\to \infty } |f_n(x_0,Y_1)|=\infty $. Therefore, there
exists $\delta >0$ such that $\lim_{n\to \infty }f_n(x,Y)=\infty $
for each $x, Y$ with $|x-x_0|<\delta $ and $|Y-Y_0|<\delta $. But
this contradicts the results obtained above, since there exists a
subsequence $f_{n_p}$ such that $\lim_{p\to \infty
}f_{n_p}(x_0,y_{n_p}(x_0))=0$ and there exists $N_0$ with
$|y_{n_p}(x_0)-Y_0|<\delta $ for each $n_p>N_0$.
\par Thus the hypothesis 1 is impossible. Therefore, there exists
$y_0\in (Y_0,Y_1)$ such that $\{ f_n: n \} $ and hence $\cal F$ also
are not normal at the point $(x_0,y_0)$. Thus for each $\eta
_1<t<\eta _2$ there exists $b>0$ such that to each $x_0\in B({\bf
K^{n-1}},0,b)$ there corresponds at least one $y_0$ such that the
family $\cal F$ is not normal at $(x_0,y_0)$. This finishes the
proof of Theorem 50.
\par {\bf 51. Note.} The point $(x,y)=(0,0)$ is the limit point of
$\sigma _n := \{ (x,y): f(x,y)=0 \} $, where $f\in \cal F$, that is,
there exists an infinite sequence $f_n$ such that for each $\epsilon
>0$ there exists $N$ for which $\sigma _n\cap B({\bf K^n},0,\epsilon
)\ne \emptyset $ for each $n>N$.
\par {\bf 52. Corollary.} {\it The family $E$ of all points $(x,y)$,
where $\cal F$ is not normal is dense in itself.}
\par {\bf Proof.} If $(x,y)$ is an accumulation point of $E$,
then $\cal F$ is not normal at $(x,y)$, hence $(x,y)\in E$. In view
of the preceding theorem each $(x,y)\in E$ is an accumulation point
of $E$.
\par {\bf 53. Remark and Definitions.} \par Remind that a point
$(x_0,y_0)$ of a $\bf K$-valued function $\psi (x,y)$ is called
regular if $\psi $ is holomorphic in a neighborhood $V$ of $(x,y)$.
If $\psi $ is holomorphic in $V\setminus \{ (x_0,y_0) \} $, but not
at $(x_0,y_0)$, then $(x_0,y_0)$ is called a singular point. \par A
characteristic surface by the definition is given as a family of all
points $(x,y)$ such that $\psi (x,y)=0$ on a domain of definition of
$\psi $, where $\psi $ is a holomorphic function, $x\in \bf
K^{n-1}$, $y\in \bf K$, $2\le n\in \bf N$. \par A point $(x_0,y_0)$
is called an ordinary or regular point if $y-y_0=\phi (x-x_0)$,
where there exists $b>0$ such that $\phi $ is the holomorphic
function on $B({\bf K^{n-1}},x_0,b)$.
\par A point $(x_0,y_0)$ is called critical if $y-y_0=\phi
((x-x_0)^{1/p})$ for some $p$ such that $1<p\in \bf N$.
\par A point $(x_0,y_0)$ is called singular of a charateristic surface,
if it is a singular point of $\psi (x-x_0,y-y_0)$ and a limit point
of regular points, where $\psi (x,y)=0$.
\par Henceforth we consider the case, when a surface $S$ is given by
$\psi (x,y)$ and at each point $(x,y)$ of $S$ the tangent space
$T_{(x,y)}S$ has the real dimension equal to $(n-1)2^r$, if
something otherwise is not stated. Suppose, that $\psi $ satisfies
conditions of the implicit function theorem at each regular point
$(x_0,y_0)$, in particular, the operator $\partial \psi
(x,y)/\partial y$ is invertible at $(x_0,y_0)$ such that $y-y_0=\phi
(x-x_0)$ in a neighborhood $Y_{x_0}\times Z_{y_0}$, where $Y_{x_0}$
and $Z_{y_0}$ are neighborhoods of $x_0$ and $y_0$ respectively.
Then $S$ is called the characteristic surface.
\par There is the following generalization of Theorem 50.
\par {\bf 54. Theorem.} {\it Let $\cal F$ be a family of functions
$f$ satisfying Condition 49(i) in $U\setminus  \{ (\xi ,\zeta ) \}
$, where $U$ is open in $\bf K^n$. Suppose, that there exists a
family of charateristic surfaces $S_{\alpha }$ given by $\psi
(x,y,\alpha )=0$ regular in a neighborhood of $(\xi ,\zeta )$ such
that $S_{\alpha }$ contains $(\xi ,\zeta )$ for each $\alpha \in
B({\bf K^{n-1}},0,b_0)$ for some $b_0>0$, where $\psi (x,y,\alpha )$
is pseudoconformal by $\alpha $ and holomorphic by $x$ and $y$. If
$\cal F$ is normal at each point in $U\setminus \{ (\xi ,\zeta ) \}
$ besides $(\xi ,\zeta )$, then on each $S_{\alpha }$ there exists
at least one point in a neighborhood $V\setminus \{ (\xi ,\zeta ) \}
$ at which $\cal F$ is not normal.}
\par {\bf Proof.} Since $\psi (x, y, \alpha )$ is holomorphic by
$\alpha $, then each point $(x,y)$ belongs to one $S_{\alpha }$ in a
neighborhood of $(\xi ,\zeta )$. In this case the family of
characteristics can be locally written as $\psi (x,y)=\alpha $. Then
consider the transformation: $x_1 = \psi (x,y) - \alpha _0$, $y_1 =
y-\zeta $. Therefore, $x=A(x_1,y_1)$ and $y=B(x_1,y_1)=y_1+ \zeta $,
where $A$ and $B$ are holomorphic functions due to conditions of
this theorem and \S 53 above. To the characteristic $\psi = \alpha $
there corresponds the plane $x_1 = \alpha -\alpha _0$. Also to each
plane $x=\beta $, where $\beta $ is in a neighborhood of zero, there
corresponds the characteristic: $\psi =\alpha _0 + \beta $. \par  To
each $f\in \cal F$ there corresponds the function $F(x_1,y_1)$
satisfying Condition 49(i) in a neighborhood of $(0,0)$. This gives
the family $\cal G$ such that $\cal F$ is normal at $(x,y)$ if and
only if $\cal G$ is normal at $(x_1,y_1)$. This gives the bijective
correspondense between the sets $E$ and $E_1$ in neighborhoods of
$(\xi ,\zeta )$ and $(0,0)$, where $\cal F$ and $\cal G$
respectively are not normal. Thus this theorem follows from Theorem
50.
\par {\bf 55. Theorem.} {\it  Let $\sigma $ be a characteristic surface
given by $\psi (x,y)=0$ such that $\sigma $ is regular at a point
$(\xi ,\zeta )$. Suppose that $\pi $ is a characteristic plane
containing $(\xi ,\zeta )$, but $\pi $ is not tangent to $\sigma $.
Then a characteristic obtained from $\sigma $ by the way of
translation into an arbitrary point $P$ in a neighborhood $V$ of
$(\xi ,\zeta )$ form a family analytic and regular in $V$.}
\par {\bf Proof.} A plane $\pi $ can be written by the equation:
$<x-\xi ,p>+(y-\zeta )=0$, where $x, p\in {\bf K^{n-1}}$, $<x,\eta
>:=\sum_{j=1}^{n-1}\mbox{ }_jx\mbox{ }_j{\tilde {\eta }}$. Then the
family of surfaces is: $\psi (x-\alpha ,y-<\alpha ,p>)=0$, where
$(\xi ,\zeta )$ is translated into $x=\xi +\alpha $ and $y= \zeta +
<\alpha ,p>$. Since $\alpha $ is arbitrary, then $\alpha $ can be
substituted on $-\alpha $ such that $\psi (x,y,\alpha )=\psi
(x+\alpha ,y+<\alpha ,p>)=0$, where $\alpha \in \bf K^{n-1}$.
\par If $\alpha =0$, then \par $(i)$
$\psi (\xi ,\zeta )=\psi (\xi ,\zeta ,0)=0$, also
\par $(ii)$ $(\partial \psi (x,y,\alpha )/\partial \alpha ).h=
(\partial \psi (x+\alpha ,y+<\alpha ,p>)/\partial x).h + (\partial
\psi (x+\alpha ,y+<\alpha ,p>)/\partial y).<h,p>$ and $(\partial
\psi (x,y,0 )/\partial \alpha ).h= (\partial \psi (\xi ,\zeta
)/\partial x).h+(\partial \psi (\xi ,\zeta )/\partial y).<h,p>$ for
each $h\in \bf K^{n-1}$.
\par The plane $\pi $ is not tangent to $\sigma $ at $(\xi ,\zeta
)$. Therefore, $(\partial \psi (\xi ,\zeta )/\partial x).h
+(\partial \psi (\xi ,\zeta )/\partial y).<h,p>\ne 0$ as the
function of $h\ne 0$, hence $(\partial \psi (\xi ,\zeta ,0)/\partial
\alpha )\ne 0$. In view of Remark and Definitions 51
\par $(iii)$ $\psi $ is regular in a neighborhood of $(\xi ,\zeta
)$,  hence $\psi $ is holomorphic by $x$ and $y$ variables. Thus
$\psi (x,y,\alpha )$ is pseudoconformal by $\alpha $ at $(\xi ,\zeta
,0)$. On $\sigma _{\alpha }$ there is the condition $d\psi
(x,y,\alpha )=0$, that is, $(\partial \psi (x,y,\alpha )/\partial
x).dx + (\partial \psi (x,y,\alpha )/\partial y).dy=0$,
consequently, \par $(iv)$ $(\partial y(x)/\partial x).h= - (\partial
\psi (x,y,\alpha )/\partial y)^{-1}. ((\partial \psi (x,y,\alpha
)/\partial x).h)$
\\ for each $h\in \bf K^{n-1}$ due to the alternativity of $\bf K$.
Since the composition of holomorphic functions is holomorphic, then
the corresponding $\phi (x,\alpha )$ is holomorphic by $x$ and
pseudoconformal by $\alpha $ (see \S 53). Since $\phi '(x-x_0)$ is
of rank $2^r$ at $0$, then $\partial \phi (x-x_0,\alpha )/\partial
x$ is of maximal rank at $(0,\alpha )$.
\par {\bf 56. Theorem.} {\it Let $\cal F$ be a family of functions
$f: U\to \bf K$ such that each $f$ satisfies Condition 49(i), where
$U$ is a domain in $\bf K^n$ satifying conditions of Remark 12.
Suppose that $E$ is a subset in $U$, where $\cal F$ is not normal.
For a marked point $z_0\in U$ it is impossible, that $E$ contains a
point $P_0$ such that the function $|z_0-z|$ by $P(z)\in E$ has a
local maximum at $P_0$.}
\par {\bf Proof.} In view of Theorem 50 $E$ is dense in itself.
At first prove that there exists a regular characteristic $\sigma $
tangent to $B:=B({\bf K^n},z_0,|OP_0|)$ at $P_0$ such that $\sigma $
has not any other common point with $B$ in a neighborhood of $P_0$,
where $O$ is the point with coordinates $z_0$. Consider the equation
$<x-\xi ,a>=y-\zeta $ of this characteristic $\sigma $, $x\in \bf
K^{n-1}$, $y\in \bf K$. Thus $y_p=\zeta _p+\sum_{j=1}^{n-1}
\sum_{i_v{\tilde i}_s=i_p}(\mbox{ }_jx_v-\mbox{ }_j\xi _v)\mbox{
}_j{\tilde a}_s$ for each $p=0,...,2^r-1$, where
$y=y_0i_0+...+y_{2^r-1}i_{2^r-1}$, $y_p\in \bf R$, $v,
s=0,1,...,2^r-1$, $r=2$ for ${\bf K}=\bf H$ and $r=3$ for ${\bf
K}=\bf O$. Then $|z_0-z|^2 =: F(x) = \sum_{j,v}(\mbox{ }_j
x_v)^2+\sum_vy_v^2$, where $P=P(z)$, consequently, \par
$F(x)=\sum_{j,v}(\mbox{ }_j x_v)^2 + \sum_p [\zeta _p + \sum_{v,
i_v{\tilde i}_s=i_p}(\mbox{ }_jx_v-\mbox{ }_j\xi _v)\mbox{
}_j{\tilde a}_s]^2$. Therefore,
\par $(\partial F/\partial \mbox{ }_jx_v)/2=\mbox{ }_jx_v
+ \sum_p[\zeta _p + \sum_j\sum_{i_w {\tilde i}_t = i_p; i_v{\tilde
i}_s=i_p}(\mbox{ }_jx_w-\mbox{ }_j\xi _w)\mbox{ }_j{\tilde a}_t]
\mbox{ }_j{\tilde a}_s$ and
\par $(\partial ^2F/\partial \mbox{ }_jx_v \partial \mbox{
}_kx_u)/2=\delta _{j,k}[\delta _{v,u} +\sum_p\mbox{ }_j{\tilde a}_t
\mbox{ }_j{\tilde a}_s]$, where $(i_u {\tilde i}_t)=i_p,$
$i_v{\tilde i}_s=i_p$, $\delta _{j,k}$ is the Kronecker function
such that $\delta _{j,k}=1$ for each $j=k$ and $\delta _{j,k}=0$ for
each $j\ne k$. The necessary condition of the minimum is: $(\partial
F/\partial \mbox{ }_jx_v)=0$ at $P_0$, where $x=\xi $, consequently,
\par $(i)$ $\mbox{ }_j\xi _v+\sum_p \zeta _p\mbox{ }_j{\tilde a}_s=0$ for each
$j=1,...,n-1$, \\   where $i_p=i_v{\tilde i}_s$. System $(i)$ of
algebraic equations has the unique solution if and only if $\zeta
\ne 0$, where $\zeta = \zeta _0i_0+...+\zeta _{2^r-1}i_{2^r-1}$,
$\zeta _v\in \bf R$, $\zeta \in \bf K$. That is, $a\in \bf K^{n-1}$
is defined uniquelly. This gives the characteristic tangent plane at
the point $P_0$ of the hypersphere with centre $z_0$ and radius
$|OP_0|$.
\par Take $\mbox{ }_jx_v=(\mbox{ }_jx_v-\mbox{ }_j\xi _v) +
\mbox{ }_j\xi _v$, then  due to $(i)$ \\
$F(x) = \{ \sum_{j,v} (\mbox{ }_j\xi _v)^2 + (\sum_p \zeta _p)^2 +
\sum_{j,v}(\mbox{ }_jx_v - \mbox{ }_j\xi _v)^2 + [\sum_{j,v,s}
(\mbox{ }_jx_v - \mbox{ }_j\xi _v)\mbox{ }_j{\tilde a}_s]^2 \} $ \\
$\ge \sum_{j,v} (\mbox{ }_jx_v)^2+ (\sum_p\zeta _p)^2$ and \\
$F(x)-F(x)|_{P_0}\ge \sum_{j,v}(\mbox{ }_jx_v - \mbox{ }_j\xi _v)^2$. \\
Thus the hyperplane $\sigma $ has not any other common point with
the hypersphere $\partial B({\bf K^n},z_0,|OP_0|)$ besides the point
$P_0$. \par Consider the characteristic plane $\pi $ such that
$OP_0\subset \pi $ and $\pi $ is given by the equation:
\par $<x,\alpha >=y\beta $, where $\alpha \in {\bf K^{n-1}}$, $\beta
\in \bf K$. Since $O\in \pi $, $P_0(\xi ,\zeta )\in \bf K^n$, then
$<\xi ,\alpha >=\zeta \beta $. The plane $\pi $ is not tangent to
$\sigma $, since $\sigma $ is contained in the hyperplane $\Pi $
tangent to the hypersphere $\partial B({\bf K^n},z_0,|OP_0|)$, but
$\Pi \perp OP_0$. The equation of $\Pi $ is: $\sum_{j,v}\mbox{
}_j\xi _v(\mbox{ }_jx_v-\mbox{ }_j\xi _v)+\sum_p \zeta _p(y_p-\zeta
_p)=0$. Associate with each $P\in \pi $ the translation $P_0P$ and
this gives the family of regular characteristics $\sigma $ in a
neighborhood of $P_0$. This family is given by the equations:
\par $(ii)$ $<g,\alpha >=h\beta $, \par $(iii)$ $<x-\xi -g,\alpha >=
y-\zeta - h$, \\
where the vector $P_0P$ has coordinates $(g,h)$, $g\in \bf K^{n-1}$,
$h\in \bf K$, $\sigma $ already does not pass through the origin
$O$. Therefore, $(ii)$ can be resolved relative to $\beta $ and this
gives the infinite family of charateristics $\sigma $ which has not
any common point with $E$, since it can be taken a shift on a vector
$z$ with $|z|>0$ from the hypersurface $\partial B({\bf
K^n},O,|OP_0|)$ into ${\bf K^n}\setminus B({\bf K^n},O,|OP_0|)$.
This contradicts common Theorem 54.
\par {\bf 57. Corollary.} {\it Let conditions of Theorem 56 be
satisfied, then the set $E$ has not any dense in itself subset $Y$
on a finite distance $dist (Y, E\setminus Y)$.}
\par {\bf Proof.} If any such $Y$ would be existent, then it would
be $P_0\in E$ with $dist (Y,P_0)>0$ contradicting Theorem 56, since
$\bf K^n$ is locally compact and $E\subset \bf K^n$.
\par {\bf 58. Corollary.} {\it Let conditions of Theorem 56 be
satisfied and $\cal F$ is normal at each point of a closed
pseudoconformal hypersurface, then $\cal F$ is normal at all points
inside this hypersurface.}
\par {\bf 59. Theorem.} {\it Let $C:=B({\bf K},0,R_1)\times ...
\times B({\bf K},0,R_{n-1})$, where $\partial G$ is a closed
pseudoconformal hypersurface in $\bf K$ and $G$ is a bounded domain,
$G\subset B({\bf K},0,R_n)$, $0<R_j<\infty $ for each $j$. If a
family $\cal F$ satisfies Condition 49(i) and $\cal F$ is normal at
each point $(i)$ $(0,y)$ with $y\in G$ or $(ii)$ $(x,y)$ with $x\in
Int (C)$ and $y\in \partial G$, then $\cal F$ is normal on $Int
(C\times G)$.}
\par {\bf Proof.} If this theorem would not be true at least at one
point $P(\xi ,\zeta )$ with $\xi \in Int (C)$ and $\eta \in Int
(G)$, where $\cal F$ is not a normal family, then $P\in E$ and we
put $|\mbox{ }_j\xi | =: R_{0,j}$. \par The family $\cal F$ is
normal at $(0,y)$ for $y\in G$, then $\cal F$ is normal in $\{
(x,y): |\mbox{ }_jx|\le R_{1,j}, y\in G \} $ for sufficiently small
$R_{1,j}>0$ for each $j=1,...,n-1$. Put $S_1 := \partial [B({\bf
K},0,R_{1,1})\times ...\times B({\bf K},0,R_{1,n-1})]$. Then by the
supposition there exists $(\xi ,\zeta )\in E$ such that
$R_{1,j}<|\mbox{ }_j\xi |<R_j$ for each $j=1,....,n-1$ and $\zeta
\in G$, where $G=Int (G)$, $cl (G) = G\cup \partial G$. Then $E$ has
no any point in $(S_1,G)$ and on its surface, consequently,
$E\subset cl \{ (x,y): R_{1,j}< |\mbox{ }_jx|<R_j, \forall
j=1,...,n-1; \zeta \in G \} $ $\bigcup \{ (x,y): |\mbox{ }_jx|=R_j$
$\mbox{for some j},$ $y\in cl (G) \} $. At the same time $E$ has no
any point $(x,y)$ with either $|\mbox{ }_jx|=R_{1,j}$ for some $j$
and $y\in G$ or $|\mbox{ }_jx|<R_j$, $y\in \partial G$ due to
supposition $(ii)$. \par Make the fractional $\bf R$-linear
pseudoconformal transformation $g$ such that \\  $X=k((\mbox{
}_1x)^{-1},...,(\mbox{ }_{n-1}x)^{-1})\in {\bf {\hat K}^{n-1}}$ and
$Y=y$ for each $\mbox{ }_jx\in \bf K$ for $j=1,...,n-1$ and $y\in
\bf K$, where $k=const >0$, $(X,Y)=g(x,y)$. If $f\in \cal F$, then
$f\circ g$ satisfies Condition 49(i) in $cl (\delta _1\times G)$,
where $\delta _1 := \{ x\in {\bf K^{n-1}}: k/R_j<|\mbox{
}_jx|<k/R_{1,j}, j=1,...,n-1 \} $. Moreover, the family ${\cal
F}\circ g$ is normal at each point $|\mbox{ }_jX|=R_{1,j}$ for some
$j=1,...,n-1$ with $Y\in \partial G$ and at each point $X\in \delta
_1$ and $Y\in
\partial G$. If $Q$ is a family of points, where ${\cal F}\circ g$
is not normal, then it is known, that on the one side $Q$ contains
the point $b := (k/\mbox{ }_1\xi ,...,k/\mbox{ }_{n-1}\xi ;\zeta
)\in (\delta _1\times G)$ and on the other side if \par $(iii)$
$(X,Y)\in Q\cap cl (\delta _1\times G)$, then there exists $1\le
j\le n-1$ with $|\mbox{ }_jX|=k/R_j$, $Y\in cl (G)$. Thus $dist
((X,Y); Q)< [\sum_{j=1}^{n-1} k^2R_j^{-2}+R_n^2]^{1/2}$ for $(X,Y)$
satisfying Condition $(iii)$, since $G\subset B({\bf K},0,R_n)$,
while $|(X,Y)| = [\sum_{j=1}^{n-1} k^2 R_{0,j}^{-2} + |\zeta
|^2]^{1/2}\ge k^2\sum_{j=1}^{n-1}R_{0,j}^{-2}$. Choose
$k^2\sum_{j=1}^{n-1} R_{0,j}^{-2}\ge [\sum_{j=1}^{n-1}k^2R_j^{-2} +
|\zeta |^2]^{1/2}$, that is, $k\ge |\zeta |[\sum_{j=1}^{n-1}
(R_{0,j}^{-2}-R_j^{-2})$, hence the distance from $Q$ to $0$ is
greater, than the distance from each $z\in Q\cap \partial (\delta
_1\times G)$ to zero in $\bf K^n$.
\par Consider the set $Q_1 := Q\cap cl (\delta _1\times G)$ which is
closed in $\bf K^n$. Therefore, there exists $z_0\in Q_1$ such that
$\max_{z\in Q_1} |z| = |z_0|$. Then $z_0\notin \partial Q_1$, since
$|z_0|>dist (\partial Q_1,0)$. Thus $Q$ has not any point from $Q_1$
in some neighborhood of $z_0$. Then $z_0$ satisfies conditions of
Theorem 56. This gives the contradiction to the supposition about
the point $P$, from which this theorem follows.
\par {\bf 60. Corollary.} {\it If a family $\cal F$ is normal at
each point $(x,y)$ with $|\mbox{ }_jx|<R_j$ for each $j=1,...,n-1$
and $y\in \partial G$ (see Theorem 59). If also at some such marked
point $(x_0,y_0)$ the family $\cal F$ is not normal, then for each
$\xi $ such that $|\mbox{ }_j\xi |<R_j$ for each $j=1,...,n-1$,
there exists at least one value $\zeta \in cl (G)$ such that $\cal
F$ is not normal at $(\xi ,\zeta )$.}
\par {\bf 61. Corollary.} {\it If a family $\cal F$ from Theorem 59
is normal $(i)$ at each point $(x,y)$ with $x\in Int (C)$ and $y\in
\partial G$; $(ii)$ at each point $(x_0,y)$ with $x_0\in Int(G)$ and
$y\in cl (G)$, where $x_0$ is a marked point, then the family $\cal
F$ is normal at each point $(x,y)$ such that $x\in Int (C)$ and
$y\in G$.}
\par {\bf Proof.} The point $x_0$ can be joined with each point
$x\in  Int (C)$ by the segment. For each infinite covering of
$[x_0,x]$ by balls there exists its finite subcovering, since the
segment $[x_0,x]$ is compact. Thus, there exists a finite family of
balls $B_1,...,B_m$ in $Int (C)$ covering the segment $[x_0,x]$
which satisfies conditions below. In view of Theorem 59 the family
$\cal F$ is normal at each point $(x,y)$ such that $x\in Int (B_1)$
and $y\in cl (G)$. Then take $x_2\in Int (B_1)\cap Int (B_2)$. The
applying the same theorem gives, that the family $\cal F$ is normal
at each $(x,y)\in Int (B_2)\times G$, partucularly for $x_3\in Int
(B_2)\cap Int (B_3)$, and this process can be continued. Thus $\cal
F$ is normal at each point $(x,y)\in Int (C)\times G$.
\par {\bf 62. Note.} Now study the problem, when a $C^2$
hypersurface $S$ defined by the equation $\omega (x,y)=0$ can belong
to the set $E$, where the family $\cal F$ terminates to be normal
(see \S 59), where $\omega $ is a real-valued function of $C^2$
class of smoothness by all variables $(\mbox{ }_1x_0,...,\mbox{
}_{n-1}x_{2^r-1},y_0,...,y_{2^r-1})$. In some domain $U$ in $\bf
K^n$ consider a point $P$ and its neigborhood $V_P$ which is crossed
by the hypersurface $S$. The hypersurface $S$ divides the domain on
two subdomains $U_- := \{ (x,y): \omega (x,y)<0 \} $ and $U_+ := \{
(x,y): \omega (x,y)>0 \} $. We can suppose, that $P$ is a regular
point of $S$. Consider the case, when
\par $\sum_{v=0}^{2^r-1} (\partial \omega (x,y)/
\partial y_v)^2>0$ at the point $P$. \par Then the question arise:
does a family $\cal F$ exist normal in the domain $U\setminus S$ and
terminating to be normal at each point of the hypersurface $S$? It
is demonstrated below, that this is impossible in general and that
such family can exist only in one of the two domains either $U_+$ or
$U_-$.
\par Study at first necessary conditions to which the hypersurface $S$
need to satisfy that a family $\cal F$ would be normal in $U_+$ and
terminating to be normal at each point of $S$. For this apply
general Theorem 54. Thus let $\cal F$ be normal in $U_+$ and let
$\cal F$ be not normal at each point of $S$. Suppose, that we can
construct the characteristic surface $\sigma $, containing point
$P$, regular at $P$, tangent at $P$ to the hypersurface $S$ and
having no any other common point with $S$ in a neighborhood $V_P$ of
$P$. Then either $(\sigma \cap V_P\setminus \{ P \} ) \subset
U_+=:U_P$ or $(\sigma \cap V_P\setminus \{ P \} ) \subset U_-=:U_P$.
Prove, that $\sigma $ need to satisfy the second inclusion, when
$\cal F$ is normal in $U_+$ and $\cal F$ is not normal at each point
of $S$.  \par In fact there exists a family of regular
characteristics around $P$ attached to $P_1\in V_P$, where $P_1$
belongs to the characteristic plane $\pi $ over $\bf R$ embedded
into $\bf K^n$, such that $\pi $ contains the normal line $\nu $
over $\bf R$ at $P$ to $S$. In view of Theorem 55 there exists this
family and it is regular at each point in a sufficiently small
neigborhood $V_P$. If a preceding translation $PP_1$ is sufficiently
small and has a nonzero projection on $\nu $ and $P_1\in U_P$, then
the obtained charateristic has not any common point with $S\cap
V_P.$ Thus if $\sigma \subset U_+$, where the family of functions
$\cal F$ is normal, then there would be characteristic surfaces of
such family in a neighborhood of $\sigma $ and having no any common
point with $S$. On these characteristics the family $\cal F$ will be
normal in $V_P\setminus \{ P \} $, since the unique points where
$\cal F$ may be not normal are boundary points and that of $S$. This
is impossible due to Theorem 54. Indeed, if a family $\cal F$ is
normal in $V_P\setminus \{ P \} $, then it would be normal on a
characteristic near to $\sigma $ besides at least at one point
nearby to $P$, where $\cal F$ terminates to be normal. Thus $\sigma
$ must be contained in $U_-$.
\par Take the characteristic surface $\sigma $ tangent to $S$ at $P$
such that $\sigma \cap (V_P\setminus \{ P \} ) =\emptyset $.
Introduce the notations: $\Delta _{1,x} := \sum_{j=1}^{n-1}
\sum_{v=0}^{2^r-1} (\partial \psi (x,y)/\partial \mbox{ }_jx_v)^2$,
$\Delta _{1,y} := \sum_{v=0}^{2^r-1}(\partial \psi (x,y)/\partial
y_v)^2$, $\Delta _{2,x} := \sum_{j=1}^{n-1} \sum_{v=0}^{2^r-1}
(\partial ^2\psi (x,y) /\partial \mbox{ }_jx_v^2)$ and $\Delta
_{2,y} := \sum_{v=0}^{2^r-1} (\partial ^2\psi (x,y)/\partial
y_v^2)$. Suppose that $\Delta _{1,y}>0$ at $P$. Then $\sigma $ is
defined by the decomposition into the noncommutative analog of the
Taylor expansion in a neighborhood of $P(\xi ,\zeta )$:
\par $(i)$ $y-\zeta = \sum_{|m|\ge 0} \{ a_m,(x-\xi )^m \} _{q(2|m|)}$, \\
where $\psi $ is holomorphic by each variable $(x,y)$ at each
regular point, in particular, at $P$ (see \S 53), also $\psi (\xi
,\zeta )=0$. Supposing that $\partial \psi (x,y)/\partial y$ is the
invertible operator such that all conditions of the implicit
function theorem are satisfied (see Theorem X.7 \cite{zorich})
rewrite $(i)$ in real variables:
\par $(ii)$ $y_p - \zeta _p = \sum_{j,v}a_{j,v,p}(\mbox{ }_jx_v-
\mbox{ }_j\xi _v) + \sum_{j,v}b_{j,v,p}(\mbox{ }_jx_v-\mbox{ }_j\xi
_v)^2 $  \par $+ 2 \sum_{j<k; v,u} c_{j,k,v,u,p} (\mbox{
}_jx_v-\mbox{ }_j\xi _v) (\mbox{ }_kx_u-\mbox{ }_k\xi _u) +...$, \\
where $a_{j,v,p}, b_{j,v,p}, c_{j,k,v,u,p}\in \bf R$, $j,
k=1,...,n-1$, $v, u=0,1,...,2^r-1$. To know the sign of $\omega $ on
$\sigma $ such that the sign in $V_P\cap \sigma $ need to be
constant it is necessary to substitute $y$ on $x$ using $(ii)$ and
insert into $\omega $, hence consider the composite function $\omega
(x,\zeta + \phi (x-\xi ))$, where $y-\zeta =\phi (x-\xi )$ on $S$ in
a neigborhood of $P$. Decompose the function $\omega (x,\zeta + \phi
(x-\xi ))$ with the help of the Taylor theorem in the Peano form up
to the second order terms:
\par $(iii)$ $\omega (x,\zeta + \phi (x-\xi )) = \omega (\xi ,\zeta)+
\sum_{j,v}(\partial \omega (\xi ,\zeta )/\partial \mbox{ }_j\xi _v)
(\mbox{ }_jx_v - \mbox{ }_j\xi _v) $ \par $+ \sum_{j,v,u} (\partial
\omega (\xi ,\zeta )/\partial \zeta _u) (\partial y_u(\xi )/\partial
\mbox{ }_j\xi _v)(\mbox{ }_jx_v - \mbox{ }_j\xi _v)$
\par $+\sum_{j,k,v,u} (\partial ^2 \omega (\xi ,\zeta )/\partial
\mbox{ }_j\xi _v\partial \mbox{ }_k\xi _u) (\mbox{ }_jx_v - \mbox{
}_j\xi _v)(\mbox{ }_kx_u - \mbox{ }_k\xi _u)/2$ \par $+
\sum_{j,k,v,u,p} (\partial ^2 \omega (\xi ,\zeta )/\partial \mbox{
}_j\xi _v\partial \zeta _p) (\partial y_p(\xi )/\partial \mbox{
}_k\xi _u)(\mbox{ }_jx_v - \mbox{ }_j\xi _v)(\mbox{ }_kx_v - \mbox{
}_k\xi _u)$
\par $+\sum_{j,k,v,u,p} (\partial \omega (\xi ,\zeta )/\partial
\zeta _p)(\partial ^2y_p(\xi )/\partial \mbox{ }_j\xi _v \partial
\mbox{ }_k\xi _u) (\mbox{ }_jx_v - \mbox{ }_j\xi _v)(\mbox{ }_kx_u -
\mbox{ }_k\xi _u)/2$
\par $+ \sum_{j,k,v,u,p,s} (\partial ^2 \omega (\xi ,\zeta )/
\partial \zeta _p\partial \zeta _s)(\partial y_p(\xi )/\partial
\mbox{ }_j\xi _v) (\partial y_s(\xi )/\partial \mbox{ }_k\xi _u)
(\mbox{ }_jx_v - \mbox{ }_j\xi _v)(\mbox{ }_kx_u - \mbox{ }_k\xi
_u)/2 + o(|x-\xi |^2)$.
\par For the marked point $P(\xi ,\zeta )$
that $\sigma \cap V_P\setminus \{ P \} $ be contained in $U_+
\setminus \{ P \} $ or in $U_- \setminus \{ P \} $ it is necessary,
that \par $(iv)$ $\partial \omega (x,y)/\partial \mbox{ }_jx_v
+\sum_p(\partial \omega (x,y)/\partial y_p)(\partial
y_p(x)/\partial \mbox{ }_jx_v)=0$ \\
at $P(\xi ,\zeta )$ for each $j=1,...,n-1$ and $v=0,1,...,2^r-1$,
which is the system of equations for real numbers $a_{j,v,p}
=\partial y_p(x)/\partial \mbox{ }_jx_v $. The family $\sigma (P_1)$
of characteristics tangent to $S$ at $P_1$ parametrized by $P_1\in
V_P\cap S$ has the envelope $S$. Thus Condition $(iv)$ need to be
satisfied at each point $P_1\in V_P\cap S$ and it can be written in
the form:
\par $(v)$ $\partial \omega (x,y)/\partial \mbox{ }_jx_v
+ (\partial \omega (x,y)/\partial y)(\partial y(x)/\partial
\mbox{ }_jx_v)=0$, \\
where $\beta := \partial \omega (x,y)/\partial y :=
\sum_{p=0}^{2^r-1} (\partial \omega (x,y)/\partial y_p) {\tilde
i}_p$, \\
hence $0\ne \beta \in \bf K$, since $\Delta _{1,y}>0$ and the
function $\omega $ is real-valued. Thus System of equations $(v)$
has the solution:
\par $(vi)$  $(\partial y(x)/\partial \mbox{ }_jx_v) = -
b(\partial \omega (x,y)/\partial \mbox{ }_jx_v)$ \\
for each $j$ and $v$, where $b=\beta ^{-1}\in \bf K$. Therefore,
\par $(vii)$ $(\partial y_p(x)/\partial \mbox{ }_jx_v) = -
b_p(\partial \omega (x,y)/\partial \mbox{ }_jx_v)$ \\
for each $p=0,1,...,2^r-1$. From Equation $(vii)$ it follows:
\par $(viii)$ $\partial ^2y_p(x)/\partial \mbox{ }_jx_v
\partial \mbox{ }_kx_u = - b_p (\partial ^2\omega (x,y)/
\partial \mbox{ }_jx_v \partial \mbox{ }_kx_u)
- (\partial b_p/\partial \mbox{ }_kx_u) (\partial \omega
(x,y)/\partial \mbox{ }_jx_v)$. \\
On the other hand, $(\partial b/\partial \mbox{ }_kx).h = (\partial
\beta ^{-1} / \partial \mbox{ }_kx).h = - b((\partial ^2 \omega
(x,y)/\partial \mbox{ }_kx \partial y).h)b)$, consequently,
\par $(ix)$ $(\partial b_p/\partial \mbox{ }_kx_u) =
- \sum_{(s,t; i_c=({\tilde i}_t{\tilde i}_s)i_p)}b_s((\partial ^2
\omega (x,y)/\partial \mbox{ }_kx_u \partial y_t))b_c)$ \\
due to the alternativity of $\bf K$: $e(ye)=(ey)e$ for each $e, y
\in \bf K$. Thus from $(viii,ix)$ it follows:
\par $(x)$ $\partial ^2y_p(x)/\partial \mbox{ }_jx_v
\partial \mbox{ }_kx_u = - b_p (\partial ^2\omega (x,y)/
\partial \mbox{ }_jx_v \partial \mbox{ }_kx_u)$ \par
$+ \sum_{(s,t; i_c=({\tilde i}_t{\tilde i}_s)i_p)}b_s((\partial ^2
\omega (x,y)/\partial \mbox{ }_kx_u \partial y_t))b_c) (\partial
\omega (x,y)/\partial \mbox{ }_jx_v)$. \\
\par Then the sufficient
condition is that the remaining quadratic form is positive or
negative definite respectively. Since $\psi $ is holomorphic, then
$y= y(x) := \zeta + \phi (x-\xi )$ is holomorphic by the $x$
variable for each regular point $P(\xi ,\zeta )$, hence $\phi $ is
locally analytic by $x$.
\par Therefore, the problem is whether the quadratic form $Q(x)$
is positive or negative definite, which is given by the equation:
\par $(xi)$  $Q(x)=\sum_{j,k=0}^{n-1}\sum_{v,u=0}^{2^r-1}
\mu ^{j,k,v,u}(\mbox{ }_jx_v - \mbox{ }_j\xi _v)(\mbox{ }_kx_u -
\mbox{ }_k\xi _u)$, where
\par $(xii)$  $\mu ^{j,k,v,u} :=
(\partial ^2 \omega (\xi ,\zeta )/\partial \mbox{ }_j\xi _v\partial
\mbox{ }_k\xi _u) /2$ \par $+ \sum_p (\partial ^2 \omega (\xi ,\zeta
)/\partial \mbox{ }_j\xi _v\partial \zeta _p) a_{k,u,p} $
\par $+\sum_p (\partial \omega (\xi ,\zeta )/\partial
\zeta _p)(\partial ^2y_p(\xi )/\partial \mbox{ }_j\xi _v \partial
\mbox{ }_k\xi _u) /2$
\par $+ \sum_{p,s} (\partial ^2 \omega (\xi ,\zeta )/
\partial \zeta _p\partial \zeta _s)a_{j,v,p} a_{k,u,s}) /2.$ \\
With the help of Equations $(vi,x)$ Equation $(xii)$ simplifies:
\par $(xiii)$  $\mu ^{j,k,v,u} :=
(\partial ^2 \omega (\xi ,\zeta )/\partial \mbox{ }_j\xi _v\partial
\mbox{ }_k\xi _u) (1 - \sum_p (\partial \omega (\xi ,\zeta
)/\partial \zeta _p)b_p) /2$
\par $ - \sum_p (\partial ^2 \omega (\xi ,\zeta )/\partial \mbox{
}_j\xi _v\partial \zeta _p) b_p (\partial \omega (\xi ,\zeta
)/\partial \mbox{ }_kx_u)$
\par $+ \sum_{p,s} (\partial ^2 \omega (\xi ,\zeta )/
\partial \zeta _p\partial \zeta _s)b_p (\partial \omega (\xi ,\zeta
)/\partial \mbox{ }_jx_v)b_s (\partial \omega (\xi ,\zeta )/\partial
\mbox{ }_kx_u)) /2$
\par $ + \sum_{(p; i_c=({\tilde i}_t{\tilde i}_s)i_p)}
(\partial \omega (\xi ,\zeta )/\partial \mbox{
}_jx_v) (\partial \omega (\xi ,\zeta )/\partial \zeta _p) (\partial
^2 \omega (\xi ,\zeta )/\partial \mbox{ }_kx_u\partial
y_t)b_sb_c)/2.$
\par We have that $\mu ^{j,k,v,u}$ form the symmetric
matrix with indices of rows $(j,v)$ and indices of columns $(k,u)$.
It is well-known from the linear algebra that the quadratic form
over $\bf R$ can be reduced to the diagonal canonical form by an
orthogonal transformation (see Theorem X.5.7 in \cite{gant}).
Therefore, the general quadratic form $Q(x)$ can be transformed to
the diagonal canonical form by the pseudoconformal transformations
$x'=g(x)$ and $y=y(g^{-1}(x'))$, since the composition of
pseudoconformal mappings is pseudoconformal (see Theorems 4-6
above). At the same time a pseudoconformal function transforms a
domain into a domain by Theorem 30.
\par The characteristic surface $\sigma $ in a neighborhood
$V_P\setminus \{ P \} $ is not contained in $U_+$ if $Q(x)$ is
negative definite, also $\sigma \cap V_P\setminus \{ P \} $ is not
contained in $U_-$ if $Q(x)$ is positive definite. Thus the
following theorem is demonstrated.
\par {\bf 63. Theorem.} {\it Let suppositions of \S 62 be satisfied.
Then the necessary condition for the existence of a family $\cal F$
normal in $U_+$ or $U_-$ and terminating to be normal on $S$ which
is contained in $\partial U_+$ or $\partial U_-$ is $Q(x)<0$ or
$Q(x)>0$ on $S$ respectively, where $Q(x)$ is given by Equations
62(iv,v,vii). For an existence of a family $\cal F$ in $U_+$ and
$U_-$ simultaneously it is necessary, that $Q(x)=0$ on $S$.}
\par {\bf 64. Note.} While proving Theorem 50 if the family $\cal
F$ is normal at each point $(x,y)$ with $y=0\in \bf K$ in a
neighborhood $U$ of $(0,0)$ besides $(0,0)$, then there exists a
sequence of functions $f_1(x,y)$, $f_2(x,y)$,... , for which each
equation $f_n(0,y)=0$ in the ball $B({\bf K},0,R)$ of sufficiently
small radius $R>0$ has at least one zero $y_n$, then $\lim_{n\to
\infty }y_n=0$. Let now $\sigma $ be a characteristic surface and
$P$ be a regular point in $\sigma $. Then this property is also
true, when $\cal F$ is normal at each point in a neighborhood $U_P$
of $P$ besides $P$ itself on a characteristic surface $\sigma $,
then $P$ is the limit point of the sequence $P_n(x_n,y_n)\in \sigma
$ at which $f_n(x_n,y_n)=0$. Thus each point from $E$ is the limit
point of the family of points $(x,y)$ such that $f(x,y)=0$.
\par {\bf 65. Theorem.} {\it Let $ \{ f_m (x,y): m\in {\bf N} \} $
be a sequence of functions satisfying Condition 49(i) in a domain
$U$ satisfying conditions of Remark 12 such that the set $\{ (x,y):
f_m(x,y)=0 \} =: T_m$ has an infinite family of limit points forming
a connected subset. Then there exists a sequence $G := \{ g_m(x,y):
n \} $ of functions satisfying Condition 49(i) and having the same
zeros as $f_m$ and $G$ is normal on $U$ besides a limit point $P$ of
$ \{ T_m: m \} $.}
\par {\bf Proof.} If a family $\{ T_m: m \} $ of sets $T_m$ has a
limit point $P_0\in U$, then the family of such limit points is
infinite forming a connected set $T$ containing $P_0$. Consider the
hypersphere $S := \partial B({\bf K^n},P_0,R)$ contained in $U$,
where $0<R<\infty $ and such that $P_0$ is the limit point of $\{
T_m\cap Int (B({\bf K^n},P_0,R)): m \} $. Then there exists an
infinite subsequence $\{ T_{m_j}: j \} $ such that each $T_{m_j}$
crosses $S$, since otherwise the function $F_{m_j}:=1/f_{m_j}$ would
be satisfying Condition 49(i) on $Int (B({\bf K^n},P_0,R))$ due to
Theorem 2.12 \cite{luladfcdv}. Therefore, $T_{m_j}\cap Int (B({\bf
K^n},P_0,R))$ is infinite for each $j\in \bf N$. Therefore, $ \{
T_{m_j}: j \} $ has not less than one limit point on each
hyperpshere $\partial B({\bf K^n},P_0,R')$ for each $R'$ such that
$0<R' \le R$.
\par To each point $(x,y)\in T_m$ there corresponds the hypersphere
$S_m(x,y) := \{ (X,Y): |X-x|^2 + |Y-y|^2 = R_m^2 \} $, where
$0<R_m<\infty $. Consider the ball $B_m(x,y)$ bounded by $S_m(x,y)$.
Then the set $C_m$ of all points $(X,Y)$ in $S_m\cup (V\setminus
B_m(x,y))$ at which $f_m(X,Y)\ne 0$ is closed. Then $\inf_{z\in C_m}
|f_m(z)|=:a_m>0$ and put $g_m(x,y):=mf_m(x,y)/a_m$. Each $g_m$ has
the same zeros as $f_m$ and satisfies Condition 49(i) on $V$. Prove
that the only points of $G:= \{ g_m: m \} $, where $G$ is not normal
are limit points of $T_m$.
\par Let $Q\in V$ is not a limit point of the sequence $ \{ T_m: m
\} $. Then there exists $m_0\in \bf N$ and the ball $B({\bf
K^n},Q,\delta )$ of the radius $\delta >0$ such that $Int (B({\bf
K^n},Q,\delta ))$ does not contain any point from $T_m$ for each
$m\ge m_0$. Then $Q$ belongs to ${\bf K^n}\setminus B({\bf
K^n},Q,t_m)$, where $0<t_m<\delta $ and $\lim_{m\to \infty }t_m=0$.
Therefore, there exists $m_1\in \bf N$ such that $B({\bf
K^n},Q,\delta /2)\cap B_m(x,y)=\emptyset $ for each $m\ge m_1$. Thus
on $\sigma := \partial B({\bf K^n},Q,\delta /2)$ the sequence $G$
converges to the infinity uniformly, hence $G$ is normal on $\sigma
$. \par Let $P(x,y)$ be a limit point of $\{ T_m: m \} $. Then there
exists a sequence $m_1<m_2<...<m_p<...$ such  that $\inf_{z\in
T_{m_p}} |z-(x,y)| =: R_{m_p}$ and $\lim_{p\to \infty }R_{m_p}=0$.
If $G$ would be normal at $P$, then from $ \{ g_{m_p}: p \} $ it
would be possible to extract a subsequence $ \{ g_{l_p}: p \} $
converging uniformly relative to the metric $\rho $ on some
hypersphere $\sigma _0 := \partial B({\bf K^n},P,R')$ of
sufficiently small radius $R'$, $0<R'<\infty $. Since $l_p$ are
chosen among $m_p$, then there exists $p_0$ such that $T_{l_p}\cap
\sigma _0\ne \emptyset $ for each $p>p_0$. This happens for
$R_{l_p}<R'$. On each $T_{l_p}$ we have $g_{l_p}=0$. \par Then it is
impossible, that the limit of the sequence $g_{l_p}$ would be the
infinity while $p$ tends to the infinity, since on $\sigma _0$ there
are points from $T_{l_p}$ for each $p>p_0$. Thus the limit function
$g=\lim_{p\to \infty }g_{l_p}$ on $\sigma _0$ satisfies Condition
49(i) due to Lemma 11, where $g(P)=0$. If it would be $g(P)=:a\ne
0$, then for each $\epsilon >0$ there exists a hypersphere $\sigma
_0$ with the centre $P$ on which $|g(z)-a|<\epsilon $. Take
$\epsilon =|a|/2$, then $|g(z)|>|a|/2$ on $\sigma _0$. Since
$g_{l_p}$ converges to $g$ relative to the metric $\rho $ on $\sigma
_0$, then for each $\epsilon _1>0$ there exists $p_0\in \bf N$ such
that $|g_{l_p}-g|<\epsilon _1$ on $\sigma _0$ for each $p>p_0$. Take
$\epsilon _1=|a|/4$. Then there exists $p_1\in \bf N$ such that on
$\sigma _0$ there is the inequality: $|g_{l_p}|>|a|/4$ for each
$p>p_1$. But this is impossible, since $\lim_{p\to \infty
}R_{l_p}=0$ and there exists $p_2\in \bf N$ such that $T_{l_p}\cap
\sigma _0\ne \emptyset $ for each $p>p_2$, consequently, on $\sigma
_0$ there are points at which $g_{l_p}(z)=0$.
\par Therefore, for each $\epsilon >0$ there exists
$\sigma =\partial B({\bf K^n},P,t)$ for some $0<t<\infty $ such that
$|g(z)|<\epsilon $ on $\sigma $. Since $g_{l_p}$ converges to $g$
relative to the metric $\rho $ on $\sigma $, then there exists
$p_0\in \bf N$ such that $|g_{l_p}(z)|<2\epsilon $ on $\sigma $ for
each $p>p_0$. But this is impossible. In fact there exists $p_1\in
\bf N$ such that $0<R_{l_p}<t<\infty $ for each $p>p_1$ and
inevitably $T_{l_p}\cap \sigma \ne \emptyset $.
\par Therefore, the sequence $T_m\cap \partial B({\bf K^n},P,R_m)=:V_m$
satisfies $\bigcap_{m=1}^{\infty }V_m= \{ P \} $. Thus there exists
$p_0\in \bf N$ such that there exists a point $(x_p,y_p)\in B({\bf
K^n},P,t)\setminus V_{m_p}$ at which $|g_{l_p}(x_p,y_p)|>l_p$ and
this contradicts the inequality $|g_{l_p}(z)|<2\epsilon $ on $\sigma
$. Thus the supposition that $g$ is normal at $Q$ leads to the
contradiction, consequently, $\{ g_m: m \} $ is not normal at $Q$.
Therefore, the set $E$ relative to $ \{ g_m: m \} $ is the family of
limit points of the sequence $ \{ T_m: m \} $.
\par {\bf 66. Example.} Consider the domain $U$ in $\bf K^2$ bounded
by the hypersurface $S$ given by the equation $R_2 = f(R_1)$, where
$R_1^2 = \sum_vx_v^2$ and $R_2^2 = \sum_vy_v^2$, $f\in C^2$. That is
each internal point in $U$ is charaterized by the condition:
$f(|x|)-|y|>0$. Suppose that $f_{k,m} :=
((a_{k,m}x^k)b_{k,m})((c_{k,m}y^m)d_{k,m}) $, where $x, y \in \bf
K$, $\lim_{k+m\to \infty
}((a_{k,m}x^k)b_{k,m})((c_{k,m}y^m)d_{k,m})=0$ for each $(x,y)\in
U$, where $a_{k,m}, b_{k,m}, c_{k,m}, d_{k,m}\in \bf K$, $0\le k,
m\in \bf Z$. Each function $f_{k,m}$ is either pseudoconformal by
$x$ and $y$ on $U$ due to Theorems 4-6 and Lemma 2 or equal to zero.
Then $ \{ f_{k,m}: k,m \} $ is the normal family on $U$, but it is
not normal on $\partial U$, since for each point $(x,y)$ in a
neighborhood $Int (B({\bf K^2},(x_0,y_0),\epsilon ))$ of
$(x_0,y_0)\in
\partial U$ with $|x|>R_1$ and $|y|>R_2$ there is the subsequence
$f_{k,m}(x,y)$ converging to the infinity while it is converging to
zero in each internal point of $U$, since $\epsilon >0$ is
arbitrary. With this $\omega (x,y)=f(|x|)-|y|$ it is possible to
calculate $Q(x)$ more concretely:
\par $\partial \omega (x,y)/\partial \mbox{ }_jx_v=f'\mbox{
}_jx_v|x|^{-1}$;
\par $\partial \omega (x,y)/\partial y_v=-y_v|y|^{-1}$;
\par $\partial ^2\omega
(x,y)/\partial \mbox{ }_jx_v\partial y_u=0$;
\par $\partial ^2\omega (x,y)/\partial \mbox{ }_jx_v\mbox{ }_kx_u=
f''\mbox{ }_jx_v\mbox{ }_kx_u|x|^{-2}-f'\mbox{ }_jx_v\mbox{ }_kx_u
|x|^{-3}+f'|x|^{-1}\delta _{j,k}$;
\par $\partial ^2\omega (x,y)/\partial y_v\partial y_u=
y_vy_u|y|^{-3}-\delta _{v,u}|y|^{-1}$ for each $j, k=1,...,n-1$, $v,
u=0,1,...,2^r-1$, where $r=2$ for ${\bf K}=\bf H$ and $r=3$ for
${\bf K}=\bf O$. In this case $n=2$ and indices $j, k$ can be
omitted, $b= - {\tilde y}|y|^{-1}$. Therefore, due to Equations
62(xi,xiii) the quadratic form $Q(x)$ has real coefficients:
\par  $\mu ^{v,u} = [f'' x_vx_u |x|^{-2} - f' x_vx_u |x|^{-3}
+f' \delta _{u,v} |x|^{-1}] y_0^2|y|^{-2} $ \par $+ \sum_{p,s}
(y_py_s|y|^{-3} - \delta _{p,s}|y|^{-1}) y_p^2y_s^2|y|^{-4}
(f')^2x_vx_u|x|^{-2}$ \\
for each $v, u=0,1,...,2^r-1$, since $1- \sum_py_pb_p |y|^{-1}=
(|y|^2 + y_0^2 - y_1^2 - y_2^2-...-y_{2^r-1}^2)|y|^{-2}=2y_0^2$,
$b_p^2=y_p^2|y|^{-2}$ for each $p$. Then the condition $Q(x)\le 0$
can be satisfied on $S\cap V_p$ with certain $f$. For example, at
the point $P(\xi ,\zeta )$ with $\xi =(R_1,0,...,0)$ and $\zeta =
(R_2,0,...,0)$ for $R_1>0$ and $R_2>0$ and due to continuity on
$S\cap V_P$ for sufficiently small neighborhood $V_P$, since
$\sum_{p,s}(y_py_s|y|^{-3} - \delta
_{p,s}|y|^{-1})y_p^2y_s^2|y|^{-4}=\sum_{p,s}(y_p^3y_s^3 -
y_p^2y_s^4) |y|^{-7}=\sum_{p<s}y_p^2y_s^2 (2y_py_s-y_p^2 -
y_s^2)|y|^{-7} <0$ for $y\ne 0$. At this point $Q(x)<0$ if
\par $(i)$ $R_2(R_1f'' + f') \le R_1(f')^2$. Denoting $ln (f(R)) =:
\Phi (t)$, where $t=ln (R)$, $f(R)>0$ on its domain, the Condition
$(i)$ takes the form: $\Phi ''(t)\le 0$.

\section{Groups of pseudoconformal diffeomorphisms.}
\par {\bf 1. Definition.} Let $M$ be a manifold on $\bf K^n$, where
$n\in \bf N$, such that its atlas is $At (M) = \{ (U_j,\phi _j): j
\} $, each $U_j$ is open in $M$, $\bigcup_jU_j=M$, $\phi _j: U_j \to
\phi (U_j)\subset \bf K^n$ is a homeomorphism for each $j$, while
each connecting mapping of charts $\phi _i\circ \phi _j^{-1}$ is
holomorphic on its domain, when $U_i\cap U_j\ne \emptyset $. Then we
call $M$ the $\bf K$ holomorphic manifold, where ${\bf K}=\bf H$ or
${\bf K}=\bf O$.
\par {\bf 2. Theorem.} {\it Let $M$ be a $\bf K$ holomorphic
manifold. Then the family $DifH (M)$ of all holomorphic
diffeomorphisms of $M$ form the group.}
\par {\bf Proof.} A mapping $f\in DifH (M)$ if and only if
$f_{i,j} := \phi _i\circ f\circ \phi _j^{-1}$ is holomorphic on its
domain for each $i, j$. Then $f_{i,j}\circ f_{j,k}=f_{i,k}$ for each
$i, j, k$. Compositions of functions are associative in the sense of
the theory of sets. Thus the proof of this theorem follows from the
first part of the proof of Theorem 2.6, since compositions of
holomorphic functions are holomorphic and the inverse $f^{-1}$ of a
holomorphic function $f\in DifH (M)$ is holomorphic.
\par For further proceedings it is necessary to prove an auxiliary
theorem about differential equations in the class of holomorphic
functions over quaternions or octonions. The problem is the
following.
\par {\bf 3.} Suppose there is given a system of partial differential
equations:
\par $(1)$ $\partial ^{n_j}u_j/\partial
t^{n_j}=F_j(t,x_1,...,x_n,u_1,...,u_N,...,\partial ^ku_l/\partial
t^{k_0}\partial x_1^{k_1}...\partial x_n^{k_n},...),$ \\
where $j, l=1,...,N$; $k_0+k_1+...+k_n=k\le n_l$; $k_0<n_l$, $N, n
\in \bf N$, $0\le k\in \bf Z$, $(x_1,...,x_n)\in U$, $x_j\in \bf K$
for each $j$, $U$ is a domain in $\bf K^n$. If $t$ is real, then
$F_j$ are $\bf K$-valued functions or functionals. If $t\in \bf K$,
then each $F_j$ is a mapping with values in the space $L_q({\bf
K}^{\otimes n_j};{\bf K})$ of all $\bf R$ $n_j$-polyhomogeneous $\bf
K$-additive operators with values in $\bf K$. Each partial
derivative $\partial ^ku_l/\partial t^{k_0}\partial
x_1^{k_1}...\partial x_n^{k_n}$ is the mapping with values in
$L_q({\bf K}^{\otimes s};{\bf K})$, where $s=k-k_0$ for $t\in \bf R$
and $s=k$ for $t\in \bf K$, ${\bf K}=\bf H$ or ${\bf K}=\bf O$.
\par For $t=t_0$ there are given initial conditions:
\par $(2)$ $(\partial ^su_j/\partial t^s)|_{t=t_0}=\phi
_j^{(s)}(x_1,x_2,...,x_n)$ \\
for each $s=0,1,...,n_j-1$ and each $(x_1,...,x_n)\in U$, where
$\partial ^0u_j/\partial t^0:=u_j$.
\par A solution of the problem $(1,2)$ is searched in a neighborhood
$V$ of the given point $(t_0,x_1^0,...,x_n^0)$ either in ${\bf
R}\times U$ or in ${\bf K}\times U$. Introduce the notation:
\par $(\partial ^{k-k_0}\phi _j^{(k_0)}/\partial
x_1^{k_1}...\partial x_n^{k_n})|_{(x_1=x_1^0,...,x_n=x_n^0)}= \phi
^0_{(j,k_0,k_1,...,k_n)}$ \\
for each $j=1,...,N$, $k_0+k_1+...+k_n=k\le n_j$. \par $(3)$.
Suppose that all functions $F_j$ are analytic by all variables and
holomorphic by $x_1,...,x_n;$  $u_1,...,u_N$, also holomorphic by
$t$, when $t\in \bf K$, in a neighborhood $Y$ of \\
$(t_0,x_1^0,...,x_n^0,...,\phi ^0_{j,k_0,k_1,...,k_n},...)$, where
analyticity by operators $A_{l,k} := \partial ^ku_l/\partial
t^{k_0}\partial x_1^{k_1}...\partial x_n^{k_n}$ is defined relative
to their suitable compositions and the convergence of a series of
$F_j$ by $\{ A_{l,k}: l,k \} $ is supposed to be relative to the
operator norm, where $A_{l,k}\in \bigcup_{m=1}^{\infty }L_q({\bf
K}^{\otimes m},{\bf K})$. Let also all functions $\phi _j^{(k)}$ be
holomorphic in a neighborhood $X$ of $(x_1^0,...,x_n^0)$.
\par {\bf 4. Theorem.} {\it The problem $1(1-3)$ has a holomorphic
solution either by all variables $(t,x_1,...,x_n)$ in some
neighborhood of $(t_0,x_1^0,...,x_n^0)$, when $t\in \bf K$, or
holomorphic by the variables $(x_1,...,x_n)$ and analytic by $t$,
when $t\in \bf R$. This solution is unique in such class of
functions.}
\par {\bf Proof.} Introduce new functions, when $n_j>1$,
such that $v_{j,p} := \partial ^pu_j/\partial t^p$ for each
$p=1,...,n_j-1$, $w_{l,k} := \partial ^ku_l/\partial t^{k_0}\partial
x_1^{k_1}...\partial x_n^{k_n}$ for each $k_0+k_1+...+k_n=k\le n_l$
and $k_0<n_l$. Then $\partial v_{j,p}/\partial t=v_{j,p+1}$ for each
$p<n_j$ and $\partial v_{j,n_j-1}/\partial t=F_j$, also $w_{j,k}=
\partial ^kv_{l,k_0}/\partial x_1^{k_1}...\partial x_n^{k_n}$ for each
$(k_0,...,k_n)$ and $l$, where
$F_j=F_j(t,x_1,...,x_n,u_1,...,u_N,...,w_{j,k},...)$. The initial
conditions become:
\par $u_j|_{t=t_0}=\phi
_j^{(0)}(x_1,x_2,...,x_n)$,
\par $v_{j,s}|_{t=t_0}=\phi _j^{(s)}(x_1,x_2,...,x_n)$ \\
for each $s=1,...,n_j-1$ and each $(x_1,...,x_n)\in U$.
\par If $\partial ( v_{j,k-k_0} - \partial ^{k-k_0}u/\partial x_1^{k_1}
...\partial x_n^{k_n})/\partial t=0$ for each $t$ and
$(x_1,...,x_n)\in U$, then $( v_{j,k-k_0} -
\partial ^{k-k_0}u/\partial x_1^{k_1} ...\partial x_n^{k_n})$ is
independent of $t$. If in addition $v_{j,k-k_0}|_{t=t_0}=\phi _j^{
(k-k_0)}(x_1,x_2,...,x_n)$ for each $(x_1,...,x_n)\in U$, then
$\partial ^{k-k_0}u/\partial x_1^{k_1} ...\partial x_n^{k_n} = \phi
_j^{ (k-k_0)}(x_1,x_2,...,x_n)$ in $U$ due to analyticity of
functions.
\par Therefore, this procedure of the introduction of new functions
$v_{j,s}$ reduces the problem $1(1-3)$ to the problem, when $n_j=1$
for each $j$. It is possible to simplify further the problem, making
substitutions: $w_j(t,x) := u_j(t,x) - \phi _j^{(0)}(x)$,
$w_{j,s}(t,x) := v_{j,s}(t,x) - \phi _j^{(s)}(x)$ for each $j, s$,
where $x=(x_1,...,x_n)$. For new functions $w_j$ and $w_{j,s}$ the
problem has the same form, but with zero initial conditions. Thus
without loss of generality consider zero initial conditions and
denote $w_j$ and $w_{j,s}$ again by $u_j$ and $v_{j,s}$.
\par Consider at first the case when each $F_j$ is $\bf R$-homogeneous
and $\bf K$ additive by each $u_j$ and each $v_{j,s}$, where $n_j=1$
for each $j$. Without loss of generality consider $t\in \bf K$,
since when $t\in \bf R$, then due to the local analyticity by $t$ it
is possible to take the extension of each $F_j$ to $F_j\in L_q$ on a
neighborhood of $t_0$ in $\bf K$. In view of Proposition 2.18
\cite{ludfov} and Theorem 18 above if a solution of the problem
exists in the considered class of holomorphic functions $u_j$,
$j=1,...,N$, then it is unique.
\par Consider now series of $F_j$, $u_j$, $v_{j,s}$ for each $j$ and
$s$ by all their variables $(t-t_0, x_1-x_1^0,..., x_n-x_n^0)$ in a
neighborhood of $(t_0,x_1^0,...,x_n^0)$. Expansion coefficients of
$u_j$ and $v_{j,s}$ can be evaluated from the system of equations by
induction in the order of increasing powers, since $u_j(t_0,x)=0$
and $v_{j,s}(t_0,x)=0$ for each $j, s$ and each $x=(x_1,...,x_n)\in
U$. Now prove the local convergence of the series.
\par It is possible to make the shift of variables such that
consider $(t_0,x_1^0,...,x_n^0)=0$ without loss of generality. Let
$g(t,x) = \sum_{|k|\ge 0} \{ (a_k,(t,x)^k) \} _{q(2|k|)}$ be
converging at each point $b = (b_0,b_1,...,b_n)$ with given values
$|b_0|$,...,$|b_n|$, where $b_l\ne 0$ for each $l=0,1,...,n$. Then
there exists a constant $S>0$ such that $|\sum_{k=(k_0,...,k_n)} \{
(a_k,(t,x)^k) \} _{q(2|k|)}|\le S$ for each $k$, $|k|:=k_0+...+k_n$.
Therefore, $|\sum_{k=(k_0,...,k_n)} \{ (a_k,(\beta _0,\beta
_1,...,\beta _n)^k) \} _{q(2|k|)}|\le S |b_0|^{-k_0} ...
|b_n|^{-k_n}$ for each $\beta _0\in {\bf K},...,\beta _n\in {\bf K}$
with $|\beta _0|=1$,...,$|\beta _n|=1$ and each integers $0\le
k_0,...,k_n<\infty $, where vectors of associators $q(*)$ indicating
on orders of brackets or multiplications are important only over
$\bf O$, $a_k=a_k(g)$ (see Formulas 13(ii-iv)). Take the function
$Q(t,x) := S (1-t/|a_0|)^{-1}...(1-x_n/|a_n|)^{-1}$. Then $Q(t,x)$
is characterized by the condition: $$|\sum_{k=(k_0,...,k_n)} \{
(a_k(g),(\beta _0,\beta _1,...,\beta _n)^k) \} _{q(2|k|)}|\le
|\sum_{k=(k_0,...,k_n)} \{ (a_k(Q),(\beta _0,\beta _1,...,\beta
_n)^k) \} _{q(2|k|)}|$$ for each $k$ and each $\beta _0,...,\beta
_n$ in $\bf K$ with $|\beta _0|=1$,...,$|\beta _n|=1$. The
multiplications $z\mapsto \beta _j^pz$ and $z\mapsto z\beta _j^p$
are pseudoconformal mappings from $\bf K$ onto $\bf K$ for each
$p\in \bf N$. Such holomorphic function $Q$ is called the majorant
of the holomorphic function $g$. In the case of variables
$A_{j,s}\in L_q$ consider majorants in the class of analytic
functions by $A_{j,s}$ and their suitable compositions, where we put
$b_{j,s}\in \bf K$ such that $|b_{j,s}| := \| A_{j,s} \| $,  since
we can take $A_{j,s}.(h_1,...,h_s) $ and $F_j.h_0$ for $F_j\in
L_1({\bf K},{\bf K})$ for arbitrary $h_0, h_1,..., h_s\in \bf K$
with $|h_0|=1$, $|h_1|=1$,...,$|h_s|=1$ as well as $b_j$. Then such
composite functions while action on $h_0, h_1,..., h_s$ have also
series expansions which are uniformly convergent series by
$h_0,h_1,...,h_s\in B({\bf K},0,1)$.
\par Let problem I be the initial problem and problem II be a
problem, in which each participating function is substituted by its
majorant. Let $U_j$, $V_{j,s}$ denotes the solution of problem II.
Then $$|\sum_{k=(0,k_1,...,k_n)}\{ (a_k(v_{j,s}),(\beta _0,\beta
_1,...,\beta _n)^k) \} _{q(2|k|)}| \le |\sum_{k=(0,k_1,...,k_n)}\{
(a_k(V_{j,s}),(\beta _0,\beta _1,...,\beta _n)^k) \} _{q(2|k|)}|$$
for each $\beta _0,...,\beta _n$ in $\bf K$ with $|\beta
_0|=1$,...,$|\beta _n|=1$ and each $j$ and $s$ due to the initial
conditions, where $V_{j,0} := U_j$ and $v_{j,0} := u_j$. In the case
$k_0>0$ coefficients $a_k(v_{j,s})$ or $a_k(V_{j,s})$ are obtained
from coefficents $a_m(v_{j,s})$ or $a_m(V_{j,s})$ respectively with
the help of addition and multiplication with $0\le m_0<k_0$.
Therefore, by induction these inequalities are true for each $k$.
Thus if problem II has a solution, the problem I has a solution.
\par Choose constants $S>0$ and $a>0$ such that the function
$S[1-(t/\alpha +x_1+...+x_n)/a]^{-1}$ for some $0<\alpha <1$ would
be the majorant for all terms of the system besides free terms, as
well as take $T[1-(t/\alpha +x_1+...+x_n)/a]^{-1}$ as majorant for
free terms, since such majorants with $S=S_j$, $a=a_j$, $T=T_j$ can
be taken for each $F_j$ and hence we can take
$S=\max_{j=1,...,n}S_j$, $T:=\max_{j=1,...,n}T_j$,
$a=\min_{j=1,...,n}a_j$, where $\alpha $ we choose below. For it
problem II can be taken as:
\par $\partial U_j/\partial t = [\sum_{j,k} \partial U_j/\partial x_k +
\sum_jU_jI + CI]S[1-(t/\alpha +x_1+...+x_n)/a]^{-1}$, \\
where $C=T/S$, while $I\in L_1$ is the constant unit operator. Seek
a solution $U_j(t,x)=U(z)$, where $z=t/\alpha +x_1+...+x_n$. Then we
get the equation:
\par $\alpha ^{-1}dU(z)/dz =(NndU(z)/dz + NUI + CI)S(1-z/a)^{-1}$. \\
Choose a number $\alpha >0$ sufficiently large such that $1/\alpha
-NnS(1-z/a)^{-1}\ne 0$ in a neighborhood of $z=0$. Then the problem
is:
\par $(\alpha ^{-1} - S(1-z/a)^{-1}Nn) dU(z)/dz=I (NU+C)$, hence
\par $(NU+C)^{-1}(dU(z)/dz) = Ig(z)$, where $g(z) =
(\alpha ^{-1} - Nn S(1-z/a)^{-1})^{-1}$ due to the alternativity of
$\bf K$. The integration of the latter equation gives the
holomorphic solution due to Theorem 18, since $g$ is holomorphic.
Now estimate its expansion coefficients. For this consider an
algebraic embedding $\theta $ of the complex field $\bf C$ into $\bf
K$. The restriction of $U(z)$ on this copy of $\bf C$ has the form:
\par $U(z)=[ \exp ((\int_0^z y(\xi )d\xi )N/C)-1]C/N$,
for each $z\in \theta ({\bf C})$, \\
where $y(z) = (\alpha ^{-1} - Nn S(1-z/a)^{-1})^{-1}CS
(1-z/a)^{-1}$. The function $S (1-z/a)^{-1}$ has the nonnegative
expansion coefficients by powers $z^k$. Hence $y(\xi )$ also has the
nonegative expansion coefficients by $\xi $. Thus $s(z)=N\int_0^z
y(\xi )d\xi /C$ also has nonnegative expansion coefficients by $z$
and inevitably $\exp (s(z))-1=s(z)+s^2(z)/2+...$ and thus $U(z)$
also have nonnegative expansion coefficients by $z\in \theta ({\bf
C})$. But $\bf K$ is alternative and hence power associative,
moreover, $\bigcup_{\theta }\theta ({\bf C})=\bf K$. Thus expansion
coefficients of $U(x_1+x_2+...+x_n)$ by degress of $x_1$,...,$x_n$
are nonnegative, consequently, $U(0,x_1,...,x_n)$ is the majorant of
zero. Therefore, functions $U_j(t,x_1,...,x_n)=U(t/\alpha
+x_1+...+x_n)$ compose the solution of problem II.
\par The problem I solved above can be reformulated in such fashion, that
there exists an inverse operator $K^{-1}$ to the differential
operator together with initial conditions which compose the operator
equation $Ku=g$ in the space of holomorphic vector-functions $u, g$.
\par Consider now the case, when $F_j$ are not $\bf R$-homogeneous
$\bf K$-additive by $u_j$ or $v_{j,s}$ for some $j, s$.  Denote the
operator corresponding to this problem by $M$ and the problem is
$Mu=g$ in the space of holomorphic vector-functions. Let $X$ be a
$\bf K$-vector space. Then we call a subset $W$ in $X$ $\bf
K$-convex if $(a_1v_1)b_1 + a_2(v_2b_2)\in W$ for each $a_1, a_2,
b_1, b_2\in B({\bf K},0,1)$ and each $v_1, v_2\in W$. Since $F_j$
are holomorphic by $(t,x)$ and locally analytic by $u_j$, $v_{j,s}$
for each $j, s$, then we can take power series expansions by these
varaibles and consider a sufficiently small $\bf K$-convex
neighborhood of the point of zero initial conditions, where
corrections of the second and higher order terms in $u_j$, $v_{j,s}$
are much smaller. This neighborhood $W$ can be chosen such that
$M=K+(M-K)$ with $\| M-K \|_W /\| K \|_W <1/2$, where $\| M \|_W :=
\sup_{0\ne u\in W}|Mu|/|Ku| $. Hence $M$ is invertible on $W$ such
that there exists $M^{-1}: M(W)\to W$. The inverse of a holomorphic
function is holomorphic (see also Theorem 2). That is functions
defining $M^{-1}$ are holomorphic by $(x,u)$ also by $t$, when $t\in
\bf K$, moreover, they are locally analytic by others variables. In
view of the proof above we can construct a solution $u$ of
$Ku(t,x)=g$ for each $(t,x)$ with $|(t,x)|<R$ for some $0<R<\infty $
such that $u\in W/2$. Therefore, we get a solution of the problem
$Mu=g$ also for each $(t,x)$ in a sufficiently small neighborhood of
$0$. In view of the proof above this leads to the existence of a
solution of the general problem.
\par {\bf 5. Note.} The problem 3 is the noncommutative holomorphic
analog of the Cauchy problem, while Theorem 4 is the noncommutative
analog of the Kovalevsky theorem (see its commutative classical form
in \S 2 \cite{petrov}).
\par {\bf 6. Definitions.} A scalar product in a vector space $X$
over $\bf H$ (that is, linear relative to the right and left
multiplications separately on scalars from $\bf H$) it is a
biadditivite $\bf R$-bilinear mapping $<*;*>: X^2\to \bf H$, such
that
\par $(1)$ $<x;x>=\alpha _0$, where $\alpha _0\in \bf R$;
\par $(2)$ $<x;x>=0$ if and only if $x=0$;
\par $(3)$ $<x;y>=<y;x>^{\tilde .}$ for each $x, y \in X$;
\par $(4)$ $<x+z;y>=<x;y>+<z;y>$;
\par $(5)$ $<xa;yb>={\tilde a}<x;y>b$ for each
$x, y, z\in X$, $a, b\in \bf H$.
\par In the case of a vector space $X$ over $\bf O$ we consider an
$\bf O$-valued function on $X^2$ such that \par $(1')\quad <\zeta
,\zeta >=a$ with $a\ge 0$ and $<\zeta ,\zeta >=0$ if and only if
$\zeta =0$, \par $(2')\quad <\zeta ,z+\xi >=<\zeta ,z>+<\zeta ,\xi
>$, \par $(3')\quad <\zeta +\xi ,z>=<\zeta ,z>+<\xi ,z>$,
\par $(4')\quad <\alpha \zeta ,z>=\alpha <\zeta ,z>=<\zeta ,\alpha
z>$ for each $\alpha \in \bf R$ and $<\zeta a ,\zeta >={\tilde
a}<\zeta ,\zeta >$ for each $a\in \bf O$, \par $(5')\quad <\zeta
,z>^{\tilde .}=<z,\zeta >$ for each $\zeta , \xi $ and $z\in X$.
\par While the representation of $X$ in the form
$X=X_0i_0\oplus X_1i_i\oplus ... \oplus X_mi_m$, where $m=3$ for
$\bf H$ and $m=7$ for $\bf O$, $X_0, ...,X_m$ are pairwise
isomorphic $\bf R$-linear spaces, $ \{ i_0,...,i_m \} $ is the
family of standard generators of the algebra ${\bf K}=\bf H$ or
${\bf K}=\bf O$, $i_0=1$, we shall suppose naturally, that
\par $(6)$ $<x_p,y_q>\in \bf R$ and
\par $(7)$ $<x_pi_p,y_qi_q>=<x_p,y_q>i_p^*i_q$
for each $x_p\in X_p, y_q\in X_q$, $p, q\in \{ 0,1,2,...,m \} $.
\par These scalar products are called also $\bf K$-Hermitian inner
products.
\par If $X$ is complete relative to the norm topology
\par $(8)$ $|x|:=<x;x>^{1/2}$,  then $X$ is called $\bf K$ Hilbert space.
\par In particular, for $X=\bf K^n$ we can take the canonical scalar
product: \\
$(9)\quad <\zeta ;z>:= (\zeta ,z)=\sum_{l=1}^n\mbox{ }^l{\tilde
\zeta }\mbox{ }^lz$, where $z=(\mbox{ }^1z,...,\mbox{ }^nz)$,
$\mbox{ }^lz\in \bf K$.
\par Let $l_2 ({\bf K})$ be the $\bf K$-vector space consisting
of all sequences $x = \{ x_n\in {\bf K}: n \in {\bf N} \} $ with the
finite norm $\| x \| <\infty $ and scalar product
$(x,y):=\sum_{n=1}^{\infty } x_ny_n^*$ with $\| x \| := (x,x)^{1/2}$
is called the Hilbert space (of separable type) over $\bf K$, where
$z^*$ denotes the conjugated $z$, $zz^* =: |z|^2$, $z\in \bf K$,
${\bf K}=\bf H$ or ${\bf K}=\bf O$.
\par Henceforth, the term a $\bf K$ manifold $N$
modelled on $X=\bf K^n$ or $X=l_2({\bf K})$ means a metric separable
space supplied with an atlas $\{ (U_j,\phi _j): j\in \Lambda _N \} $
such that:
\par $(i)$ $U_j$ is an open subset of $N$ for each $j\in \Lambda _N$
and $\bigcup_{j\in \Lambda _N}U_j=N$, where $\Lambda _N\subset \bf
N$;
\par $(ii)$ $\phi _j: U_j\to \phi _j(U_j)\subset X$ is a
$C^{\infty }$-diffeomorphism for each $j$, where $\phi _j(U_j)$ is a
$C^{\infty }$-domain in $X$;
\par $(iii)$ $\phi _j\circ \phi _m^{-1}$ is a $\bf K$ biholomorphic
mapping from $\phi _m(U_m\cap U_j)$ onto $\phi _j(U_m\cap U_j)$
while $U_m\cap U_j\ne \emptyset $. When $X=l_2({\bf K})$ it is
supposed, that $\phi _j\circ \phi _m^{-1}$ is Frech\'et (strongly)
$C^{\infty }$-differentiable for each $j$ and $m$ and certainly
holomorphic by each $\bf K$ variable.
\par A family $ \{ A_s: s \in S \} $ of subsets $A_s$ of a
topological space $X$ is called locally finite, if for each point
$x\in X$ there exists an open neighborhood $U$ of $x$ such that the
set $\{ s: A_s\cap U\ne \emptyset \} $ is finite, where $S$ is a
set.
\par A topological space $X$ is called paracompact, if each its
open covering can be refined by a locally finite open covering.
\par {\bf 7. Proposition.} {\it Let $M$ be a $\bf K$ holomorphic
manifold, where ${\bf K}=\bf H$ or ${\bf K}=\bf O$. Then there
exists a tangent bundle $TM$ which has the structure of the $\bf K$
holomorphic manifold such that each fibre $T_xM$ is the vector space
over $\bf K$.}
\par {\bf Proof.} Since $\bf K$ is the algebra over $\bf R$, then
it has the real shadow, which is the Euclidean space $\bf R^{2^r}$,
where $r=2$ for quaternions, $r=3$ for octonions. Therefore, $M$ has
also the structure of real manifold. Since each $\bf K$ holomorphic
mapping is infinite differentiable (see Theorems 2.15 and 3.10
\cite{ludoyst,ludfov}), then there exists its tangent bundle $TM$
which is $C^{\infty }$-manifold such that each fibre $T_xM$ is a
tangent space, where $x\in M$, $T$ is the tangent functor. If $At
(M)= \{ (U_j,\phi _j): j \} $, then $At (TM)= \{ (TU_j, T\phi _j): j
\} $, $TU_j=U_j\times X$, where $X$ is the $\bf K$ vector space on
which $M$ is modelled, $T(\phi _j\circ \phi _k^{-1})=(\phi _j\circ
\phi _k^{-1}, D(\phi _j\circ \phi _k^{-1}))$ for each $U_j\cap
U_k\ne \emptyset $. Each transition mapping $\phi _j\circ \phi
_k^{-1}$ is $\bf K$ holomorphic on its domain, then its (strong)
differential coincides with the superdifferential $D(\phi _j\circ
\phi _k^{-1})= D_z(\phi _j\circ \phi _k^{-1})$, since ${\tilde
\partial } (\phi _j\circ \phi _k^{-1})=0$. Therefore, $D(\phi
_j\circ \phi _k^{-1})$ is $\bf R$-linear and $\bf K$-additive, hence
is the automorphism of the $\bf K$ vector space $X$. Since $D_z(\phi
_j\circ \phi _k^{-1})$ is also $\bf K$ holomorphic, then $TM$ is the
$\bf K$ holomorphic manifold.
\par {\bf 8. Definitions.} A $C^1$-mapping $f: M\to N$ is called
an immersion, if $rang (df|_x: T_xM\to T_{f(x)}N) = m_M$ for each
$x\in M$, where $m_M := dim_{\bf R}M$. An immersion $f: M\to N$ is
called an embedding, if $f$ is bijective.
\par {\bf 9. Theorem.} {\it Let $M$ be a compact
$\bf K$ holomorphic manifold, $dim_{\bf K}M=m<\infty $, where ${\bf
K}=\bf H$ or ${\bf K}=\bf O$. Then there exists a $\bf K$
holomorphic embedding $\tau : M\hookrightarrow {\bf K}^{2m+1}$ and a
$\bf K$ holomorphic immersion $\theta : M\to {\bf K}^{2m}$
correspondingly. Each continuous mapping $f: M\to {\bf K}^{2m+1}$ or
$f: M\to {\bf K}^{2m}$ can be approximated by $\tau $ or $\theta $
relative to the norm $\| * \|_{C^0}$. If $M$ is a paracompact $\bf
K$ holomorphic manifold with countable atlas on $l_2({\bf K})$, then
there exists a holomorphic embedding $\tau : M\hookrightarrow
l_2({\bf K})$.}
\par {\bf Proof.} Let at first $M$ be compact.
Since $M$ is compact, then it is finite dimensional over $\bf K$,
$dim_{\bf K}M=m\in \bf N$, such that $dim_{\bf R}M=2^rm$ is its real
dimension. Choose an atlas $At' (M)$ refining initial atlas $At (M)$
of $M$ such that $({U'}_j, \phi _j)$ are charts of $M$, where each
${U'}_j$ is $\bf K$ holomorphic diffeomorphic to an interior of the
unit ball $Int (B({\bf K^m},0,1))$, where $B({\bf K^m},y,r) := \{
z\in {\bf K^m}: |z-y|\le r \} $ (see Theorem 47). In view of
compactness of $M$ a covering $\{ {U'}_j, j \} $ has a finite
subcovering, hence $At' (M)$ can be chosen finite. Denote for
convenience the latter atlas as $At (M)$. Let $(U_j, \phi _j)$ be
the chart of the atlas $At (M)$, where $U_j$ is open in $M$, hence
$M\setminus U_j$ is closed in $M$.
\par  Consider the space ${\bf K^m}\times {\bf R}$
as the $\bf R$-linear space $\bf R^{2^rm+1}$. The unit sphere
$S^{2^rm}:=S ({\bf R}^{2^rm+1},0,1) := \{ z\in {\bf R}^{2^rm+1}:$
$|z|=1 \} $ in ${\bf K^m}\times \bf R$ can be supplied with two
charts $(V_1, \phi _1)$ and $(V_2, \phi _2)$ such that $V_1:=
S^{2^rm}\setminus \{ 0,...,0, 1 \} $ and $V_2:=S^{2^rm}\setminus \{
0,...,0, - 1 \} $, where $\phi _1$ and $\phi _2$ are stereographic
projections from poles $\{ 0,...,0, 1 \} $ and $ \{ 0,...,0, -1 \} $
of $V_1$ and $V_2$ respectively onto $\bf K^m$. Since $z^* = -
(2^r-2)^{-1} \sum_{p=0}^{2^r-1}(i_pz)i_p $ in $\bf K^m$, then $\phi
_1\circ \phi _2^{-1}$ (in the $z$-representation) is the $\bf K$
holomorphic diffeomorphism of ${\bf K^m}\setminus \{ 0 \} $, but
certainly with neither right nor left superlinear superdifferential
$D_z (\phi _1\circ \phi _2^{-1})$, where $i_0,...,i_{2^r-1}$ are the
standard generators of $\bf K$. Thus $S^{2^rm}$ can be supplied with
the structure of the $\bf K$ holomorphic manifold.
\par Therefore, there exists a $\bf K$ holomorphic mapping
$\psi _j$ (that is, locally $z$-analytic \cite{ludoyst}) of $M$ into
the unit sphere $S^{2^rm}$ such that $\psi _j: U_j\to \psi _j(U_j)$
is the $\bf K$ holomorphic diffeomorphism onto the subset $\psi
_j(U_j)$ in $S^{2^rm}$. The using of such $\psi _j$ is sufficient,
where $\psi _j$ can be considered as components of a holomorphic
diffeomorphism: $\psi : M\to (S^{2^rm})^n$ with $n$ equal to number
of charts. In the case of $\partial M\ne \emptyset $ and $U_j\cap
\partial M\ne \emptyset $ an image $\psi _j(U_j)$
has a nonvoid boundary and Theorem 2.39 can be used. There is
evident embedding of ${\bf K^m}\times \bf R$ into $\bf K^{m+1}$.
Then the mapping $\psi (z):=(\psi _1(z),...,\psi _n(z))$ is the
embedding into $(S^{2^rm})^n$ and hence into ${\bf K}^{n{m+1}}$,
since the rank $rank [d_z\psi (z)]=2^rm$ at each point $z\in M$,
because $rank [d_z\psi _j(z)]=2^rm$ for each $z\in U_j$ and
$dim_{\bf K}\psi (U_j)\le dim_{\bf K}M=m$. Moreover, $\psi (z)\ne
\psi (y)$ for each $z\ne y\in U_j$, since $\psi _j(z)\ne \psi
_j(y)$. If $z\in U_j$ and $y\in M\setminus U_j$, then there exists
$l\ne j$ such that $y\in U_l\setminus U_j$, $\psi _j(z)\ne \psi
_j(y)=x_j$.
\par Let $M\hookrightarrow {\bf K}^N$ be the $\bf K$
holomorphic embedding as above. There is also the $\bf K$
holomorphic embedding of $M$ into $(S^{2^rm})^n$ as it is shown
above, where $(S^{2^rm})^n$ is the $\bf K$ holomorphic manifold as
the product of $\bf K$ holomorphic manifolds. Consider the bundle of
all $\bf K$ straight lines ${\bf K}x$ in ${\bf K}^N$, where $x\in
{\bf K}^N$, $x\ne 0$, ${\bf K}x$ is the $\bf K$ vector space of
dimension $1$ over $\bf K$, which has the real shadow isomorphic
with $\bf R^{2^r}$. By our definition they compose the
noncommutative analog of the projective space ${\bf K}P^{N-1}$,
since $\bf K$ is alternative. Fix the standard orthonormal base $ \{
e_1,...,e_N \} $ in ${\bf K}^N$ and projections on $\bf K$-vector
subspaces relative to this base $P^L(x):=\sum_{e_j\in L}x_je_j$ for
the $\bf K$-linear span $L=span_{\bf K} \{ e_i:$ $i\in \Lambda _L \}
$, $\Lambda _L\subset \{ 1,...,N \} $, where $x=\sum_{j=1}^Nx_je_j$,
$x_j\in \bf K$ for each $j$, $e_j=(0,...,0,1,0,...,0)$ with $1$ at
$j$-th place. In this base consider the $\bf K$-Hermitian scalar
product $<x,y> := \sum_{j=1}^Nx_jy^*_j$. Let $l\in {\bf K}P^{N-1}$,
take a $\bf K$-hyperplane denoted by ${\bf K}^{N-1}_l$ and given by
the condition: $<x,[l]>=0$ for each $x\in {\bf K}^{N-1}_l$, where
$0\ne [l]\in {\bf K}^N$ characterises $l$. Take $\| [l] \| =1$. Then
the orthonormal base $\{ q_1,...,q_{N-1} \} $ in ${\bf K}^{N-1}_l$
and together with $[l]=:q_N$ composes the orthonormal base $\{
q_1,...,q_N \} $ in ${\bf K}^N$. This provides the $\bf K$
holomorphic projection $\pi _l: {\bf K}^N\to {\bf K}^{N-1}_l$
relative to the orthonormal base $ \{ q_1,...,q_N \} $. The operator
$\pi _l$ is $\bf K$ left and also right linear (but certainly
nonlinear relative to $\bf K$), hence $\pi _l$ is $\bf K$
holomorphic.
\par To construct an immersion it is sufficient, that each
projection $\pi _l: T_xM\to {\bf K}^{N-1}_l$ has $ker [d(\pi
_l(x))]= \{ 0 \} $ for each $x\in M$. The set of all $x\in M$ for
which $ker [d(\pi _l(x))] \ne \{ 0 \} $ is called the set of
forbidden directions of the first kind. Forbidden are those and only
those directions $l\in {\bf K}P^{N-1}$ for which there exists $x\in
M$ such that $l'\subset T_xM$, where $l'=[l]+z$, $z\in {\bf K}^N$.
The set of all forbidden directions of the first kind forms the $\bf
K$ holomorphic manifold $Q$ of $\bf K$ dimension $(2m-1)$ with
points $(x,l)$, $x\in M$, $l\in {\bf K}P^{N-1}$, $[l]\in T_xM$. Take
$g: Q\to {\bf K}P^{N-1}$ given by $g(x,l):=l$. Then $g$ is $\bf K$
holomorphic.
\par Each paracompact manifold $A$ modelled
on $\bf K^p$ can be supplied with the Riemann manifold structure
also. Therefore, on $A$ there exists a Riemann volume element. In
view of the Morse theorem $\mu (g(Q))=0$, if $N-1>2m-1$, that is,
$2m<N$, where $\mu $ is the Riemann volume element in ${\bf
K}P^{N-1}$. In particular, $g(Q)$ is not contained in ${\bf
K}P^{N-1}$ and there exists $l_0\notin g(Q)$, consequently, there
exists $\pi _{l_0}: M\to {\bf K}^{N-1}_{l_0}$. This procedure can be
prolonged, when $2m<N-k$, where $k$ is the number of the step of
projection. Hence $M$ can be immersed into ${\bf K}^{2m}$.
\par Consider now the forbidden directions of the second type:
$l\in {\bf K}P^{N-1}$, for which there exist $x\ne y\in M$
simultaneously belonging to $l$ after suitable parrallel translation
$[l]\mapsto [l]+z$, $z\in {\bf K}^N$. The set of the forbidden
directions of the second type forms the manifold $\Phi
:=M^2\setminus \Delta $, where $\Delta := \{ (x,x):$ $x\in M \} $.
Consider $\psi : \Phi \to {\bf K}P^{N-1}$, where $\psi (x,y)$ is the
straight $\bf K$-line with the direction vector $[x,y]$ in the
orthonormal base. Then $\mu (\psi (\Phi ))=0$ in ${\bf K}P^{N-1}$,
if $2m+1<N$. Then the closure $cl (\psi (\Phi ))$ coinsides with
$\psi (\Phi )\cup g(Q)$ in ${\bf K}P^{N-1}$. Hence there exists
$l_0\notin cl (\psi (\Phi ))$. Then consider $\pi _{l_0}: M\to {\bf
K}_{l_0}^{N-1}$. This procedure can be prolonged, when $2m+1<N-k$,
where $k$ is the number of the step of projection. Hence $M$ can be
embedded into ${\bf K}^{2m+1}$.
\par The approximation property follows from compactness of $M$
and the noncommutative analog of the Stone-Weierstrass theorem (see
also Theorem 2.7 \cite{ludoyst,ludfov}.)
\par Let now $M$ be a paracompact $\bf K$ holomorphic manifold
with countable atlas on $l_2({\bf K})$. Spaces $l_2({\bf K})\oplus
{\bf K}^m$ and $l_2({\bf K})\oplus l_2({\bf K})$ are isomorphic as
$\bf K$ Hilbert spaces with $l_2({\bf K})$. Take an additional
variable $z\in \bf K$, when $z=j\in \bf N$, then it gives a number
of a chart. Each $TU_j$ is $\bf K$ holomorphically diffeomorphic
with $U_j\times l_2({\bf K})$. Consider $\bf K$ holomorphic
functions $\psi $ on domains in $l_2({\bf K})\oplus l_2({\bf
K})\oplus \bf K$. Then there exists a $\bf K$ holomorphic mapping
$\psi _j: M\to l_2({\bf K})$ such that $\psi _j: U_j\to \psi
_j(U_j)\subset l_2({\bf K})$ is a $\bf K$ holomorphic
diffeomorphism. Then the mapping $(\psi _1,\psi _2,...)$ provides
the $\bf K$ holomorphic embedding of $M$ into $l_2({\bf K})$.
\par {\bf 10. Theorem.} {\it Let $M$ be a $\bf K$ pseudoconformal
locally compact manifold with and atlas $At (M)= \{ (U_j,\phi _j) :
j \} $ such that transition mappings $\phi _j\circ \phi _k^{-1}$ of
charts with $U_j\cap U_k\ne \emptyset $ are pseudoconformal. Then
$M$ is orientable.}
\par {\bf Proof.} Since $M$ is locally compact, then $M$ is finite
dimensional over $\bf H$ such that $dim_{\bf K}X=m<\infty $, where
$X=T_xM$ for $x\in M$. Consider open subsets $U$ and $V$ in $X$ and
a function $\phi : U\to V$ which is pseudoconformal. Write $\phi $
in the form $\phi = (\mbox{ }_1\phi ,..., \mbox{ }_m\phi )$, where
$\mbox{ }_v\phi \in \bf K$ for each $v=1,...,m$. Use Theorems 2.4
and 2.5 to present $\phi '$ as compositions of left and right $\bf
K$ linear operators.
\par Consider a right superlinear operator $A$
on $X$, then $A(h) =A(e_k)h$ for each $k=1,...,m$ and each $h\in
{\bf K}$, where $e_k=(0,...,0,1,0,...,0)\in \bf K^m$ with $1$ on the
$k$-th place. If $A$ is left superlinear, then $A(h)=hA(e_k)$ for
each $k=1,...,m$ and each $h\in {\bf K}$. Moreover, left and right
superlinear operators are certainly $\bf R$-linear and $\bf
K$-additive.
\par Consider now real local coordinates in $M$:
$ \{ \mbox{ }_{v,l}\phi : v= 1,...,m; l= 0,1,2,2^r-1 \} $ such that
$\mbox{ }_v\phi = \sum_{l=0}^{2^r-1} \mbox{ }_{v,l} \phi i_l $, then
the determinant of change of local coordinates between charts is
positive: $det ( \{
\partial \mbox{ }_{v,l}\phi /
\partial \mbox{ }_{w,p}z \}_{v,l;w,p} )>0$, since
$det (A)>0$ for each nondegenerate left or right superlinear
operator $A$, because ${\bf H}={\bf C}\oplus {\bf C}i_2$, ${\bf
O}={\bf H}\oplus {\bf H}i_4$, where $\{ i_0,i_1,...,i_{2^r-1} \} $
are standard generators of $\bf K$ (see also Theorem 2.1.3.7
\cite{lusmldg05}).
\par {\bf 11. Proposition.} {\it Let $M$ be a $\bf K$
paracompact holomorphic manifold on $X=\bf K^n$ or on $X=l_2({\bf
K})$, where $n\in \bf N$, ${\bf K}=\bf H$ or ${\bf K}=\bf O$. Then
$M$ can be supplied with a $\bf K$ holomorphic Hermitian metric.}
\par {\bf Proof.} Mention that a continuous mapping $\pi $
of a Hausdorff space $A$ into another $B$ is called an $\bf
K$-vector bundle, if
\par $(i)$ $A_x:=\pi ^{-1}(x)$ is a vector space $X$ over $\bf K$
for each $x\in B$ and
\par $(ii)$ there exists a neighbourhood $U$ of $x$ and a homeomorphism
$\psi : \pi ^{-1} (U)\to U\times X$ such that $\psi (A_x)= \{ x \}
\times X$ and $\psi _x: A_x\to X$ is a $\bf K$-vector isomorphism
for each $x\in B$, where $\psi ^x:=\psi \circ \pi _2$, $\pi _2:
U\times X\to X$ is the projection, $(U,\psi )$ is called a local
trivialization.
\par If $\pi : A\to B$ is an $\bf K$-vector bundle, then
an $\bf K$-Hermitian metric $g$ on $A$ is an assignment of an $\bf
K$-Hermitian inner product $<*,*>_x$ to each fibre $A_x$ such that
for each open subset $U$ in $B$ and $\xi ,\eta \in C^{\infty }(U,A)$
the mapping $<\xi ,\eta >: U\to \bf K$ such that $<\xi , \eta
>(x)=<\xi (x), \eta (x)>_x$ is $C^{\infty }$. An $\bf K$-vector
bundle $A$ equipped with an $\bf K$-Hermitian metric $g$ we call an
$\bf K$-Hermitian vector bundle. If $A$ is paracompact, then each
its (open) covering contains a locally finite refinement (see
\cite{eng}). Take a locally finite covering $\{ U_j: j \} $ of $B$.
By the supposition of this proposition $X=\bf K^n$ or $X=l_2({\bf
K})$. Choose a subordinated real partition of unity $\{ \alpha _j: j
\} $ of class $C^{\infty }$: $\quad \sum_j\alpha _j(x)=1$ and
$\alpha _j(x)\ge 0$ for each $j$ and each $x\in B$. Therefore, there
exists a frame $\{ e_k: k \in \Lambda \} $ at $x\in B$, where either
$\Lambda = \{ 1,2,...,n \} $ or $\Lambda =\bf N$, such that $e_k\in
C^{\infty }(U_j,A)$ for each $k$, $\{ e_k: k \in \Lambda \} $ are
$\bf K$-linearly independent at each $y\in U_j$ relative to left and
right mulitplications on constants $a_k, b_k$ from $\bf K$, that is,
$\sum_k a_ke_kb_k=0$ if and only if $\sum_k |a_kb_k|=0$, where $x\in
U_j$. Define $<\xi ,\eta >^j_x:=\sum_k \xi (e_k)(x) (\eta
(e_k)(x))^*$ and $<\xi ,\eta >_x := \sum_j \alpha _j(x) <\xi ,\eta
>^j_x$. If $\xi ,\eta \in C^{\infty }(U,A)$, then the mapping
$x\mapsto <\xi (x), \eta (x)>_x$ is $C^{\infty }$ on $U$. Since
$|<\xi ,\eta >_x|\le \sum_j \alpha _j(x) (<\xi ,\xi >^j_x)^{1/2}
(<\eta ,\eta >^j_x)^{1/2}<\infty $ for each $x\in B$, hence the $\bf
K$-Hermitian metric is correctly defined. For each $b\in \bf R$ we
have $<\xi b, \eta b>_x=|b|^2<\xi ,\eta >_x$, since $\bf R$ is the
centre of the algebra $\bf K$ over $\bf R$. Take in particular $\pi
: TM\to M$ and this provides an $\bf K$-Hermitian metric in $M$.
\par In view of Theorem 9 there exists a $\bf K$ holomorphic embedding
$\tau $ of $M$ into the corresponding $\bf K$ vector space $X$
either finite dimensional over $\bf K$ or isomorphic with $l_2({\bf
K})$. Take in $X$ a $\bf K$ Hermitian metric and consider its
restriction on $\tau (M)$. Since $\tau $ is the holomorphic
diffeomorphism of $M$ onto $\tau (M)$, then the metric $g$ in $M$
can be chosen $\bf K$ holomorphic and Hermitian, that is, $<\xi ,
\eta >(x)$ is holomorphic by the $x$ variable.
\par {\bf 12. Theorem.} {\it Let $M$ be a $\bf K$ holomorphic
manifold, then there exists an open neighborhood ${\tilde T}M$ of
$M$ in $TM$ and an exponential $\bf K$ holomorphic mapping $\exp :
{\tilde T}M\to M$ of ${\tilde T}M$ on $M$.}
\par {\bf Proof.} It was shown in \S \S 10 and 11
that each $\bf K$ holomorphic manifold has also the structure of the
Riemann manifold. Therefore, the geodesic equation
\par  $(i)$  $\nabla _{\dot c}{\dot c}=0$ with initial conditions
$c(0)=x_0$, ${\dot c}(0)=y_0$, $x_0\in M$, $y_0\in T_{x_0}M$ \\
has a unique $C^{\infty }$-solution, $c: (- \epsilon , \epsilon )
\to M$ for some $\epsilon >0$. For a chart $(U_j, \phi _j)$
containing $x$, put $\psi _j(\beta ) = \phi _j\circ c(\beta )$,
where $\beta \in ( - \epsilon , \epsilon )$. Consider an $\bf
K$-Hermitian metric $g$ in $M$, then $g$ is $\bf K$ holomorphic in
local $\bf K$ coordinates in $M$ in accordance with Proposition 11,
$\partial g(z)/\partial {\tilde z}=0$, where $g(z)(*,*)=<*,*>_z$ is
the $\bf K$-Hermitian inner product in $T_zM$, $z\in M$, where $\phi
_l(z)\in X$ is denoted for convenience also by $z$. Consider
real-valued inner product induced by $g$ of the form $G(z)(x,y) :=
(g(z)(x,y) + (2^r-2)^{-1} \{ - g(z)(x,y) + \sum_{l=1}^{2^r-1} i_l
(g(z)(x,y) i_l^*) \} )/2$. Then $G(z)$ is also $\bf K$ holomorphic:
$\partial G(z)(x,y)/\partial {\tilde z}=0$ for each $x, y$. In real
local coordinates $G$ can be written as: $G(z)(\eta ,\xi
)=\sum_{l,p} G^{l,p}(z)\eta _l \xi _p$, where $x = (\mbox{
}_1x,...,\mbox{ }_mx)$, $\mbox{ }_vx=\sum_{l=2^r(v-1)+1}^{2^rv} \eta
_l i_{l-2^r(v-1)-1}$ and $\mbox{ }_vy=\sum_{l=2^r(v-1)+1}^{2^rv} \xi
_l i_{l-2^r(v-1)-1}$ for each $v=1,...,m$. Then the Christoffel
symbols are given by the equation $\Gamma ^a_{b,c} =
(\sum_{l=1}^{2^rm} G^{a,l} (\partial _bG_{c,l} +
\partial _c G_{l,b} - \partial _l G_{b,c})/2$ for each $a, b,
c=1,...,2^rm$. Using expression of the Christoffel symbol $\Gamma $
through $g$, we get that $\Gamma (z)$ is $\bf K$ holomorphic:
\par  $(ii)$  $\partial  \Gamma ^a_{b,c}(z)/\partial {\tilde z}=0$
for each $a, b, c=1,...,2^rm$. \\
Thus differential operator corresponding to the geodesic equation
has the $\bf K$ holomorphic form:
\par $(iii)$ $d^2 c(t)/dt^2 + \Gamma (c(t)) (c'(t), c'(t))=0$. \\
In view of Theorem 4 the mapping ${\tilde T}V_1\times (-\epsilon ,
\epsilon ) \ni (z_0,y_0;\beta )\mapsto \psi _j(\beta ;x_0,y_0)$ is
$\bf K$ holomorphic by $(x_0,y_0)$, since components of $y_0$ can be
expressed through $\bf R$-linear combinations of $\{ i_ly_0i_l: l=0,
1, 2, ...,2^r-1 \} $, where $0<\epsilon $, $z_0=\phi _j(x_0)\in
V_1\subset V_2\subset \phi _j(U_j)$, $V_1$ and $V_2$ are open,
$\epsilon $ and ${\tilde T}V_1$ are sufficiently small, that to
satisfy the inclusion $\psi _j(\beta ;x_0,y_0)\in V_2$ for each
$(z_0,y_0;\beta )\in {\tilde T}V_1\times (-\epsilon , \epsilon )$.
\par Then there exists $\delta >0$ such that
$c_{aS}(t)=c_S(at)$ for each $a\in (- \epsilon , \epsilon )$ with
$|aS(\phi _j(q))|<\delta $ since
$dc_S(at)/dt=a(dc_S(z)/dz)|_{z=at}$. The projection $\tau :=\tau _M:
TM\to M$ is given by $\tau _M(s)=x$ for each vector $s\in T_xM$,
$\tau _M$ is the tangent bundle. For each $x_0\in M$ there exists a
chart $(U_j,\phi _j)$ and open neighbourhoods $V_1$ and $V_2$, $\phi
_j(x_0)\in V_1\subset V_2\subset \phi _j(U_j)$ and $\delta >0$ such
that from $S\in TM$ with $\tau _MS=q\in \phi _j^{-1}(V_1)$ and
$|S(\phi _j(q))|<\delta $ it follows, that the geodesic $c_S$ with
$c_S(0)=S$ is defined for each $t\in (-\epsilon ,\epsilon )$ and
$c_S(t)\in \phi _j^{-1}(V_2)$. Due to paracompactness of $TM$ and
$M$ this covering can be chosen locally finite \cite{eng}.
\par This means that there exists an open neighbourhood
${\tilde T}M$ of $M$ in $TM$ such that a geodesic $c_S(t)$ is
defined for each $S\in {\tilde T}M$ and each $t\in (-\epsilon ,
\epsilon )$. Therefore, define the exponential mapping $\exp :
{\tilde T}M\to M$ by $S\mapsto c_S(1)$, denote by $\exp _x:=\exp
|_{{\tilde T}M\cap T_xM}$ a restriction to a fibre. Then $\exp $ has
a local representation $(x_0,y_0)\in {\tilde T}V_1\mapsto \psi
_j(1;x_0,y_0)\in V_2\subset \phi _j(U_j).$ From Equations $(i -
iii)$ it follows that $\exp $ is $\bf K$ holomorphic from ${\tilde
T}M$ onto $M$.
\par {\bf 13. Theorem.} {\it Let $M$ be a compact holomorphic
manifold on $\bf K^n$ without boundary, where $0<n\in \bf N$. \par
$(i)$ If $v$ is a holomorphic vector field on $M$, its flow $\eta
_t$ is a one parameter subgroup of $DifH(M)$, where $t\in \bf R$.
\par $(ii)$ The curve $t\mapsto \eta _t$ is locally analytic.
\par $(iii)$ The map $Exp: T_{id}DifH(M)\to DifH(M)$, $v\mapsto \eta _1$
is holomorphic.}
\par {\bf Proof.} Remind, that a holomorphic vector field on a
holomorphic manifold $M$ is a holomorphic mapping $v: M\to TM$ such
that $\pi \circ v=id$, where $TM$ denotes the tangent bundle, $\pi :
TM\to M$ is the natural projection, $\pi (x,y)=x$ for each $x\in M$
and $y\in T_xM$, $id(x)=x$ for each $x\in M$. Consider a flow $\eta
_t$ of a vector field $v$, which is characterized by the Equations
$(iv-vi)$: \par $(iv)$ $\eta _0=id$, \par $(v)$ $\eta _t\circ \eta
_p=\eta _{t+p}$ for each $t, p$, \par $(vi)$ $(\partial \eta
_t/\partial t)|_{t=0}=V$, where $t, p \in \bf R$. \\
Differentiating Equation $(v)$ by $t$ gives: \par $(vii)$ $\partial
\eta _p/\partial p= v\circ \eta _p$ for each $p$. \\
This is the autonomous Cauchy problem. Therefore, its solution
satisfies Condition $(v)$. In view of Theorem 4 if $v$ is
holomorphic, then its flow $\eta _t$ is holomorphic by $x\in M$ and
locally analytic by $t$. \par This solution is at first constructed
for each $|t|<\epsilon $, where $\epsilon >0$ is some positive
number. For sufficiently small $\epsilon >0$ each $\eta _t$ belongs
to $DifH(M)$. For general $t$ write for each $m\in \bf N$: $\eta _t=
(\eta _{t/m})^m = \eta _{t/m}\circ ... \circ \eta _{t/m}$. For large
$m$ the number $t/m$ is sufficiently small $|t/m|<\epsilon $ and
hence $\eta _{t/m}\in DifH(M)$. Since $DifH(M)$ is the group, then
$\eta _t\in DifH(M)$. \par The curve $\eta _t$ is analytic for
sufficiently small $t$. On the other hand, $\eta _t=\eta _{t-\theta
}\circ \eta _{\tau }$, hence the curve $t\mapsto \eta _t$ is locally
analytic by $t$, since the composition in $DifH(M)$ is holomorphic.
Then Equation $(v)$ also shows, that $Exp$ is holomorphic.
\par {\bf 14. Definition.} Let $M$ be two $\bf K$ holomorphic manifolds
on a $\bf K$ vector $X$ space either $\bf K^m$ or $l_2({\bf K})$
such that $N$ is a submanifold in $M$. Then $N$ is foliated relative
to $M$ if there exists a $\bf K$ holomorphic atlas $At(N) = \{
(V_k,\psi _k): k \} $ of $N$ and if there is a $\bf K$ holomorphic
atlas $At (M)= \{ (U_j,\phi _j): j \} $ of $M$ such that the
restriction of $M$ on $N$ gives $At (N)$ and
\par $(i)$ $\phi _{i,j}: \phi _j(U_i\cap U_j)\to X$
are of the form $\phi _{i,j}(x,y)=(\alpha _{i,j}(x),\beta
_{i,j}(y))$, where $(x,y)$ are local $\bf K$ coordinates in $M$, $x$
are local $\bf K$ coordinates in $N$, $\alpha _{i,j}(x)\in X_1$,
$\beta _{i,j}(y)\in X_2$, $X_1\oplus X_2=X$. If $\partial M\ne
\emptyset $ then in addition:
\par $(ii)$ $\partial N\subset \partial M$ and for each boundary
component $M_0$ of $M$ and $U_i\cap M_0\ne \emptyset $ we have $\phi
_i: (U_i\cap M_0)\to Y$ and $\phi _i: (U_i\cap M_0\cap \partial
N)\to Y_1$, where $Y := \{ x\in X: Re (x_1)=0 \} $, $Y_1 = \{ x\in
X_1: Re (x_1)=0 \} $, $x=(x_1,x_2,...)\in X$, $x_j\in \bf K$ for
each $j$.
\par {\bf 15. Theorem.} {\it Let $M$ be a compact holomorphic manifold
on $\bf K^n$ without boundary and $N\subset M$ be its closed
holomorphic submanifold on $\bf K^b$ with $0\le b\le n$ (possibly
zero dimensional) without boundary such that $N$ is foliated
relative to $M$, where $0<n\in \bf N$. Let also $DifH_N (M) := \{
f\in DifH (M): f(N)\subset N \} $ and $DifH_{N,p}(M) := \{ f\in
DifH(M): f(x)=x \forall x\in N \} $. Then $DifH(M)$, $DifH_N(M)$ and
$DifH_{N,p}(M)$ are holomorphic Lie groups with Lie algebras $Vect
(M)$, $Vect _N(M)$ and $Vect _{N,p}(M)$ consisting of the
holomorphic vector fields on $M$, those tangent to $N$ and zero at
points of $N$ respectively such that $DifH_{N,p}(M)\subset
DifH_N(M)\subset DifH(M)$.}
\par At first proof two lemmas.
\par {\bf 16. Lemma.} {\it Let $M$ and $N$ be as in Theorem 15, then
there exists an $\bf K$ Hermitian metric $g$ on $M$ such that $N$ is
totally geodesic. There is $\epsilon >0$ such that if $v\in T_xM$
and $\| v \| <\epsilon $, where $x\in N$, then $v\in T_xN$ if and
only if $\exp_xv\in N$.}
\par {\bf Proof.} Let at first $N$ be the zero section of a vector
bundle $\pi : E\to N$. In view of Proposition 11 there exists a
product $\bf K$ holomorphic Hermitian metric $G$ in $E$ such that
$G$ is constant on fibers, where fibers are perpendicular to $N$.
Then $G$ makes $N$ totally geodesic in $E$, since the vector field
providing the geodesic spray is tangent to $TN$. This means that the
Christoffel symbols $\Gamma ^i_{j,k}$ are zero if $i=b+1,...,n$,
$j=1,...,b$, where $\mbox{ }_1z,...,\mbox{ }_mz$ are coordinates for
$N$ and $\mbox{ }_1z,...,\mbox{ }_nz$ those for $E$ (see also \S
4.10 \cite{luladfcdv}). In general let $V$ be a tubular neighborhood
chart for $N$ and consruct $G$ as above.
\par Since $N$ is foliated in $M$, then $TN$ is foliated in $TM$.
Therefore, there exists a $\bf K$ Hermitian holomorphic metric $u$
in $M$ which has as restriction a $\bf K$ holomorphic metric $g$ in
$N$. This $g$ makes $N$ totally geodesic as it equals $G$ on a
neighborhood of $N$.
\par The last statement follows from the first, since there exists
$\epsilon >0$ such that if $\| v \| <\epsilon $, $v\in T_xM$, then
$x$ and $\exp_xv$ are joined by a unique geodesic in a given
neighborhood.
\par {\bf 17. Lemma.} {\it There is a chart $F: U\to X$ such that
$F(U)=V_1\times V_2$, where $V_1$ and $V_2$ are open subsets of
separable Hilbert spaces over $\bf K$, $V_1\subset X$ and such that
$F(U\cap DifH_N(M))=V_1$, where $X:= Vect _N(M)$. A similar
statement is true for $DifH_{N,p}(M)$ with $X_p := Vect _{N,p}(M)$.}
\par {\bf Proof.} Mention that $Vect _N(M)= \{ v\in Vect (M):
v(x)\in T_xN \forall x\in N \}$, $Vect _{N,p} = \{ v\in Vect (M):
v(x)=0 \forall x\in N \}$. Put $V_1=F(U)\cap X$, where $F$ is the
exponential chart on $DifH(M)$ for the metric of Lemma 16. Assume
that $v\in F(U)$ implies $\| v \| <\epsilon $. Let $V_2$ be a
neighborhood in the orthogonal complement of $X$ in $Vect (M)$. Then
taking restriction in a case of necessity put $F(U)=V_1\times V_2$.
In view of Lemma 16 $F(U\cap DifH_N(M))=V_1$. For $DifH_{N,p}(M)
\subset DifH(M)$ it can be used any $\bf K$ Hermitian holomorphic
metric, since if $v\in T_xM$, $v=0$ if and only if $\exp_xv=x$ for
$\| v \| < \epsilon $. Then repeat the proof for $X_p$ instead of
$X$. Using the above chart for $DifH_N(M)$ also shows, that
$DifH_{N,p}(M)$ is the submanifold of $DifH_N(M)$ at $id$.
\par {\bf Proof} of Theorem 15 with the help of Lemmas 16 and 17
finishes by the observation, that charts at $f\in DifH(M)$ are
obtained by right translation and since the uniform spaces being
submanifolds near $id \in DifH(M)$, hence they are submanifolds.
\par {\bf 18. Theorem.} {\it Let $M$ be a compact holomorphic manifold
on $\bf K^n$ with boundary $\partial M$ and $DifH_p(M) := \{ f\in
DifH(M): f(x)=x \forall x\in \partial M \} $. Then $DifH(M)$ and
$DifH_p(M)$ are holomorphic Lie groups with Lie algebras $Vect
_{\partial M}(M)$, $Vect _{\partial M,p}(M)$ the holomorphic vector
fields tangent to $\partial M$ and those zero on $\partial M$
respectively such that $DifH_p(M)\subset DifH(M)$. Moreover,
$DifH(M)$ and $DifH_p(M)$ can be supplied with uniformities relative
to which they are complete.}
\par {\bf Proof.} If $f\in DifH(M)$, then $f$ is holomorphic, hence
it is a $C^{\infty }$ mapping, consequently, $f(\partial M)=\partial
M$ and $f(M\setminus \partial M) = M\setminus \partial M$. To finish
the proof we need the following lemma.
\par {\bf 19. Lemma.} {\it Let $M$ be embedded into its double
covering holomorphic manifold $\hat M$. Then $DifH(M)\subset
H(M,{\hat M})$ is a submanifold, as well as $DifH_p\subset DifH(M)$,
where $H(M,{\hat M})$ denotes the set of all holomorphic mappings
from $M$ into $\hat M$.}
\par {\bf Proof.} In accordance with Proposition 11 and Lemma 16
there exists a $\bf K$ Hermitian holomorphic metric on $\hat M$ such
that $\partial M$ becomes totally geodesic. Suppose $f\in
DifH(M)\subset H(M,{\hat M})$, and choose an exponential chart $\exp
: U\to H(M,{\hat M})$ in a neighborhood $U\subset T_fH(M,{\hat M})$
of $(f,0)$. Consider $X_f(M):= \{ V\in H(M,T{\hat M}): V$
$\mbox{covers}$ $f$ $\mbox{and}$ $V(x)\in T_x\partial M$  $\forall
x\in \partial M \} $, then $X_f(M)$ is the closed uniform subspace
in the tangent bundle $H(M,T{\hat M})=TH(M,{\hat M})$. The proof of
Lemma 17 shows, that $\exp $ maps $U\cap X_f(M)$ onto a neighborhood
of $f$ in $DifH(M)$, that defines the submanifold chart, since
$\partial M$ is totally geodesic. The case $DifH_p(M)$ can be proven
similarly.
\par {\bf 20. Corollary.} {\it If $f\in DifH(M)$ is in the component
of the identity, then $f$ has an extension $y\in DifH({\hat M})$.}
\par {\bf Continuation of the proof of Theorem 18.}
Topologize the family $H(M,P)$ of all holomorphic mappings from a
holomorphic manifold $M$ on $\bf K^n$ into a holomorphic manifold
$P$ on $\bf K^p$ as the family of locally analytic functions $f$,
satisfying the condition $\partial f_{i,j}/\partial \mbox{ }_kz=0$
for each $\bf K$ variable $\mbox{ }_kz$, where $f_{i,j} = \phi
_i\circ f\circ \psi _j^{-1}$, $At (M) = \{ (U_i,\phi _i): i \} $,
$At (P) = \{ (W_j,\psi _j): j \} $. Using charts and components of
functions it is sufficient to describe this topology of $H(U,{\bf
K^p})$ for a domain $U$ in $\bf K^n$, where $U$ is not necessarily
open. Namely, consider subspaces $X(z_0,R)$ of all $f\in H(U,{\bf
K^p})$ with a finite norm
\par $\| f \|_{z_0,R} := \sup_{(c=0,1,2,...; z\in B({\bf K},\mbox{
}_1z_0,R_1)\times ... \times B({\bf K},\mbox{ }_nz_0,R_n))} |\phi
_{f,c}(z-z_0)|<\infty $, \\
where $\phi _{f,c}$ is given by Formulas 2.13(ii,iii),
$R=(R_1,...,R_n)$, $0<R_1,...,R_n<\infty $. The spaces $X(z_0,R)$
and $H(U,{\bf K^p})$ are $\bf K$ vector spaces, hence they have real
shadows, which are $\bf R$ linear spaces. Then the finest locally
convex topology on $H(U,{\bf K^p})$ relative to which all natural
embeddings $\theta _{z_0,R}: X(z_0,R) \hookrightarrow H(U,{\bf
K^p})$ are continuous is called the inductive limit topology (see
Example 5.10.5 \cite{nari}). It is the desired topology on the space
of holomorphic mappings. If $U$ is compact, then there exists a
countable dense subset $\{ z_{0,k}: k\in {\bf N} \} $ in $U$ and
take $R(k)= (1/k,...,1/k)$. Then there are $\bf K$ vector subspaces
$Y_m$ in $H(U,{\bf K^p})$ such that $\bigcup_{m=1}^{\infty }Y_m =
H(U,{\bf K^p})$ and $Y_m\subset Y_{m+1}$ for each $m\in \bf N$. Thus
$H(U,{\bf K^p})$ is supplied with the strict inductive limit
topology. For example, take integers $0<k_1<k_2<...$ and $Y_m :=
\bigcup_{v=1}^m \bigcap_{s=1}^{k_v} X(z_{0,s},R(v))$, where $k_v$,
$R(s)$ and $z_{0,s}$ are chosen such that $\bigcup_{s=1}^{k_v}
B({\bf K},\mbox{ }_1z_{0,s},R_1(v))\times ... \times B({\bf
K},\mbox{ }_nz_{0,s},R_n(v))\supset U$ for each $v$. Each $Y_m$ is
complete, since $H(B({\bf K},\mbox{ }_1z_0,R_1)\times ... \times
B({\bf K},\mbox{ }_nz_0,R_n),{\bf K^p})\cap X(z_0,R)$ is the Banach
space, hence each $\bigcap_{s=1}^{k_v} X(z_{0,s},R(v))$ is complete,
and inevitably $H(M,{\bf K^p})$ is complete for compact $M$ by
Theorem (12.1.6) \cite{nari}. This topology on $H(M,P)$ induces the
topology on $DifH(M)$.
\par Then $DifH(M)$ is the topological group, since the composition
$(f,g)\mapsto f\circ g$ and
the inversion $f\mapsto f^{-1}$ are continuous relative to such
topology. This topology defines the corresponding uniformity on
$DifH(M)$ generated by left shifts $fV$ of neighborhoods $V$ of
$id$, where $V\in {\cal B}(id)$, $f\in DifH(M)$, ${\cal B}(id)$ is
the base of topology at $id$ (see \S 8.1 \cite{eng}). Since $M$ is
compact and $H(M,M)$ is complete, then $DifH(M)$ and $DifH_p(M)$
also are complete. This finishes the proof of Theorem 18.
\par It is necessary to note, that relative to the metric $\rho $
the group $DifH(M)$ has the completion topologically and
algebraically isomorphic with the group $Diff^1(M_{\bf R})$ of all
$C^1$-diffeomorphisms of the real shadow $M_{\bf R}$ of $M$. Thus
$DifH(M)$ and $Diff^k(M_{\bf R})$ are different groups for each
$k\in \bf N$.
\par {\bf 21. Theorem.} {\it The groups $DifH(M)$, $DifH_N(M)$,
$DifH_{N,p}(M)$ of Theorem 15 and $DifH(M)$ and $DifH_p(M)$ of
Theorem 18 admit holomorphic exponential mappings. That is, if $v\in
Vect(M)$ or $v\in Vect _N(M)$, or $v\in Vect _{N,p}(M)$, or $v\in
Vect _{\partial M}(M)$, or $v\in Vect _{\partial M,p}(M)$, then the
flow $\eta _t$ of $v$ is a one parameter locally analytic by $t\in
\bf R$ subgroup of the corresponding group of holomorphic
diffeomorphisms.}
\par {\bf Proof.} Theorem 21 for $DifH(M)$, $DifH_N(M)$,
$DifH_{N,p}(M)$ is evident from Theorem 15, since the flow lies in
$DifH(M)$ and hence in $DifH_N(M)$, $DifH_{N,p}(M)$, when $v\in Vect
_N(M)$, or $v\in Vect _{N,p}(M)$. For groups $DifH(M)$ and
$DifH_p(M)$, when $\partial M\ne \emptyset $, take the double
covering $\hat M$ of the manifold $M$, if $M$ is not orientable.
This $\hat M$ can be chosen holomorphic, when $M$ is holomorphic.
Then a holomorphic vector field $v\in Vect _{\partial M}(M)$ or
$v\in Vect _{\partial M,p}(M)$ on $M$ has a holomorphic extension on
$\hat M$ such that $v$ is tangent to $\partial M$ or zero on
$\partial M$ respectively. Thus the flow $\eta _t: M\to \hat M$ is
holomorphic by $x\in M$ and locally analytic by $t\in \bf R$ and
leaves $M$ invariant or $\partial M$ fixed correspondingly.
Therefore, $M\setminus \partial M$ is also invariant. The
restriction of $\eta _t$ to $M$ finishes the proof of Theorem 21.
\par {\bf 22. Remark.} Since $Vect (M)$ are infinite dimensional
manifolds, then the groups $DifH(M)$ are infinite dimensional
holomorphic Lie groups in accordance with Theorems 15, 18, 21. Apart
from them groups $DifP(M)$ are finite dimensional for compact $M$,
which is proved below.
\par {\bf 23. Theorem.} {\it The family of all pseudoconformal
diffeomorphisms $DifP(M)$ of a compact manifold $M$ over $\bf K$
form the topological metrizable group, which is complete relative to
its metric and locally compact.}
\par {\bf Proof.} In view of Theorem 2.6
the family $DifP(M)$ form the group relative to the composition of
functions, since compositions of functions are asssociative. The
metric $\rho $ on $P(M)$ induces the left invariant metric $\phi
(f,g):=\rho (g^{-1}f,id)$ on $DifP(M)$ (see also Definitions 2.8),
where $P(M,N)$ denotes the family of all pseudoconformal functions
from the pseudoconformal manifolds into a pseudoconformal manifold
$N$ over $\bf K$, $P(M):=P(M,M)$, $id(z):=z$ for each $z\in M$. The
tangent bundle $TP(M)$ of $P(M)$ is isomorphic with $P(M,TM)$, where
$\pi \circ v=f\in P(M)$ for $v\in P(M,TM)$, $\pi M\to TM$ is the
natural projection, $v$ is the vector field on the manifold $M$
along $f$. Since $M$ is the pseudoconformal manifold, then $TM$ is
the holomorphic manifold, that is, connecting mappings of its charts
are holomorphic. Therefore, $P(M,TM)$ is the holomorphic bundle.
Consider the closed subbundle $P(M|TM)$ consisting of all $v\in
P(M|TM)$ such that $\pi \circ v=id$. In view of Theorem 4 there
exists the exponential mapping $\exp : U_{TM}\to M$ of the manifold
$M$, which induces the exponential mapping $Exp : U_{TP(M)}\to PM$,
where $U_{TM}$ is a neighborhood of $M$ in $TM$ and $U_{TP(M)}$ is a
neighborhood of $P(M)$ in $TP(M)$. The Cauchy problem for
differential equations in the class of holomorphic functions has a
unique holomorphic solution by Theorem 4. To each $v\in P(M|TM)$ in
a neighborhood of the zero section there corresponds the one
parameter subgroup $g^t$ in $DifP(M)$ such that $(dg^t/dt)|_{t=0}
=v$, where $t\in \bf R$ in the nonperiodic case or $t\in
[-c,c]\subset \bf R$ for some period $c>0$ such that $g^c=id$ in the
periodic case, where $g^0:=id$. Thus $T_{id}DifP(M)$ is the
neighborhood of zero section in $P(M|TM)$ (see also Theorem 21
above). The completeness of $B$ and sequential completeness of the
ball $B := \{ f\in P(M): \rho (f,id)\le r \} $ for $0<r< \max
(1,diam (M))/2$ follows from Lemma 2.11, since $id'=I$ is the unit
operator, since the limit function $f$ of any sequence $f_n\in B$
can not be with degenerate derivative $f'$, where $diam
(M):=\max_{x, y\in M} d(x,y)$, $d$ is the metric in $M$,
consequently, $T_{id}P(M)$ is complete. Since $P(M|TM)$ is closed in
$TP(M)$, then $P(M|TM)$ is complete.
\par Recall necessary topological definitions and theorems.
A topological space $X$ is called a $T_1$-space if for each two
different points $x_1\ne x_2\in X$ there exists an open subset
$U\subset X$ such that $x_1\in U$ and $x_2\notin U$. A topological
space $X$ is called a Hausdorff or $T_2$-space if  for each two
different points $x_1\ne x_2\in X$ there exist two open subsets
$U_1$ and $U_2$ in $X$ such that $x_1\in U_1$ and $x_2\in U_2$ and
$U_1\cap U_2=\emptyset $. A topological space $X$ is called regular
if $X$ is a $T_1$-space and for each point $x\in X$ and each closed
subset $F\subset X$ such that $x\notin F$ there exist open subsets
$U_1$ and $U_2$ in $X$ such that $x\in U_1$ and $F\subset U_2$ and
$U_1\cap U_2=\emptyset $. \par A topological space $X$ is called
Lindel\"of if it is regular and from each its open covering can be
chosen a countable subcovering. A metric space $(X,\rho )$
consisting of a topological space $X$ and a metric $\rho $ in it is
called complete if each Cauchy sequence in $X$ is converging in $X$.
A subset $A$ in a topological space $X$ is called dense in $X$ if
the closure $cl (A)$ of $A$ is equal to $X$. A density $d(X)$ of a
topological space $X$ is defined as the least cardinal number $card
(A)$ among all subsets $A$ dense in $X$. If $d(X)\le \aleph _0 :=
card ({\bf N})$, then $X$ is called separable. \par A topological
space $X$ is called compact if from each its open covering can be
chosen a finite subcovering. Then $X$ is called compactum if $X$ is
compact and Hausdorff. A topological space $X$ is called locally
compact if for each point $x\in X$ there exists a neighborhood $U$
of $x$ such that $cl (U)$ is the compact subspace in $X$.
\par A topological space $X$ is called countable compact if from
each its countable open covering can be chosen a finite subcovering.
A topological space $X$ is called sequentially compact if each
sequence of points in $X$ has a converging subsequence.
\par In accordance with Theorem 3.10 \cite{eng} a topological space
$X$ is compactum if and only if it is Lindel\"of and countable
compact.
\par Theorem 3.10.30 \cite{eng} states that each sequentially
compact space is countably compact.
\par By Corollary 3.3.11 \cite{eng} a Hausdorff space $X$
is locally compact if and only if it is homeomorphic to an open
subspace of some compactum.
\par On the other hand, the metric space $(P(M),\rho )$
is separable, consequently, it is Lindel\"of together with
$DifP(M)$. Then due to Theorems 3.10.1, 3.10.30 and Corollary 3.3.11
\cite{eng} $DifP(M)$ is locally compact.
\par {\bf 24. Theorem.} {\it  A group $DifP(M)$ of pseudoconformal
diffeomorphisms has the structure of an analytic manifold with a
number of real parameters no more than $n(7n+4)$ for ${\bf K}=\bf H$
or $n(29n+8)$ for ${\bf K}=\bf O$, where $n$ is the dimension of $M$
over $\bf K$.}
\par {\bf Proof.} Take a domain $U$
in $\bf K^n$, since with the help of charts of an atals $At (M)$ a
general case can be reduced to it. Consider local coordinates
$\mbox{ }_1z,...,\mbox{ }_nz$ in $U$, where $\mbox{ }_pz\in \bf K$
for each $p=1,...,n$. If $f\in DifP(U)$, then this diffeomorphism
gives the transformation of local coordinatesзадает преобразование
локальных координат $\zeta =f(z)$. To a point $O\in U$ with
coordinates $z$ there is posed a point $f(z)=:S$ with coordinates
$\zeta $, that gives the shift of the space $\bf K^n$ on the vector
$OS$, characterized by $n$ quaternion or octonion respectively
coordinates or by $2^rn$ real coordinates. At the same time there
are given $n^2$ partial derivatives $\partial \mbox{
}_pf(z)/\partial \mbox{ }_jz$. Then relative to the metric $\rho $
in $DifP(U)$ each diffeomorphism $f$ is characterized by $n$
quaternion or octonion correspondingly parameters of a shift and
values of its $n^2$ partial derivatives for each $z\in U$. In
accordance with Theorems 2.17 and 2.18 $f(z_0)$ and $f'(z)$ restore
a function $f$, since the integral along a rectifiable path of a
holomorphic function or of its derivative in a domain $U$ satisfying
conditions of Remark 2.12 depends only on an initial and final
points and not on a path by Theorem 2.11 \cite{ludoyst,ludfov}. From
Theorems 2.4 and 2.5 it follows, that each partial derivative
$\partial \mbox{ }_pf(z)/\partial \mbox{ }_jz$ belongs to the group
$G := S^3\otimes S^3\otimes (0,\infty )$ for quaternions or $G :=
(S^7\otimes (S^7\otimes S^7))\otimes S^7\otimes (0,\infty )$ for
octonions, where $S^7:=\{ z\in {\bf O}: |z|=1 \} $ is the
multiplicative, but nonassociative group, hence the group $G$ is
characterized by either $7$ for quaternions or $29$ for octonions
real parameters.
\par Thus the group $DifP(M)$ has no more, than $4n+7n^2=n(7n+4)$
for quaternions or $8n+29n^2=n(29n+8)$ for octonions real
parameters.
\par The manifold $M$ is holomorphic, since the connecting mappings
$\phi_i\circ \phi _j^{-1}$ are pseudoconformal, hence holomorphic.
It was proved above, that the exponential mapping $\exp_z(v)$ is
holomorphic by $z\in M$ and $v\in Vect (M)$, $\exp : {\tilde T}M\to
M$ is given on the corresponding neighborhood ${\tilde T}M$ of a
manifold $M$ in the tangent bundle, $\exp_z:=\exp |_{{\tilde T}M\cap
T_zM}$, $Vect (M)$ is the set of holomorphic vector fields on $M$
(see Theorems 15, 18 and 21). \par Recall, that if $\pi : Q\to J$ is
a vector bundle and $Y: M\to J$ is a morphism, then a section along
$Y$ is a morphism $\eta : M\to Q$ with $\pi \circ \eta =Y$, where
$M$ and $J$ are manifolds and for each $x\in J$ the fibre $\pi
^{-1}(x)$ is a Banach space (for example, $\bf K^n$ or $l_2({\bf
K})$ in this paper). If $f: M\to N$ is a differentiable mapping,
then there is defined its tangential $Tf: TM\to TN$. Consider the
right multiplication: $R_g: f\mapsto f\circ g$, $R_g: DifH(M)\to
DifH(M)$ for $g\in DifH(M)$, $R_g: DifP(M)\to DifP(M)$ for $g\in
DifP(M)$. Therefore, $TR_g: V\mapsto V\circ g$, where $V$ is a
vector field, $V_s\in T_sDifH(M)$ or $V_s\in T_sDifP(M)$
correspondingly for each $s\in DifH(M)$ or $s\in DifP(M)$
respectively. The right multiplication is holomorphic. Then it is
possible to call a vector field right invariant if
$TR_g(V_s)=V_{s\circ g}$ for each $s, g$ as above. Therefore, the
space of right invariant vector fields on $DifH(M)$ is isomorphic
with $T_{id}DifH(M)$, which is seen by evaluation at $id$. This
means that for each $V\in T_{id}DifH(M)$ there exists $v\in Vect
(M)$ such that $V_f=v\circ f$ for each $f\in DifH(M)$.
\par Thus, the tangent bundle $TP(M)$ is isomorphic with the family
$H_P(M,TM)$ of all $\bf K$ holomorphic sections $F$ from $M$ into
$TM$ with pseudoconformal mappings $\pi \circ F: M\to M$, where $\pi
: TM\to M$ is the natural projection. While for each $V\in Vect
(DifP(M))$ there exists $v\in Vect (M)$ such that $V_f=v\circ f$ for
each $f\in DifP(M)$, where $V_f\in T_f DifP(M)$, since as in the
real case also it is sufficient to consider right invariant vector
fields (see the real case in \cite{ebmars,lurim1}). Then on a
neighborhood ${\tilde T}DifP(M)$ of the manifold $DifP(M)$ in the
tangent bundle $T DifP(M)$ there exists the induced mapping
$Exp_f(V) := \exp_{f(z)}(v\circ f)$, which is holomorphic, where
$f\in DifP(M)$.
 Therefore, the mapping $Exp$ induces the real
locally analytic atlas of the manifold $DifP(M)$.
\par {\bf 25. Theorem.} {\it  The group $DifP(M)$ is the analytic
Lie group over $\bf R$.}
\par {\bf Proof.} It remains to verify that the mapping
$DifP(M)^2\ni (f,g)\mapsto g^{-1}f\in DifP(M)$ is locally analytic
relative to the real parameters $t_1,...,t_p$, providing locall
coordinates in $DifP(M)$ as the real analytic manifold, $p:=dim_{\bf
R}DifP(M)$. Let $q(t): U\to U$ be a $\bf K$ holomorphic function on
a domain $U$ in $\bf K^n$, depending analytically on the real
parameter $t$ in a neighborhood $B$ of a point $t_0\in \bf R$ and
consider the differential equation $dq(t)/dt=\Psi (q,t)$ with the
initial condition $q(t_0)=q_0$, then it has a unique locally
analytic solution (see Theorem 4) $q(t)$ in an interval
$(-a+t_0,a+t_0)$ with $a>0$, if $\Psi (q,t)$ is locally analytic by
$z\in U$ and $t\in B$. Vice versa, if $q(t)$ is locally analytic,
then $\Psi $ is also locally analytic. Then one parameter subgroups
in $DifP(M)$ are locally analytic. Since $g\circ f(z):=g(f(z))$,
while $f$ and $g$ are locally analytic by $\bf K$ variables, then
the mapping $(g,f)\mapsto gf$ is locally analytic relative to
$t_1,...,t_p$. On the other hand, $f^{-1}\circ f=id$, consequently,
the mapping $f\mapsto f^{-1}$ is also locally analytic. Thus, the
mapping $\phi =\phi (u_1,...,u_p;v_1,...,v_p)$ is locally analytic,
where $(u_1,...,u_p)$, $(v_1,...,v_p)$, $\phi =(\phi _1,...,\phi
_p)$ are local coordinates corresponding to $f, g$ and $g^{-1}f$
respectively.
\par {\bf 26. Definitions.} A domain $U$ in $\bf K^n$ we call
a perfect domain of $\cal P$ transformations, if $U\subset cl (Int
(U))$ and for each $a, b\in Int (U)$ there exists a pseudoconformal
mapping $f: U\to U$ such that $f(a)=b$ and $f(b)=a$, where $n\in \bf
N$.
\par A domain $U$ in $\bf K^n$ is called simply connected, if
$\partial U$ is connected.
\par {\bf 27. Examples.} The ball $B({\bf K^2},0,1)$ and
the polydisc $B({\bf K},0,1)^2$ are perfect domains of $\cal P$
transformations.
\par {\bf 28. Proposition.} {\it If $U$ is a bounded open perfect
domain of $\cal P$ transformations, then $U$ is a domain of
existence of a normal family $\cal F$ and $U$ is pseudoconvex.}
\par {\bf Proof.} Consider the family $Q(U,U)$ of all
mappings $f: U\to U$ such that $f$ for each $\xi \in U$ by each
variable $\mbox{ }_vz$ is either pseudoconformal or $\partial f(\xi
)/\partial \mbox{ }_v\xi =0$. Then $P(U,U)\subset Q(U,U)$, where
$P(U,U)$ is nonvoid due to Definition 26. Choose in $Q(U,U)$ a
bounded relative to the metric $\rho $ (see \S 2.8) subfamily $\cal
F$. A domain of definition of this family $\cal F$ is $V$ such that
$U\subset V$. Consider the subfamily $S$ of $f\in \cal F$ such that
for each point $x$ in $U$ there exists its open neighborhood $W$ and
$f\in S$ such that $f(W)$ is open in $\bf K^n$. Suppose that $U\ne
V$ and $cl (U)\subset Int (V)$, then for each $P_0\in \partial U$
there exists a neighborhood $G$ of $z$ such that $G\subset V$, that
is, $P_0\in Int (V)$. There exists a sequence $P_n\in U$ such that
$\lim_{n\to \infty }P_n=P_0$. \par For each $n$ choose $f_n\in S$
with $f_n(P_n)=z_0$, where $z_0$ is a marked point in $U$. Without
loss of generality suppose that $z_0=0$, since it is possible to use
shift on a vector to new system of coordinates. Extract from $S$ a
converging subsequence in $\cal F$ such that $\lim_{n\to \infty
}f_n=f_0$ realtive to $\rho $ on $V$. Since $P_0\in Int (V)$, then
for each $\epsilon >0$ there exist $N\in \bf N$ and $\delta >0$ such
that $|f_n(z)-f_n(P_0)|<\epsilon $ for each $z$ with $|z-P_0|<\delta
$ and $n>N$. But $f_n(P_n)=0$, hence for each $\epsilon >0$ there
exists $N\in \bf N$ such that $|f_n(P_n)-f_n(P_0)| =
|f_n(P_0)|<\epsilon $ for each $n>N$, consequently, $\lim_{n\to
\infty }f_n(P_0)=0$ and inevitably $f_0(P_0)=0$. On the other hand,
$P_0\in \partial U$ and $f_0$ is holomorphic, hence $f_0(P_0)\notin
U$. This leads to the contradiction, hence $U=V$.
\par The last statement follows from Theorems 4.7 and 4.8
\cite{luladfcdv}.
\par {\bf 29. Proposition.} {\it If $U$ is a perfect domain of $\cal P$
transformations, then $U$ is simply connected.}
\par {\bf Proof.} This follows from Proposition 28 and Theorem
2.56 and Corollary 2.58.
\par {\bf 30. Theorem.} {\it Let $U$ be a bounded connected domain in
$\bf K^n$ containing point $0$. Suppose that $f: U\to U$ is a
pseudoconformal function onto $U$ such that $f(0)=0$ and $f'(0)=I$,
where $I$ is the unit operator. Then $f=id$ is the identity
functions, $id(x)=x$ for each $x\in U$.}
\par {\bf Proof.} Since $f$ is holomorphic, then it has
a power series expansion in a ball with centre $0$ of positive
radius $0<R<\infty $ given by Formulas 2.13(ii,iii). Therefore, in
accordance with Theorems 2.4 and 2.5 $f(z)$ can be written in the
form: $f(z)=z+P(z) + o(|z|^{\alpha })$, where $P(z)$ is a polynomial
$P: U\to \bf K$ of degree $\alpha \ge 2$ with $P'(0)=0$, $\lim_{z\to
0} o(|z|^{\alpha })/|z|^{\alpha }=0$. Then $f\circ f(z) = z + 2P(z)
+ o(|z|^{\alpha })$, hence by induction $f^{\circ k}(z) = z + k P(z)
+ o(|z|^{\alpha })$, where $f^{\circ k}(z) := f\circ f^{\circ
(k-1)}(z)$ for each $k>1$, $f^{\circ 2}(z) := f\circ f(z)$,
$f^{\circ 1}(z) := f(z)$. Therefore, \par $(i)$ $\{ (a_{f^{\circ
k},m}, (z-z_0)^m) \} _{q(2v)}= k
\{ (a_{f,m}, (z-z_0)^m) \} _{q(2v)}$, \\
for each $|m|\ge 2$, where $a_{f,m,j}\in \bf K$, $m =
(m_1,...,m_v)$, $m_j=0,1,2,...$, $a_{f,m}=(a_{f,m,1},...,
a_{f,m,v})$, $|m| := m_1+...+m_v$, $v\in \bf N$, $\{ (a_{f,m},z^m)
\} _{q(2v)}:= \{ a_{f,m,1}z^{m_1}...a_{f,m,v}z^{m_v} \} _{q(2v)}$.
There exist $C$ and $R$ and $R_1$ such that $0<C<R<R_1<\infty $ and
$P := B({\bf K},0,C)^n\subset B({\bf K^n},0,R)\subset U\subset
B({\bf K^n},0,R_1)$. Since $f(U)\subset U$ and $U$ is bounded, then
$|f^{\circ k}(z)|<R_1$ for each $z\in P$. In view of Theorem 3.10
and Remarks 3.34 \cite{ludfov} and Theorem 2.28 and Equality $(i)$
there exists a constant $0<Y<\infty $ such that \par $(ii)$
$\sup_{z\in P}|\phi _{f^{\circ
k},s}(z)| \le k \sup_{z\in P}|\phi _{f,s}(z)| \le Y R/C^s$ \\
for each $s=2,3,...$, where
\par $\phi _{f,s}(z) = \sum_{|m|=s} \{ (a_{f,m}, z^m) \} _{q(2v)}$. \\
If $\sup_{z\in P}|\phi _{f,s}(z)|\ne 0$ for some $s\ge 2$, then
$(ii)$ leads to the contradiction, since $k\in \bf N$ can be
arbitrary large. Therefore, $f(z)=id(z)=z$ on $P$. Since $U$ is
bounded and connected, then $cl (U)$ is the connected compact in
$\bf K^n$. The using of shifts in $\bf K^n$ on vectors and the
continuing of this proof gives a finite covering of $U$ by polydics
$P_j$ in each of which $f=id$ and inevitably $f=id$ on $U$.
\par {\bf 31. Remark.} If $X_l$, $X_r$ and $X$ are vector spaces
over $\bf K$ relative to multiplications on constants from the left,
right or on both sides, then generally $X_l$, $X_r$ and $X$ may be
pairwise distinct, for example, when their basic vectors are
$x_n:=z^n$, where $z\in U$, $U$ is a domain in $\bf K$. But $X_l$ is
isomorphic with $X_r$ in the way $av\mapsto va$ for each $a\in \bf
K$ and $v$ in the real shadow over $\bf R$ of $X_l$, hence $X$ is
also isomorphic to these spaces $X_l$ and $X_r$.
\par {\bf 32. Example.} Let $M\in \bf K$, $Re (M)=0$, $|M|=1$.
Consider the domain $V:= \{ Z\in {\bf K^{n^2}}: Re (M^*Z)>0 \} $,
where ${\bf K}=\bf H$ or ${\bf K}=\bf O$, $Z$ is a matrix with
entries $Z_{i,j}\in \bf K$, $Re (Z):= A$ such that $A_{i,j}=Re
(Z_{i,j})$ for each $i, j$. A real matrix $A$ is positive definite,
if it is symmetric $A^T=A$ and $x^TAx>0$ for each vector column
$x\ne 0$ in $\bf R^n$, where $Z^T$ denotes the transposed matrix of
$Z$. A Hermitian conjugate matrix is $Z^* := {\tilde Z}^T$, where
each ${\tilde Z}_{i,j}$ is the conjugate of $Z_{i,j}\in \bf K$. Then
$Re (Z)=(Z+Z^*)/2$. A matrix $Z$ is Hermitian if $Z^*=Z$. A
Hermitian matrix $Z$ is positive definite, if $x^*(Zx)>0$ for each
vector column $x\ne 0$ in $\bf K^n$. If $A$ and $B$ are two
Hermitian matrices, then by the definition $A<B$ if $B-A>0$. Put
$B_n := \{ W\in {\bf K^{n^2}}: WW^*<E \} $, where $E$ is the unit
matrix, $E=diag (1,...,1)$. We call $V$ the generalized upper
$M$-semispace.
\par If $Z\in V$, then rows of the matrix $Z+ME$ are independent
over $\bf K$. Suppose the contrary, that $(Z+ME)w=0$ for some
nonzero vector column $w\in \bf K^n$. Thus $Zw=-Mw$ and
$w^*Z^*=w^*M$, hence $w^*((M^*Z+Z^*M)w)/2=-w^*w$ due to the
alternativity of $\bf K$, but $Re (Z)>0$, consequently, $w=0$. This
contradicts our supposition, hence rows of the matrix $Z+ME$ are
independent over $\bf K$. Therefore, there exists $(Z+ME)^{-1}$ and
inevitably \par $(i)$ $W(Z):=(Z+ME)^{-1}(Z- M E)$ \\
is the holomorphic function. Henceforts, in this section 32 let
$n\ge 1$ for ${\bf K}=\bf H$ and $n=1$ for ${\bf K}=\bf O$.  Then
$(Z+ME)(Z^*-ME) - (Z-ME)(Z^*+ME)= -2ZM +2MZ^* = - 2(ZM) - 2 (ZM)^*=
- 4 Re (ZM)$, hence $E-WW^*= - 4 (Z+ME)^{-1}(Re (ZM))(Z^*-ME)^{-1}$,
where $- Re (ZM) = Re (MZ^*)= Re (M^*Z)$, consequently, $E-WW^*$ is
positive definite  (see also Formulas 2.5(i-iv)).
\par There are identities: $Z(E-W)=Z(E-(Z+ME)^{-1}(Z-ME))=
Z-Z(E-2(Z+ME)^{-1}ME)=2Z((Z+ME)^{-1}ME)$, $(Z+ME)W=Z-ME$,
$ZW=Z-ME-MW$, $Z(E-W)=M(E+W)$, hence $Z=(M(E+W))(E-W)^{-1}$. The
matrix $(E-W)$ is nondegenerate. If it would be $(E-W)x=0$ for some
vector column in $\bf K^n$, then $Ex=Wx$ and $x^*=x^*W^*$,
consequently, $x^*(Wx)=x^*x$, which for $x\ne 0$ is in contradiction
with $WW^*<E$. If $WW^*<E$, then $(Z+ME)^{-1}$ is nondegenerate and
$Re (M^*Z)>0$ due to the formulas above. Hence $W: V\to B_n$ is the
pseudoconformal diffeomorphism. The same mapping transforms
$\partial V = \{ Z: Re (M^*Z)=0 \} $ onto $\partial B_n= \{ W:
WW^*=E \} $. Since $e_j^*(E-WW^*)e_j=1-\sum_k|W_{j,k}|^2>0$ for each
$j$, then $|W|^2:=\sum_{j,k}|W_{j,k}|^2<n$, that is, $|W|<n^{1/2}$
for each $W\in B_n$.
\par {\bf 33. Definition.} A $\bf K$  homogeneous norm in $\bf K^n$
is a mapping $ \| * \| : {\bf K^n}\to [0,\infty )$ such that
\par $(1)$ $ \| z +w \| \le \| z \| + \| w \| $ for each $z$ and
$w\in \bf K^n$,
\par $(2)$ $\| cz \| = |c| \| z \| $ for each $z\in \bf K^n$ and
$c\in \bf K$,
\par $(3)$ $\| z \| =0$ if and only if $z=0$.
\par {\bf 34. Theorem.} {\it Let $U$ be a domain in $\bf K^n$ and
let $f: U\to \bf K^m$ be a pseudoconformal mapping. Soppose also
that $ \| * \| $ is a $\bf K$ homogeneous norm in $\bf K^m$. If $\|
f'(z).h \| $ attains its maximum at some point $\alpha \in U$ for
each $h\in \bf K^n$, then \par $(1)$ components $\mbox{
}_j(f'(z).h)$ of the mapping $f'(z).h$ are linearly dependent over
$\bf K$ in $U$ for each $h\in \bf K^n$,
\par $(2)$ $\| f'(z).h \| =C$ is constant in $U$ by $z$, where
$C=const \in \bf R$ may depend on $h$.}
\par {\bf Proof.} Let $\beta =f'(\alpha ).h$ and $B:=B({\bf
K^m},0,\| \beta \| )$ be the ball in $\bf K^m$ relative to the
considered norm, where $\beta $ may depend on $h$. Consider $h\in
\bf K^n$ with $|h|=1$, since $\| f'(z).(bh) \| =|b| \| f'(z).h \| $
for each $b\in \bf R$. From properties 33(1-3) of the norm it
follows, that the interior $Int (B)$ is open and $\bf K$-convex in
$\bf K^m$. If $\| \beta \| =0$, then $f'=0$ on $U$, so consider the
case $\| \beta \| >0$. Since $\beta \in
\partial B$, then there exists a real hyperplane $\pi $ of support
to $\partial B$ such that $\beta \in \pi $ and $(\pi \cap B)\subset
\partial B$. Write this hyperplane in the form: $Re \xi (w)=\beta
_0$, where $\beta _0:= Re (\beta )$, $\xi (w)=\sum_{k=1}^m\mbox{
}_ka\mbox{ }_kw$, $\mbox{ }_ka, \mbox{ }_kw\in \bf K$, $w=(\mbox{
}_1w,...,\mbox{ }_mw) \in \bf K^m$. Choose the function $\xi $ such
that $Re \xi (w)<\beta _0$ for each $w\in Int (B)$. The function
$\exp (\xi \circ f(z))$ is pseudoconformal as the composition of
pseudoconformal functions. Then $|g(z)|=\exp [Re (\xi \circ
(f'(z).h))]\le \exp (\beta _0)$, moreover, $|g(\alpha )|=\exp (\beta
_0)$. In view of Theorem 2.28 $\xi \circ (f'(z).h)=const $ in $U$ by
$z$, hence the components $\mbox{ }_k(f'(z).h)$ are linearly
dependant over $\bf K$, that is, they satisfy the equation:
$\sum_{k=1}^m \mbox{ }_ka\mbox{ }_k(f'(z).h)=C_1 \in \bf K$ in $U$
by $z$ for each $h\in \bf K^n$, where $C_1=const $ may depend on
$h$. Therefore, $f'(z).h\in \pi \cap
\partial B$, hence $\| f'(z).h \| = \| \beta \| $ for each $z\in U$.
This finishes the proof of Theorem 34.
\par {\bf 35. Remark.} If the ball
$B$ in the considered norm is strictly convex, then $\pi \cap
\partial B$ is the point and $f'(z).h=C$ on $U$ by $z$. For example,
this is the case, when $ \| z \| ^2:= \sum_{k=1}^m|\mbox{ }_kz|^2$.
But in general, for example, $\| z \| := \max_{k=1,...,m} |\mbox{
}_kz|$ relative to the polydics norm the intersection $\pi \cap
\partial B$ is not a point.
\par {\bf 36. Theorem.} {\it Let $B_1:=B({\bf K^n},0,1)$ and
$B_2:=B({\bf K^m},0,1)$ be two balls in $\bf K^n$ and $\bf K^m$
relative to $\bf K$ homogeneous norms $ \| * \| _1$ and $ \| * \|
_2$ respectively. Suppose that $f: Int (B_1)\to Int (B_2)$ is a
pseudoconformal mapping with $f(0)=0$. Then for each $z\in B_1$ the
inequality $\| f(z) \| _2\le \| z \| _1$ is satified for each $z\in
B_1$.}
\par {\bf Proof.} Take the $\bf K$ straight line $\bf l$ given by the
equation: $z=\zeta z_0$, where $\zeta \in \bf K$, $z_0\in \partial
B_1$, that is, $ \| z_0 \|_1=1$. Then ${\bf l}\cap Int (B_1)=Int
B({\bf K},0,1)$. In view of Theorem 34 $ \| (\partial f(z)/\partial
\mbox{ }_jz).h \|$ may achieve its maximum for each $h\in \bf K$ on
$({\bf l}\cap \partial B_1)$ only while varying $\mbox{ }_jz$ for
${\bf l}={\bf K}e_j$, where $e_j=(0,,,0,1,0,...,0)$ with $1$ on the
$j$-th place. Since $f$ is pseudoconformal, then \par $(i)$
$|(\partial \mbox{ }_kf(z)/\partial \mbox{ }_jz).h| = \|
\partial \mbox{ }_kf(z)/\partial \mbox{ }_jz\| \| h \| $ for each
$j=1,...,n$, $k=1,...,m$, $z\in Int (B_1)$ and each $h\in \bf K$. In
view of Theorem 3.20 \cite{rudin} if $F: [a,b]\to \bf R^p$ is a
continuous function on $[a,b]:= \{ x\in {\bf R}: a\le x\le b \} $
differentiable on $(a,b):= \{ x\in {\bf R}: a<x<b \} $, where $p\in
\bf N$, then there exists $x\in (a,b)$ such that $|f(b)-f(a)|\le
(b-a)|f'(x)|$. Consider an $\epsilon $ net of real segments embedded
into $Int (B_1)$ along vectors $e_ji_v$, $v=0,...,2^r-1$, and apply
the aforementioned theorem to $f$. Since $\epsilon >0$ is arbitrary,
$f$ is continuously differentiable and $f(Int (B_1))\subset Int
(B_2)$, then due to Equation $(i)$ the inequality $|f'(z)|\le 1$ on
$Int (B_1)$ is satisfied.
\par Thus $\| f(z)-f(0) \|_2 \le \int_0^1 |f'(tz).z|dt\le
\max_{t\in [0,1]}|f'(tz))| \| z \|_1$ with $z$ in $Int (B_1)$, hence
$\| f(z)-f(z) \|_2 \le (\max_{\xi \in
\partial B_1}|f'(\xi )|) \| z \|_1$ for each $z\in Int (B_1)$.
\par {\bf 37. Theorem.} {\it Let $B_n := \{ z\in {\bf K^n}: |z|\le 1 \} $
be the unit ball in $\bf K^n$,  ${\bf K}=\bf H$ or ${\bf K}=\bf O$,
$<z,a> := \sum_{j=1}^n\mbox{ }_jz\mbox{ }_j{\tilde a}$ is the
Hermitian scalar product in $\bf K^n$, $|z|^2 := <z,z>$. Denote $\pi
_a(z):=<z,a>a/|a|^2$, $\zeta _a(z):=z-\pi _a(z)$,
$S_a(z):=(1-<z,a>)^{-1} (a-\pi _a(z)-(1-|a|^2)^{1/2} \zeta _a(z))$,
where $a\in Int (B_n)\setminus \{ 0 \} $, $S_0(z)=id(z)=z$. Let also
$\bf U_r$ (or $\bf U_l$) be the group of all right (or left) $\bf K$
linear unitary mappings of $\bf K^n$ on $\bf K^n$, that is, $<U(z
),U(\xi )>=<z,\xi >$ for each $z, \xi \in \bf K^n$ and each $U\in
\bf U_r$ (or $U\in \bf U_l$ correspondingly). Then $DifP(Int (B_n))$
consists of all mappings of the form either $g_a(*) := S_a\circ
(U_1(*)\circ U_2)$ with $U_1\in \bf U_r$ and $U_2\in \bf U_l$ over
$\bf H$ or $g_a(*):=S_a\circ ((U_1\circ (U_2(*)\circ U_3))\circ
U_4)$ with $U_1, U_2\in U_r$ and $U_3, U_4\in U_l$ over $\bf O$
respectively}.
\par {\bf Proof.} The mapping $S_a(z)$ is pseudoconformal
on $Int (B_n)$. Since $|<z,a>|<1$ for each $z\in B_n$, then
$(1-<z,a>)^{-1}$ is correctly defined. From $\pi _a(a)=a$ and $\zeta
_a(a)=0$ it follows, that $S_a(a)=0$. For convenience direct the
axis $e_n$ along $a$, that can be achieved by the $\bf R$ linear
pseudoconformal mapping. Then $\pi _a(z)=(0',y)$, where $z=(x,y)$,
$y=\mbox{ }_nz$, $x=z'=(\mbox{ }_1z,...,\mbox{ }_{n-1}z)$ and $\zeta
_a(z)=(x,0)$. Therefore, $S_a(z)=w=(w',\mbox{ }_nw)$ with
$w'=-(1-|a|^2)(1-y\mbox{ }_n{\tilde a})^{-1}x$, $\mbox{ }_nw=
(1-y\mbox{ }_n{\tilde a})^{-1}(\mbox{ }_na-y)$, hence $|w|<1$ for
each $|z|<1$. If $|z|=1$, then $|w|=1$, consequently, $S_a(Int
(B_n))= Int (B_n)$ and $S_a(\partial B_n)=\partial B_n$.
\par On the other hand, $S_a\circ S_a(z)=z$ for each
$z\in Int (B_n)$, since $<z,a>$ and $(1-<z,a>)^{-1}$ commute, hence
the equality $S_a(z)=S_a(\eta )$ implies $z=\eta $. Thus $S_a$ is
the pseudoconformal diffeomorphism.
\par Let now $f\in DifP(B_n)$, take $a\in Int (B_n)$ such that
$f(a)=0$. Then $g:=f\circ S_a^{-1}$ and $g^{-1}$ are pseudoconformal
diffeomorphisms with $g(0)=0$.  In view of Theorem 36 $|g(z)|\le
|z|$ and $|g^{-1}(\eta )|\le |\eta |$ for each $z, \eta \in Int
(B_n)$. The combining these inequalities with $\eta =g(z)$ gives
$|g(z)|=|z|$ for each $z\in Int (B_n)$. \par For $\xi \in \partial
B_n$ and $z\in Int (B_n)$ consider the function $G(v):=g(v\xi
)v^{-1}$ on $Int (B({\bf K},0,1))$ such that $G(0)=0$. Then $G(v)$
is pseudoconformal and $|G(v)|=|g(v\xi )|/|v|=1$. In view of Theorem
34 $|G'(v).h|=C(h)$ for each $h\in \bf K$, where $C\in \bf R$ is
independent from $v$. Thus $G$ is $\bf R$-linear and $\bf K$
additive along each $\bf K$ straight line $\{ z=v \xi : v\in {\bf K}
\} $. In view of Formulas 2.13(i-iii) $g(z)$ is $\bf R$-linear and
$\bf K$ additive. Then due to Theorems 2.4 and 2.5 applied by each
variable $\mbox{ }_jz$ and each component of $\mbox{ }_kg$ as well
as along each $\xi \in Int (B_n)$ the function $g$ has the form
$g(*)=U_1(*)\circ U_2$ with $U_1\in \bf U_r$, $U_2\in \bf U_l$ for
${\bf K}=\bf H$ and $g(*)= (U_1\circ (U_2(*)\circ U_3))\circ U_4$
with $U_1, U_2\in U_r$ and $U_3, U_4\in U_l$ for ${\bf K}=\bf O$.
\par {\bf 38. Theorem.} {\it Let $P_n:= \{ z: \| z \| <1 \} $ be
an open polydisc in $\bf K^n$ relative to the norm $\| z \| :=
\max_{k=1,...,n} |\mbox{ }_kz|$, where ${\bf K}=\bf H$ or ${\bf
K}=\bf O$. Then $DifP(P_n)$ consists of all pseudoconformal
diffeomorphisms \par $(i)$ $\mbox{ }_jz\mapsto (1- \mbox{ }_{\sigma
(j)}\zeta \mbox{ }_{\sigma (j)} {\tilde b})^{-1} [\mbox{ }_{\sigma
(j)}b- \mbox{ }_{\sigma (j)}\zeta ]$, \\
where $\mbox{ }_j\zeta :=(\mbox{ }_jc_3((\mbox{ }_jc_1\mbox{
}_jz)\mbox{ }_jc_2))\mbox{ }_jc_4$, $\mbox{ }_jc_k\in \bf K$,
$|\mbox{ }_jc_k|=1$, $\mbox{ }_jb\in \bf K$ and $|\mbox{ }_jb|<1$
for each $j=1,...,n$, $k=1,...,4$; over $\bf H$ take $\mbox{
}_jc_k=1$ for $k=3,4$ and the order of multiplication is not
essential;  $\sigma : \{ 1,...,n \} \to \{ 1,...,n \} $ is a
transposition.}
\par {\bf Proof.} For each ball $B_1 := B({\bf K},0,1)$ the group
$DifP (Int (B_1))$ was calculated in Theorem 37. Let $f\in DifP
(P_n)$ and $a=f(0)$, where $g$ is chosen in accordance with Formula
$(i)$, equations $\mbox{ }_jb = (\mbox{ }_jc_3((\mbox{ }_jc_1\mbox{
}_ja)\mbox{ }_jc_2))\mbox{ }_jc_4$ define $b=(\mbox{ }_1b,...,\mbox{
}_nb)$. Then the function $F:=g\circ f$ belongs to $DifP(P_n)$. In
view of Theorem 36 applied to $F$ and $F^{-1}$ there is satisfied
the equality $\| F(z) \| = \| z \| $ for each $z\in P_n$. Since $ \|
F \| = \max_{j=1,...,n} \sup_{z\in P_n} |\mbox{ }_jF(z)|$ and
$F(P_n)=P_n$, then there exists a neighborhood $V$ in $P_n$ such
that $\| F(z) \| = |\mbox{ }_jF(z)|$. Let $\| z \| = |\mbox{ }_kz|$
in $V$, hence $|\mbox{ }_jF(z)| = |\mbox{ }_kz|$ in $V$. In view of
Theorem 36 applied to $\mbox{ }_jF$ by the variable $\mbox{ }_kz$
there is the equality $\mbox{ }_jF(z)=\exp (M(z))\mbox{ }_kz$, where
$M(z)\in \bf K$, $|M(z)|=1$, but due to Note 35 $M(z)=:\mbox{}_kM$
is constant. Since $F$ is pseudoconformal and characterized by its
local power series expansions, then this property $\mbox{ }_jF(z) =
\exp (\mbox{ }_kM)\mbox{ }_kz$ is accomplished also in $P_n$. On the
other hand, $F$ is the homeomorphism, hence $k = \sigma (j)$ is the
index transposition.
\par {\bf 39. Theorem.} {\it If $n>1$, then there is not any
pseudoconformal diffeomorphism of $Int (B_n)$ onto $P_n$.}
\par {\bf Proof.} This follows from Theorems 37 and 38, since
$DifP(Int(B_n))\ne DifP(P_n)$. It can also be proved using Theorem
36. If there would be a pseudoconformal diffeomorphism $f: Int
(B_n)\to P_n$, then taking in a case of necessity a composition with
$g\in DifP(P_n)$ it is possible to suppose, that $f(0)=0$. Then $\|
f(z)\| = \| z \| $ for each $z\in Int (B_n)$. Therefore, there would
be a pseudoconformal diffeomorphism of the sphere $ \{ z: |z|=1/2 \}
$ on the nonsmooth surface $\{ z: \| z \| =1/2 \} $, but the latter
is impossible.
\par {\bf 40. Theorem.} {\it Let $M\in \bf K$, $|M|=1$,
$Re (M)=0$, $n\in \bf N$, $Mat_n ({\bf H})$ be the algebra of
$n\times n$ square matrices with elements from ${\bf K}= \bf H$, for
${\bf K}=\bf O$ take $n=1$ and $Mat_1({\bf O})=\bf O$. Let  $V := \{
Z\in Mat_n({\bf K}): Re (M^*Z)>0 \} $ be the domain, where $Re (Z)
:= (Z+Z^*)/2$, $Z^*$ is the Hermitian conjugate quaternion matrix or
octonion, $Re (A)>0$ means the positive definiteness of the matrix
$A$ (see \S 32). Denote by $L_r$ (or $L_l$) the set of all right (or
left respectively) $\bf K$ linear mappings from $\bf K^n$ on $\bf
K^n$ (with alternativity for octonions instead of associativity for
quaternions). Then $DifP(V)$ consists of all bijective epimorphisms
of the form either \par $(1)$ $F (*):= W\circ (P_1(*)\circ P_2)$
over quaternions with $P_1\in L_r$ and $P_2\in L_l$, or \par $(2)$
$F (*):= W\circ ((P_1\circ (P_2(*)\circ P_3))\circ P_4)$ over
octonions with $P_1, P_2\in L_r$, $P_3, P_4\in L_l$, where $W(Z) :=
(Z+D)^{-1}(Z+B)$, $B, D\in Mat_n({\bf K})$, satisfying conditions:
\par $(i)$ $B(M^*D^*)=-(DM)B^*$, $(ii)$ $M^*D^*+MB^*=EM^*$, $(iii)$
$M^*P_1=-P_2^*M$ for quaternions, $M^*(P_1\circ P_2)=-(P_3\circ
P_4)^*M$ for octonions, where $E$ denotes the unit $n\times n$
square matrix.}
\par {\bf Proof.} Consider at first conditions, when the matrix function
$W= (ZC+D)^{-1}(ZA+B)$ preserves the sign $Re (M^*Z)$. We have $2 Re
(M^*Z) = M^*((ZC+D)^{-1}(ZA+B)) + (A^*Z^*+B^*)(C^*Z^*+D^*)^{-1})M
=M^*(((ZC+D)^{-1}[(ZA+B)(M^*(C^*Z^*+D^*)) + ((ZC+D)M)(A^*Z^*+B^*)]
(C^*Z^*+D^*)^{-1})M)$ due to the alternativity. Then
$[(ZA+B)(M^*(C^*Z^*+D^*)) + ((ZC+D)M)(A^*Z^*+B^*)]=
[Z(A(M^*C^*)+(CM)A^*)Z^* +
(ZA)(M^*D^*)+(DM)(A^*Z^*)+B(M^*(C^*Z^*))+((ZC)M)B^* + B(M^*D^*)
+(DM)B^*]$, hence $Re (M^*W)>0$, if \par $(iv)$
$A(M^*C^*)=-(CM)A^*$, $(v)$ $B(M^*D^*)=-(DM)B^*$, $(vi)$
$A(M^*D^*)+(CM)B^*=EM^*$. Under these conditions $Re
(M^*W)=M^*((ZC+D)^{-1}[Re (M^*Z)](C^*Z^*+D^*)^{-1})M$, since $Re
(M^*Z)=Re (ZM^*)=Re (MZ^*)=Re (Z^*M)$, hence $Re (M^*Z)$ and $Re
(M^*W)$ are simultaneously positive or zero. If $Re (M^*W)>0$, then
$W$ is nondegenerate, hence diffeomorphisms of the form given in the
condition of this theorem form the group $G$ relative to
compositions of mappings. \par Consider a subgroup $G_0$ of mappings
$W=A(ZA^*)+BA^*$, then $A(ZA^*)=W-BA^*$ and in view of the
alternativity for each $W\in V$ there exists $Z\in V$ satisfying
this equation, hence $G_0$ and inevitably $G$ acts transitively on
$V$. Since $A(M^*C^*)$ and $B(M^*D^*)$ and $(DM)B^*$ and
$C^{-1}(AM^*)$ are skew Hermitian, since $A(M^*C^*)=-(CM)A^*$, then
$(DC^{-1})(A(M^*D^*))+(DM)B^* = (DC^{-1})M^*$, hence $(DC^{-1})M^*$
is skew Hermitian. The domain $V$ is preserved, when Condition
$(iii)$ is imposed.
\par Let $f\in DifP(V)$, then there exists $g_1\in G_0$ such that
$f_1(ME)=ME$ for $g_1\circ f =: f_1$. Then either
$f_1'(ME).h=(U_1((U_2(B.h))U_3)U_4$ over octonions for each $h\in
\bf O$, where $B$ is a positive definite operator, $U_1,...,U_4\in
S_7\subset \bf O$ (see Theorem 2.5), or $f_1'(ME).h=U_1(B.h)U_2$
over quaternions for each $h\in Mat_n({\bf H})$, where $B$ is a
Hermitian matrix and $U_1, U_2$ are a unitary matrices in
$Mat_n({\bf H})$ (see also Theorem 2.4), since each matrix $n\times
n$ over quaternions is also a $2n\times 2n$ matrix over complex
numbers. If $U_j\ne ME$, then choose $A_j$ in accordance with
Formula 32(i). Then take $g_2$ in accordance with Formulas $(1)$ or
$(2)$ over quaternions or octonions respectively such that either
$g_2'(ME).h=U_1hU_2$ for each $h\in Mat_n({\bf H})$ over
quaternions, or $g_2'(ME).h=(U_1((U_2h)U_3))U_4$ for each $h\in \bf
O$ over octonions. The preserving domain $V$ implies Condition
$(iii)$. Then $f_2:=f_1\circ g_2^{-1}$ has $f_2'(ME)=B$. Take the
pseudoconformal diffeomorphism from $F$ from $V$ onto $B_n$ given by
Formula 32(i). Then $F(ME)=0$ and $F'(ME).h=M^*h/2$ for each $h\in
Mat_n({\bf K})$. Thus $g:=F\circ f_2\circ F^{-1}$ is the
pseudoconformal diffeomorphism of $B_n$ such that $g(0)=0$ and
$g'(0).h=B.h$ for each $h\in Mat_n({\bf K})$. Since $B>0$, then all
eigenvalues $a_1>0,...,a_n>0$ are positive. For $m$-th iteration of
the mapping $g$: $g^m := g\circ g^{m-1}$, where $g^0=id$, $g^1=g$,
$g^2=g\circ g$, eigenvalues of the matrix $(g^m)'(0)$ are
$a_1^m>0,...,a_n^m>0$. For each $c\in \bf K$ with $|c|>1$ there
exists $\lim_{m\to \infty }c^m=\infty $, or for $|c|<1$ there exists
$\lim_{m\to \infty }c^m=0$. If $|c|=1$ and $c\ne 1$, then the
sequence $\{ c^m: m\in {\bf N} \} $ has not a limit. In view of
Theorem 2.13 the family $g^m$ is normal, but it can converge only
when $a_1=1,...,a_n=1$. Thus we can consider, that $g'(0)=E$.
Therefore, by Theorem 30 $g(Z)=Z$ for each $Z$, hence $f_2=id$ and
inevitably $f=g_1^{-1}\circ g_2\in G$.
\par {\bf 41. Definition.} We say that a group $G$ of pseudoconformal
transformations of a domain $U$ in $\bf K^n$ has a property $[P]$,
if for each sequence $f_1,...,f_n,... \in G$ such that $f_j$
converges to $id$ relative to the metric $\rho $ on each its
canonical closed compact subset, where $f_j\ne id $ for each $j$,
there exists a subsequence $f_{j_k}$ with $j_1<j_2<...$ and $j_k\in
\bf N$ such that sequence $m_k(f_{j_k}-id)$ converges relative to
$\rho $ on each canonical closed compact susbet to a function $g$
such that $g$ is not identically zero and for each $p$ and $s$ and
each $z\in U$ either $\partial \mbox{ }_pg/\partial \mbox{ }_sz$ is
zero or $\mbox{ }_pg$ is pseudoconformal by $\mbox{ }_sz$ at $z$.
\par {\bf 42. Theorem.} {\it Let $U$ be a bounded domain in $\bf K^n$
satisfying conditions of Remark 2.12. Then $DifP(U)$ has Property
41$[P]$.}
\par {\bf Proof.} Take a polydisc $B_{z,R}=B({\bf K},\mbox{ }_1z,R_1)
\times ... \times B({\bf K},\mbox{ }_nz,R_n)$ in $Int (U)$, where
$z=(\mbox{ }_1z,...,\mbox{ }_nz)\in U$, $R=(R_1,...,R_n)$,
$0<R_p<\infty $ for each $p$. Choose for each $f_j$ an integer
$m_j$. Take a sequence $m_j(f_j-id)$ such that it is uniformly
equibounded on each canonical closed compact subset $S$ in $U$ and
it has not $id$ as a limit function of some its subsequence. Then
take a locally finite countable covering of $U$ by polydiscs
$B_{z(j),R(j)}$. Put $U_k:=\bigcup_{j=1}^kB_{z(j),R(j)}$. Prove by
induction that it is possible to construct the family with the
demanded property. For this prove, that if a family consists of
functions not identically zero on $D_k$ and uniformly equibounded on
$D_k$, then it satisfies these properties on $D_{k+1}$. For this it
is sufficient to show that it is uniformly equibounded on
$B_{z(k+1),R(k+1)}$. For each transformation $f_j$ it is possible to
choose $t_j\in \bf N$ such that $t_j(f_j-id)$ is uniformly
equibounded on $B_{z(k+1),R(k+1)}$ and does not contain a function
identical to zero. Then $\sup_{j\in \bf N} (m_j/t_j)< \infty $. In
the contrary case the family $t_j(f_j-id)=(t_j/m_j) (m_j(f_j-id))$
is uniformly equibounded on $D_k\cap B_{z(k+1),R(k+1)}$ and it would
have in this intersection a limit function equal identically to
zero, but this is impossible. But then $\sup_{j\in \bf N} (m_j/t_j)<
\infty $, then $m_j(f_j-id)=(m_j/t_j) (t_j (f_j-id))$ is uniformly
equibounded on $B_{z(k+1),R(k+1)}$. In view of Theorem 2.13 then $\{
f_j : j \} $ has a convergent subsequence, while $DifP(U)$ is
complete relative to the compact-open topology of the metric $\rho $
in accordance with Theorem 23, hence $DifP(U)$ has property 41$[P]$.

 \vspace{1cm}

\par {\bf Acknowledgement.} The author is sincerely grateful
to Prof. Jean-Pierre Bourguignon for drawing my attention to groups
of diffeomorphisms of complex manifolds at IHES in November 1997 -
January 1998 and to Prof. Fred van Oystaeyen for helpful discussions
on noncommutative geometry over quaternions and octonions at
Mathematical Department of Antwerpen University in 2002 and 2004 and
both for hospitality.

\end{document}